\documentclass[preprint,3p,times]{elsarticle}

\usepackage{textcomp}
\usepackage{gensymb}
\usepackage{amsmath}
\usepackage{algorithm,algorithmic}
\usepackage{color}
\usepackage{color}
\usepackage{graphicx,subfigure}
\usepackage{array}
\usepackage{tabularx}
\usepackage{longtable}
\usepackage{ulem}
\usepackage{cancel}
\usepackage{url}
\definecolor{airforceblue}{rgb}{0.36, 0.54, 0.66}

\usepackage{array}
\usepackage{color}
\usepackage{tabularx}
\usepackage{graphicx}
\usepackage{amsmath}
\usepackage{amssymb}
\usepackage{amsfonts}	
\usepackage{moreverb}
\usepackage{dsfont}
\usepackage{grffile}
\usepackage{bm}
\usepackage{multirow}
\usepackage{soul}

\usepackage[textsize=footnotesize backgroundcolor=yellow!70, bordercolor=orange]{todonotes}

 \usepackage{amsmath, amstext, amssymb, amsthm, latexsym, color}
 \usepackage{algorithmic, algorithm, epsfig, subfigure, booktabs}

\newcommand{\bea}[0]{\begin{eqnarray}}
\newcommand{\eea}[0]{\end{eqnarray}}
 \newcommand{\xx}{\boldsymbol{x}}

\newcommand{\Int}{\displaystyle\int}
\newcommand{\Sum}{\displaystyle\sum}

\newtheorem{theorem}{Theorem}[section]
\newtheorem{remark}[theorem]{Remark}

\usepackage{multicol}

\usepackage{hyperref} 
\hypersetup{
    colorlinks=true,                          
    linkcolor=blue, 
    citecolor=red, 
    urlcolor=blue  } 
\usepackage[hyperpageref]{backref} 
\RequirePackage[hyperpageref]{backref}
\backreffrench
\renewcommand*{\backref}[1]{}  
\renewcommand*{\backrefalt}[4]{
 \ifcase #1 %
  \relax
  \or
  (Cited page~#2.)%
  \else
  (Cited page~#2.)%
  \fi}

\newcommand{\U}{\boldsymbol{U}}

\newcommand{\F}{\boldsymbol{F}}

\newcommand{\bn}{\mathbf{n}}
\newcommand{\dpar}[2]  {\frac{\partial #1}{\partial #2}}

\newcommand{\tens}[1]{{\mathbf{#1}}}

\newfont{\numerikEleven}{ecrm1000}
\newfont{\numerikTen}{cmss10}
\newfont{\numerikNine}{cmss9}
\newfont{\numerikEight}{cmss8}
\newfont{\numerikSeven}{cmss7}
\newfont{\numerikSix}{cmss6}

\journal{Journal of Computational Physics}

\begin{document} 
\begin{frontmatter}
\title{  ``\textit{A Posteriori}'' Limited High Order and Robust Schemes for Transient Simulations of Fluid Flows in Gas Dynamics} 

\author[zurich,polimi]{Paola Bacigaluppi$^{*}$}
\ead{paola.bacigaluppi@polimi.it}
\cortext[cor1]{Corresponding author}

\author[zurich]{R{\'e}mi Abgrall}
\ead{remi.abgrall@math.uzh.ch}

\author[losalamos]{Svetlana Tokareva}
\ead{tokareva@lanl.gov}


\affiliation[zurich]{
organization={Institute of Mathematics, University of Zurich},
addressline={Winterthurerstrasse 190}, 
city={Zurich},
postcode={8057}, 
country={Switzerland}}

\affiliation[polimi]{
organization={Department of Aerospace Science and Technology, Politecnico di Milano},
addressline={Via Privata Giuseppe La Masa 34}, 
city={Milan},
postcode={20156}, 
country={Italy}}

\affiliation[losalamos]{
organization={Applied Mathematics and Plasma Physics Group, Theoretical Division, Los Alamos National Laboratory},
addressline={PO Box 1663}, 
city={Los Alamos},
postcode={87545}, 
country={USA}}



\begin{abstract}
In this paper, we propose a novel approximation strategy for time-dependent hyperbolic systems of conservation laws for the Euler system of gas dynamics that aims to represent the dynamics of strong interacting discontinuities. The goal of our method is to allow an approximation with a high-order of accuracy in smooth regions of the flow, while ensuring robustness and a non-oscillatory behaviour in the regions of steep gradients, in particular across shocks.

Following the Multidimensional Optimal Order Detection (MOOD) (\cite{Clain2011,Diot2012}) approach,  a candidate solution is computed at a next time level via a high-order accurate explicit scheme (\cite{Abgrall2017,AbgrallHO2018}). A  so-called detector determines if the candidate solution reveals any spurious oscillations or numerical issue and, if so, only the troubled cells are locally recomputed via a more dissipative scheme. 
This allows to design a family of ``a posteriori'' limited, robust and positivity preserving, as well as high accurate, non-oscillatory and effective scheme.
Among the detecting criteria of the novel MOOD strategy, two different approaches from literature, based on the work of \cite{Clain2011,Diot2012} and of \cite{Vilar2018}, are investigated.
Numerical examples in 1D and 2D, on structured and unstructured meshes, are proposed to assess the effective order of accuracy for smooth flows, the non-oscillatory behaviour on shocked flows, the robustness and positivity preservation on more extreme flows.
\end{abstract}

\end{frontmatter}


\section{Introduction} \label{sec:introduction}

Flows displaying strongly interacting discontinuities have been intensively investigated in gas dynamics. Several mathematical models can be used to approximate the corresponding physics. The Euler equations in multiple dimensions count among these models, with hyperbolicity property being one its most important features. Euler equations allow us to represent physical flow structures such as rarefactions, contact discontinuities and shocks. Even though there exist countless approximation strategies to tackle the solution to this problem, each method designed so far encounters some restrictions such as providing a highly accurate approximation while at the same guaranteeing the robustness, i.e. the capability to provide a solution even in case of extremely tough initial conditions. Moreover, methods must not be excessively expensive, both in terms of computational time and required computer memory.  In reality, there is a trade-off between accuracy and robustness, as an accurate method might affect negatively the robustness and computational efficiency, while on the other hand, a robust and efficient approach might be limited to low accuracy.
The proposed methodology, based on Residual Distribution (RD) schemes (see \cite{Abgrall2001}, \cite{Ricchiuto2007}), represents a good compromise between these requirements, as it is designed for high order accuracy and guarantees at the same time an excellent robustness. Computational efficiency of the method is not discussed in detail in the present work, but, in general, Residual Distribution schemes can be implemented for modern parallel computing architectures. For example, their approximation stencil is compact by construction, allowing for easy application of domain decomposition approach. As for the memory efficiency, the design principle of Residual Distribution methods considered here leads to a diagonal mass matrix even for high order of accuracy, hence there is no need to store any sparse matrices unlike in classical finite element methods. 

Although the construction of first-order, robust and stable Residual Distribution (RD) schemes for steady state problem has been achieved in the 80s, the construction of  very high order accurate and robust RD schemes for unsteady problems is more recent (see \cite{Abgrall2010,AbgrallENUMATH2015,Abgrall2017,AbgrallHO2018}). 
In this paper, we consider a RD formulation based on a finite element approximation of the solution as a globally continuous piecewise polynomial.
Further, we follow \cite{AbgrallENUMATH2015,Abgrall2017}, where we have shown how one can solve a scalar version of a hyperbolic system with a method that approximates the spatial term using the RD approach, without having to solve a large linear system with a sparse mass matrix and its extension to systems was achieved in \cite{AbgrallHO2018}. This reformulation allows to avoid any mass matrix "inversion" while solving an explicit scheme. This is achieved by first approximating the time operator in a consistent way with the spatial term. A priori, this would lead either to an implicit method in case of nonlinear schemes, as done in order to avoid spurious oscillations on discontinuous solutions, or at least the inversion of a sparse but non diagonal matrix.  This apparent difficulty can be solved by applying a Deferred Correction-like time-stepping method (\cite{Abgrall2017,AbgrallHO2018}) inspired by \cite{shu-dec,Minion2} among others, but recast in a different way, and the use of proper basis functions. It has been demonstrated in \cite{AbgrallENUMATH2015,Abgrall2017} that Bernstein polynomials are a suitable choice, but this is not the only possible one. The idea to use Bernstein polynomials as shape functions instead of the more typical Lagrange polynomials, has been discussed in \cite{Abgrall2010, Abgrall2017} and applied to the context of high order Residual Distribution schemes.

The essence of the present work is to apply to the high order explicit Residual Distribution approach, presented in \cite{AbgrallHO2018}, a blending, which is designed as an ``a posteriori'' limiter, in order to ensure a highly accurate representation of the solution in the areas of smooth flows, while ensuring a non-oscillatory behaviour across strong interacting discontinuities. The main features of the designed scheme should also ensure the overall robustness and allow for a \textit{fail-safe} numerical solution for any problem under consideration.
The proposed ``a posteriori'' limiting strategy is obtained via a Multidimensional Optimal Order Detection (MOOD) method by considering a candidate solution for the next time iteration given by a non-dissipative scheme for the spatial terms, such as the standard Galerkin method, which is characterized by highly accurate approximations in smooth regions of the flow, but is not robust and displays numerical oscillations in case of strong interacting discontinuities. The ``a posteriori'' limiting detects then the cells which display any physically non admissible solutions and those cells are \textit{locally} re-approximated by discarding the candidate solution and recomputing the solution via a more dissipating spatial numerical scheme.

In literature, the MOOD approach has first emerged in a finite volume context (\hspace{1sp}\cite{Clain2011,Diot2012,Diot2013}), with the main idea to approximate solutions via high order polynomial reconstruction and in case the detecting criteria evidence any troubled cell, the polynomial degree in the associated cell is decremented, and the solution is locally recomputed.
In a second approach, the MOOD strategy has been extended to the finite element discontinuous Galerkin schemes (\hspace{1sp}\cite{Dumbser2014,Dumbser2016}), where the leading idea has been to apply unlimited discontinuous Galerkin schemes with a high approximation degree for a candidate solution, and, detected troubled cells are re-evaluated by discarding the candidate solution. In those cells, sub-cells are introduced and a more robust second order total variation diminishing (TVD) finite volume scheme is applied to update to the next time-step the sub-cell averages within the troubled DG cells. The new sub-grid data at the next time level are then gathered back into a valid cell-centred DG polynomial of degree N by using a classical conservative and higher order accurate finite volume reconstruction technique.
Recently,  in \cite{Vilar2018}, a further method similar to this last one has been proposed. There the DG reconstructed flux on the sub-cell boundaries is substituted locally, in case the detection criteria are activated, by a robust first-order or second-order TVD numerical flux.

The idea of the proposed method, while certainly inspired by these previous works, is, nevertheless on a different level. The main difference, indeed, is given by two different traits, as for instance, the novel methodology consists in the re-computation of the local troubled cell by taking the very same element typology, i.e. the polynomial degree of the considered shape functions is kept the same, and, moreover, the cell is recomputed at a global level, in the sense that we do not recompute the sub-cell values via a different approach, but simply compute the whole cell with a more dissipating scheme. By construction the guarantee of conservation is guaranted.

Concerning the detection criteria itself, we have, furthermore investigated two different approaches within this work: the first one  based on \cite{Clain2011,Diot2012,Diot2013,Dumbser2016}, while the second has been taken after \cite{Vilar2018}.

To this end, this manuscript is organized as follows. 
In Section \ref{sec:Euler} we present briefly the considered model equations and in Section \ref{sec:RD} we recall the overall framework of the considered discretizations techniques. First, we recall the basics of the RD schemes. In particular we summarize the leading traits of this approach for the steady case along its extension to the unsteady case for high order of accuracy as in \cite{AbgrallHO2018}. Section \ref{sec:MOOD} describes the generic idea of the ``a posteriori'' Multidimensional Optimal Order Detection Method. In this section we describe the detailed detection procedure, along with the two different considered detection strategies inspired by \cite{Clain2011} and \cite{Vilar2018}. Finally we provide in Section \ref{sec:numerics} several numerical benchmark problems, both in one- and two-dimensions to assess the accuracy and overall robustness of the proposed methodology and to investigate the major differences given by the two different detection criteria. The considered 2D test cases are both for structured and unstructured meshes.
Last, in the final section, we provide some conclusive remarks along with the perspectives of this work.


%
\section{Modelling Equations}
\label{sec:Euler}

Let us start by considering a classical non-linear system for hyperbolic conservations laws that describe unsteady compressible flows in multi-dimensions. It reads
\begin{equation} \label{eqn.pde.nc} 
  \frac{\partial \U}{\partial t} + \nabla \cdot \tens{\F}(\U)  = \bm{0} , \qquad \xx \in \Omega \subset \mathbb{R}^2, \; t \in \mathbb{R}_0^+, 
\end{equation} 
with appropriate initial and boundary conditions. Here $\xx=(x,y)$ is the coordinate vector within the computational domain $\Omega$. 
In particular, let us choose the Euler equations, such that $\U=[\rho, \rho \mathbf{u}, E]^T$ is the vector of unknown conserved variables, 
$\tens{\F} = [\rho \mathbf{u}, \rho \mathbf{u}\otimes \mathbf{u}+ P \mathbf{I},\mathbf{u} ( E + P )]^T$ is the conservative non-linear flux tensor depending on $\U$.
Within this modelling equations, we have $\rho$ that denotes the mass density, $\mathbf{u}=(u,v)$ the velocity vector, $P$ the fluid pressure, $E$ the total energy, $\mathbf{I}$ the $2 \times 2$ identity matrix and $\mathbf{u} \otimes \mathbf{u}$ is the dyadic product of the velocity vector with itself.
The perfect gas law is set as equation of state (EOS) to close system \eqref{eqn.pde.nc}, such that $P = (\gamma-1)\left( E - \frac{1}{2}\rho \mathbf{u}^2 \right)$, with $\gamma$ the ratio of specific heats. 
The sound speed is defined as $c=\sqrt{\gamma P/\rho}$.\\
Physically admissible states
are those such that $\rho>0$ and $P>0$.\\
 
%
%
\section{Discretization Strategy}  \label{sec:RD}
\subsection{Spatial Discretization: Galerkin Finite Element-type}
In order to provide a complete overview of the numerical approximation strategy adopted throughout this work, we recall hereafter an overview on the finite-element-type Residual Distribution scheme, following \cite{mariohdr,AbgrallHO2018}. The reader may also refer to \cite{SWjcp,Abgrall2006,Ricchiuto2007} for further details on the construction of generic residual distribution schemes.

\begin{figure}[H]
  \begin{center}
    \includegraphics[width=0.25\textwidth]{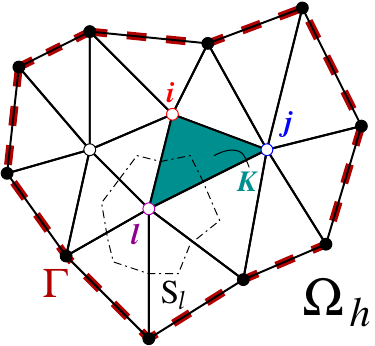}
    \caption{ \label{fig:notation}
      Notation for domain $\Omega_h$, its boundary $\Gamma$, cells and median dual cells.
    }
  \end{center}
\end{figure}

Let us introduce first the solution approximation space $V_h$ that follows the classical Galerkin Finite Element approach given by globally continuous polynomials of degree $k$, such that
\begin{equation}\label{eq:Vh}
V_h= \left\{ \U \in  C^0(\Omega_h), \; \; \U_{|K} \in \mathbb{P}^k, \; \; \forall K \in \Omega_h  \right\}.
\end{equation}
$\Omega_h$ corresponds to the set of conformal, non-overlapping elements of characteristic length $h$, obtained from the discretization of the computational domain $\Omega$.  
The generic element is called $K$ and the volume of a cell $|K|$, which is for 1D the length of the cell, while in 2D the area. 
We denote by $S_{\sigma}$ the standard median dual cell obtained by joining the gravity centres of the elements in $K_\sigma$ with the mid-points of the edges emanating from $\sigma$ whose area is given by
\begin{equation}\label{S_sigma}
|S_\sigma|=\frac{1}{N_{DoF}}\sum_{K, \sigma \in K} |K|\,
\end{equation}
 with the total number of DoFs in one cell is $N_{DoF}$. See figure~\ref{fig:notation} for an illustration,  where $\sigma$ corresponds to the degree of freedom $l$ and hence $S_l$.
Further, in section \ref{subsec:detection} we denote by $\mathcal{V}_K$ the set of all neighbour cells and by $\mathcal{W}_K$ the set of neighbour cells sharing one edge.\\
The time domain $[0,T]$ is approximated by a set of time intervals $[t_n,\,t_{n+1}]$ with $\Delta t =t_{n+1}-t_n$ the time step.

The numerical solution $\U_h^n \simeq \U(\xx,t^n)$ is represented by
\begin{equation} \label{approx_uh}
  \U_h^n(\xx) = \sum_{\sigma \in \Omega_h}  \U_\sigma^n \, \varphi_\sigma(\xx) , 
  \quad \xx \in \Omega,
\end{equation}
that corresponds, on each element, to a linear combination of the  shape functions 
$\varphi_\sigma \in V_h$,
which are assumed to be continuous within the elements and on the faces of the elements
with coefficients $\U_\sigma^n$ to be determined by a numerical method. 
In this work we employ Bernstein basis functions $\left\{ \varphi_{\sigma} \right\}_{\sigma}$ of order $k$. 
While one could observe that using Bernstein polynomials has the drawback that not all degrees of freedom $\U_{\sigma}^n$ in the expansion \eqref{approx_uh} will represent the solution values at certain nodes, one of the main advantages offered by this type of shape functions is their positivity on a cell $K$ that will enforce, the positivity of the mass matrix, as we shall see in the coming section.

For completeness, we recall the definition of Bernstein polynomials on triangular elements:
\begin{multicols}{2}
\begin{itemize}
	\item Linear ('$\mathcal{B}^1$'): 
	\begin{equation*}
	\varphi_1 = x_1, \ \varphi_2 = x_2, \ \varphi_3 = x_3. \quad \quad \quad \quad
	\end{equation*}
	
	\item Quadratic ('$\mathcal{B}^2$'):
	\begin{align*}
	& \varphi_1 = x_1^2, \ \varphi_2 = x_2^2, \ \varphi_3= x_3^2, \\
	& \varphi_4 = 2 x_1 x_2, \ \varphi_5 = 2 x_2 x_3, \ \varphi_6 = 2 x_1 x_3.
	\end{align*}

	\item Cubic ('$\mathcal{B}^3$'):
	\begin{align*}
	& \varphi_1 = x_1^3, \ \varphi_2 = x_2^3, \ \varphi_3= x_3^3, \\
	& \varphi_4 = 3 x_1^2 x_2, \ \varphi_5 = 3 x_1 x_2^2, \ \varphi_6 = 3 x_2^2 x_3, \\
	& \varphi_7 = 3 x_2 x_3^2, \ \varphi_8 = 3 x_1 x_3^2, \ \varphi_9 = 3 x_1^2 x_3, \\
	& \varphi_{10} = 6 x_1 x_2 x_3.
	\end{align*}
	
\end{itemize}
	
	\begin{figure}[H]
	\centering
\includegraphics[scale=0.3]{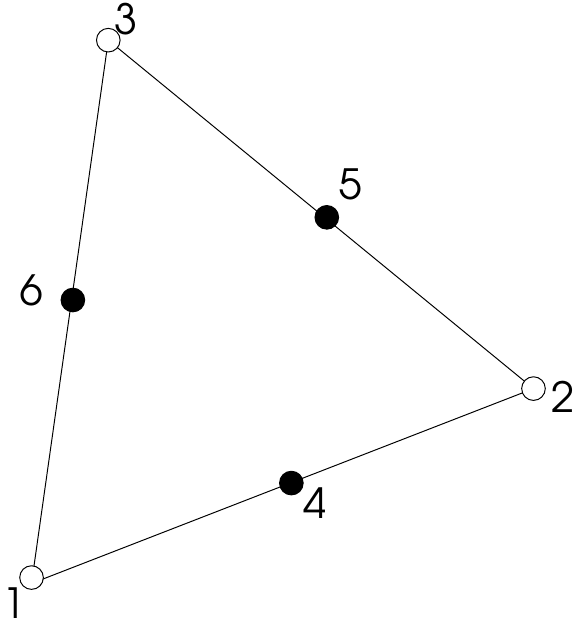}\\
\includegraphics[scale=0.3]{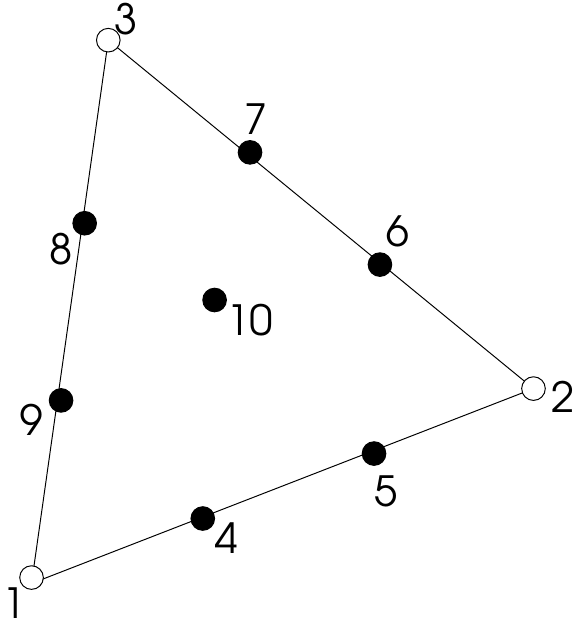}
\caption{Nomenclature of the DoFs within a $\mathcal{B}^2$ (upper triangle) and a $\mathcal{B}^3$ (lower triangle) element. }
\label{Bernstein_dof}
\end{figure}
\end{multicols}
\noindent Here, the barycentric coordinates are defined by $x_1$, $x_2$, $x_3.$
The location of the degrees of freedom for the third and fourth order are, moreover, shown in Fig.~\ref{Bernstein_dof}. 

The approximation of $\tens{F}(\U_h)$ in cell $K$ can be done in two possible ways, as explained and commented in \cite{AbgrallHO2018}, too. One can either evaluate the values of the flux at the DoFs from the data $U_h$, defining $\tens{\F}(\U_h)$ as:
\begin{equation}
\label{flux:approximation}
\tens{\F}(\U_h)\approx \sum_{\sigma\in K}\tens{\F}_\sigma \varphi_\sigma,
\end{equation}
where $\tens{\F}_\sigma$ is the degree of freedom for the flux in the Bernstein basis, which leads to a quadrature-free implementation since the integrals of the shape functions and/or gradients can be evaluated explicitly. Alternatively, one can define $\tens{\F}(\U_h)$ as the flux evaluated for the local value of $\U_h$ at the quadrature point, since both approaches are formally equivalent from the accuracy point of view.

\subsection{Spatial Residual Distribution scheme}
The idea of working with Residual Distribution schemes resides in its flexibility of treating informations with the advantage of increasing only locally within a cell the degrees of freedom in a DG fashion, while keeping the amount of informations at the same level of finite element schemes. To allow an insight, let us consider the steady version of system \eqref{eqn.pde.nc}. Its integral on a cell $K\in \Omega_h$ of $\U$ is defined as
\begin{equation} \label{eq:residual}
  \Phi^K \left( \U \right) =   \Int_{K} \nabla \cdot  \F(\U)  \, d\xx = \Int_{\partial K} \F(\U) \cdot \bm{n}  \, ds,
\end{equation}
and we define it hereafter as the \textit{total residual} $\Phi^K$. 
Looking at \eqref{eq:residual} from  a different perspective, the total residual is composed by the sum of the contributions $\phi_{\sigma}^K$ from each degree of freedom $\sigma$ within the cell $K$ and the following conservation property holds true
\begin{equation} \label{eq:fluctuation}
  \Phi^K \left( \U_h \right)  = \Sum_{\sigma \in K} \phi_\sigma^K, \; \forall K\in \Omega_h .
\end{equation}
Plugging the definition \eqref{eq:fluctuation} within the considered system of equations for the steady case, i.e. $\nabla \cdot \F(\U)=\mathbf{0}$, we have that \begin{equation}\label{residualSys}
\Phi^K \left( \U_h \right)=0.\end{equation} 
In particular, in case the considered degree of freedom $\sigma$ would belong to the physical boundary $\Gamma$, \eqref{residualSys} would be split into two contributions, one for the internal and one for the boundary, to be \begin{equation} \label{eq:phiconsBC}
  \Sum_{K, \sigma\in K} \phi_\sigma^K + \Sum_{\gamma\in \Gamma, \sigma\in \gamma} \phi_\sigma^\gamma = 0, \qquad \forall \sigma \in \Gamma,
\end{equation}
where $\gamma$ is any edge on the boundary $\Gamma$ of $\Omega_h$ (see \cite{AbgrallHO2018} for further details and figure~\ref{fig:notation} for an illustration).

\subsubsection*{In practice:}
Having provided the context, let us now take the above definitions for Residual Distributions from bottom-up. 
The actual way, the strategy is carried out, can be summarized as follows. First of all we chose a numerical spatial scheme to determine $\Phi^K$. An RD numerical scheme is entirely determined by the strategy with which one
distributes the the total residual amongst the degrees of freedom $\sigma$ .
This strategy 
in \eqref{eq:fluctuation} is defined by means of a distribution coefficients $\beta_\sigma$ as
\begin{equation} \label{eq:beta}
  \phi_\sigma^K = \beta_\sigma^K \Phi^K ,  \; \forall\sigma \in K .
\end{equation}
Hence one specific RD scheme is defined for each and every set of parameters $\beta_{\sigma}$.
The conservation property \eqref{eq:fluctuation} enforces that
\begin{equation} \label{eq:beta_cons}
  \Sum_{\sigma \in K} \beta_\sigma^K  = 1 .
\end{equation}
The final RD scheme results from collecting the residuals $\phi_\sigma^K$
from cells surrounding the point associated to the specific DoF $\sigma$, that is
\begin{equation} \label{eq:phicons}
\Sum_{K, \sigma\in K} \phi_\sigma^K = 0, \qquad \forall \sigma \in \Omega_h,
\end{equation}
which allows to compute the unknown coefficients of the polynomial solution 
$\U_\sigma$ in \eqref{approx_uh}. 

\begin{figure}
  \begin{center}
    \includegraphics[width=0.3\textwidth]{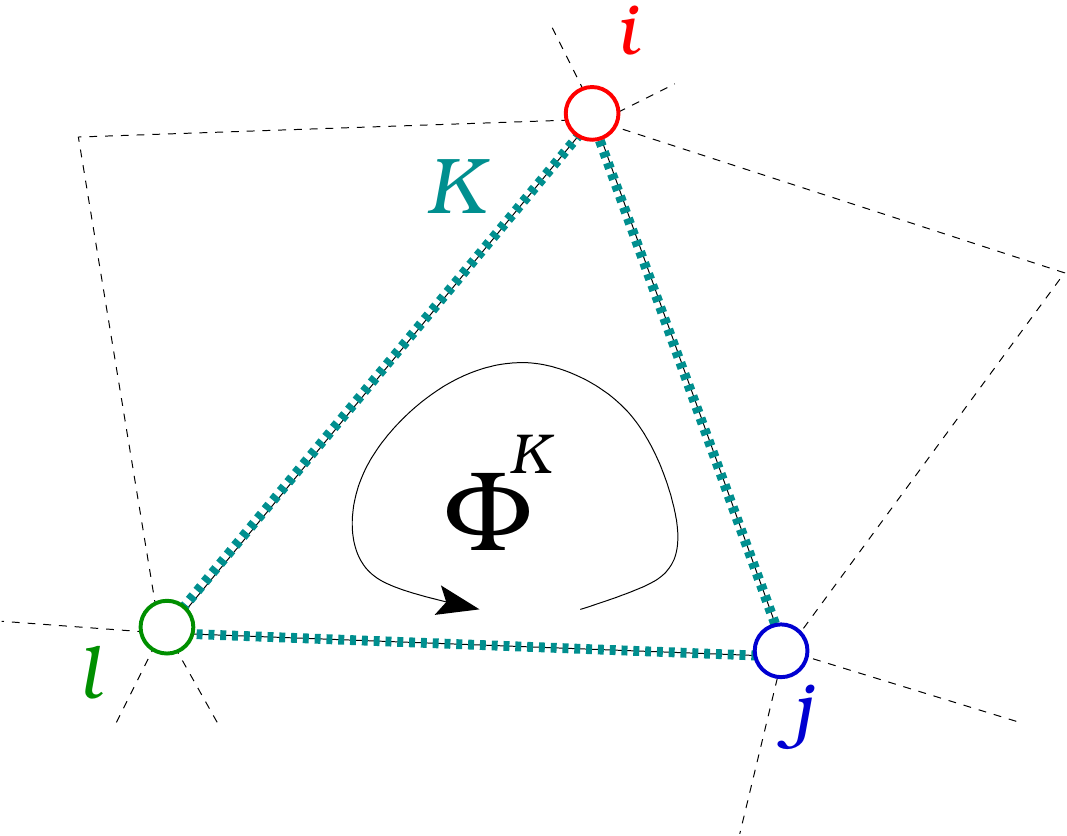}
        \includegraphics[width=0.3\textwidth]{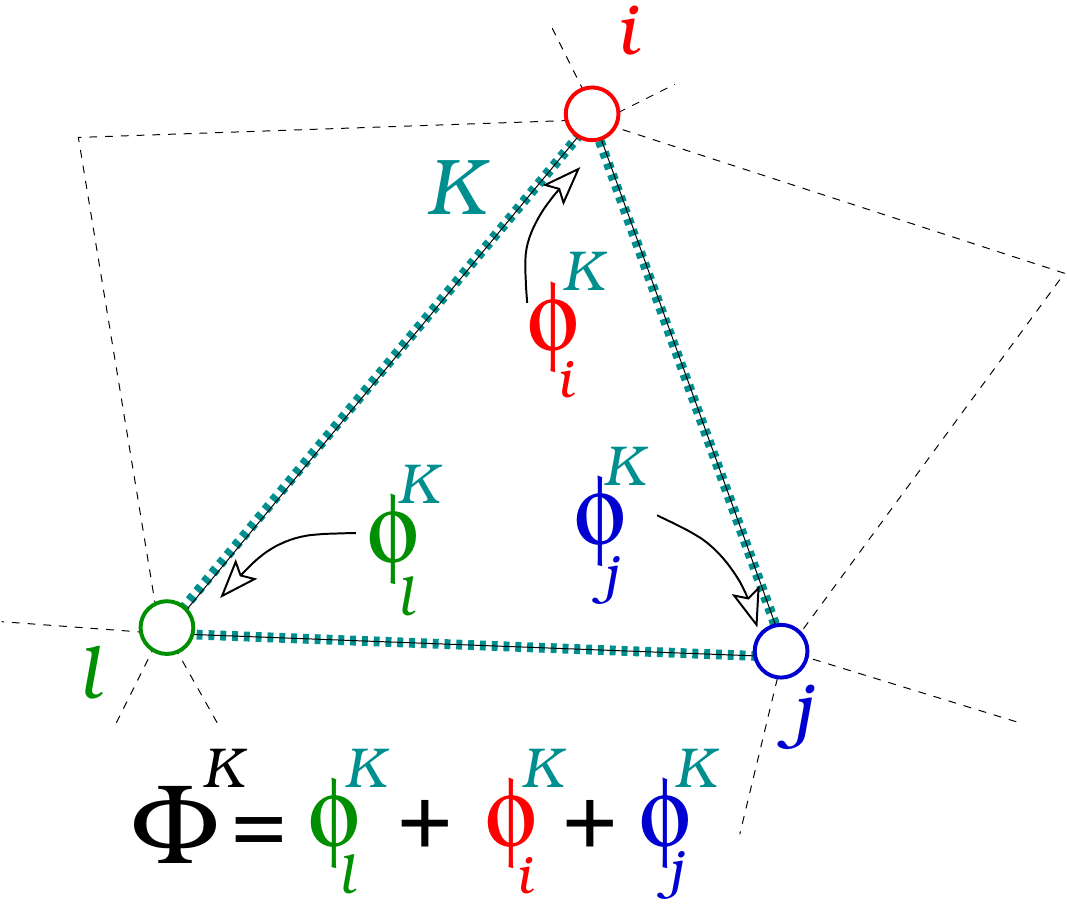}
            \includegraphics[width=0.3\textwidth]{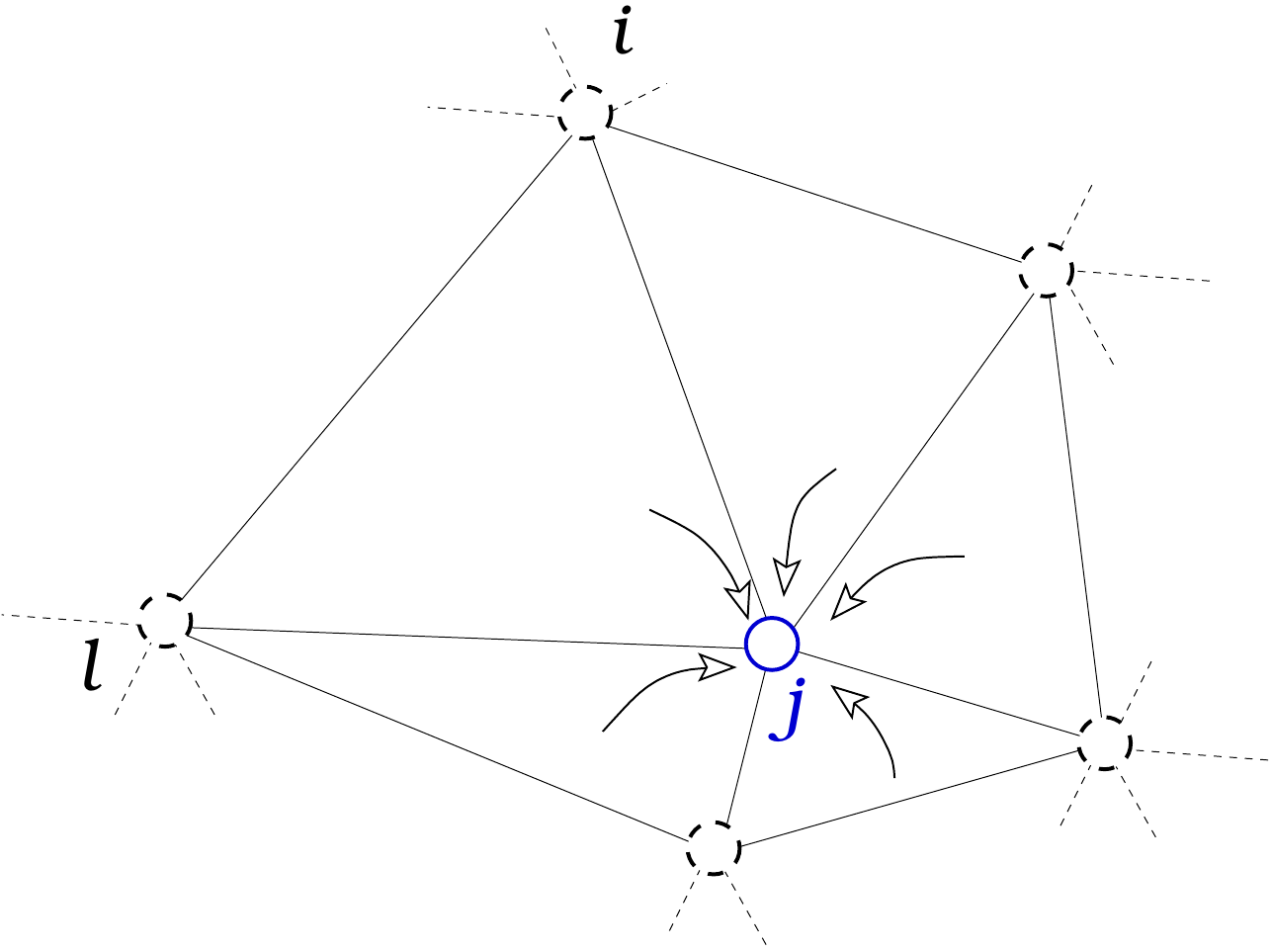}
    \caption{ \label{fig:RD}
      Illustration of the three steps of the residual distribution approach: Compute total residual, Distribute
      them amongst the degrees of freedom and finally Gather them.
    }
  \end{center}
\end{figure}

In figure~\ref{fig:RD} we illustrate where $\Phi^K$ is actually defined, and how it is composed, where, for simplicity we only plot the degrees of freedom which coincide with the vertices of $K$. Finally,  the gathering of the residuals around the DoF is shown.

\subsection{Time Discretization: Deferred Correction-like Method} \label{ssec.TimeRD}
To provide a complete frame within which the novel strategy of this paper is collocated, we briefly recall the construction of the high-order accurate time discretisation for residual distribution schemes.  This timestepping approach is similar but not identical to the Deferred Correction methodology of \cite{dutt2000,Minion2,shu-dec} and was first proposed in \cite{Abgrall2017} and then extended to multidimensional systems in \cite{AbgrallHO2018}, with recent studies showing analogies to e.g. ADER (\cite{han2021dec}) and Runge-Kutta schemes (\cite{AbgrallMeledo2022}), extending studies also towards stabilization requirements (\cite{michel2021spectral,michel2022spectral}).

Starting by considering the numerical solution at a discrete time $t^n$ denoted by $\U^n$, we search for the 
solution at $t^{n+1}$ denoted by $\U^{n+1}$.
We split each interval $[t^n, t^{n+1}]$ into sub-timesteps 
$t^n\equiv t^{n,0} < t^{n,1} < \ldots < t^{n,m} < \ldots < t^{n,M} \equiv t^{n+1}$.
For the $m$-th subinterval $[t^{n,m},t^{n,m+1}]$, we introduce the correction indices
$r= 0, \ldots, R$ and further denote the solution at correction index $r$ of the sub-timestep $m$ 
by $\U^{n,m,r}$. In addition we denote the solution vector of the $r$-th corrections by
 system \eqref{eqn.pde.nc} can be formally integrated on $[t^n, t^{n+1}]$ as
\begin{equation}\label{eq:Ucorr}
\U^{(r)} = \left( \U^{n,0,r}, \ldots, \U^{n,M,r}  \right).
\end{equation} 

Let us proceed  within the time step $[t^n, t^{n+1}]$ as follows:
\begin{enumerate}
\item We initialize for $m=1,\ldots, M$: $ \U^{n,m,0} =  \U^n$;
\item for each correction $r=0,\ldots,R-1$, we iterate for $m=1,\ldots, M$, such that,\\
knowing $\U^{n,m,r}$, we evaluate $\U^{n,m,r+1}$ as the solution of 
\begin{equation}\label{eq:defcor}
\mathcal{L}^1(\U^{n,m,r+1})=\mathcal{L}^1(\U^{n,m,r})-\mathcal{L}^2(\U^{n,m,r})
\end{equation}
\item set $\U^{n+1}=\U^{n,M,R}.$
\end{enumerate}

In \cite{Abgrall2017} it has been proved that, under some assumption on the differential operators
$\mathcal{L}_\Delta^1,\mathcal{L}_\Delta^2$ depending on a parameter $\Delta$, the deferred correction-like method is convergent, and after $R$ iterations
the error is smaller than $\nu^R||\U^{(0)}-\U_{\Delta}^{\star}||$, where $\nu=\frac{\alpha_2}{\alpha_1} \Delta <1$ 
is a real constant depending on the operators $\mathcal{L}_\Delta^1,\mathcal{L}_\Delta^2$.\\
The $\alpha_1$ and $\alpha_2$ are assumed to be the parameters of the operators such that
$\mathcal{L}_\Delta^2$ has a unique root $\U_\Delta^\star$ (such that $\mathcal{L}_\Delta^2(\U_\Delta^\star)=0$), 
$\mathcal{L}_\Delta^1$ is coercive with coercivity constant $\alpha_1$, 
and $\mathcal{L}_\Delta^2-\mathcal{L}_\Delta^1$ is uniformly Lipschitz continuous with Lipschitz constant $\alpha_2 \Delta$. 

Let us first recall the low order differential operator $\mathcal{L}^1$ then 
the high order differential operator $\mathcal{L}^2$.
To do so the starting point consists in recalling that
 system \eqref{eqn.pde.nc} can be formally integrated on $[t^n, t^{n+1}]$ as
\begin{equation}
\U(\xx,t^{n+1}) = \U(\xx,t^{n}) + \Int_{t^n}^{t^{n+1}} \nabla\cdot \tens{F}( \U (x,t) ) \, dt,
\end{equation}
and the solution be approximated with a suitable quadrature rule
\begin{equation}
\U(\xx,t^{n+1}) \simeq \U(\xx,t^{n}) + \Delta t  \Sum_{l} \omega_l \, \nabla\cdot \tens{F}( \U (x,t_l) ).
\end{equation}

\subsubsection{Low Order Differential Operator $\mathcal{L}^1$}
For any $\sigma \in K$, we define $\mathcal{L}^1$ as
\begin{equation}\label{eq:L1}
\mathcal{L}^1_\sigma ( \U^{(r)} ) =\mathcal{L}^1_\sigma ( \U^{n,1,r}, \ldots , \U^{n,M,r})
\end{equation}
where in particular we have for a generic sub-timestep $m$ at a correction $r$ that
\begin{equation}
\mathcal{L}^1_\sigma (\U^{n,m,r})=|C_\sigma| \left(  \U_\sigma^{n,m,r} - \U_\sigma^{n,0,r} \right) + \Sum_{K | \sigma\in K} \Int_{t^{n,0}}^{t^{n,m}} \mathcal{I}_0 \left( \phi_\sigma^K(\U^{(r)}) ,t\right)\,dt, 
\end{equation}

where $\mathcal{I}_0$ represents any first order piecewise-constant interpolant under the following notation
\begin{equation}
 \phi_\sigma^K(\U^{(r)}) = \big( \phi_\sigma^K(\U^{n,1,r}),  \ldots , \phi_\sigma^K(\U^{n,M,r}) \big).
\end{equation} 
As it stands, system \eqref{eq:L1} is a time implicit system. An explicit in time version is
obtained by considering $\mathcal{I}_0$ as being a simple approximation of $\U^{n,0}$ for all $m$, so that \eqref{eq:L1} becomes
\begin{equation}\label{eq:L1b}
\mathcal{L}^1_\sigma (\U^{n,m,r})=
|C_\sigma| \left(  \U_\sigma^{n,m,r} - \U_\sigma^{n,0,r} \right) + \Delta t \, \xi_m \, \Sum_{K | \sigma\in K} \phi_\sigma^K(\U^{(0)})  
\end{equation}
with $\xi_m$ satisfying for $m=1,\ldots,M$  $t_{n,m}=t_n+\xi_m \Delta t$ and $0=\xi_0<\ldots<\xi_m<\xi_{m+1}<\ldots \xi_M=1$ within the considered time interval $[t_n, t_{n+1}]$.\\
The coefficients $|C_\sigma|$ play the role of the dual cell measure and and in order for (\ref{eq:L1}) to be solvable we have to satisfy the constraint
\begin{equation}\label{Csigma_mood}
|C_\sigma| = \Int_K \varphi_\sigma( \xx) \, d\xx >0,
\end{equation}
This requirement has a direct consequence on the choice of the polynomial basis $\left\{ \varphi_\sigma \right\}_\sigma$. 
For instance, the classical Lagrange basis on simplex is disqualified as it does not verify \eqref{Csigma_mood} for $k>1$, and this is the reason why we consider Bernstein polynomials \cite{Bernstein1912,Lorentz_book_Bersntein53} for high order approximations. Indeed, Bernstein basis functions verify $\varphi_\sigma(\xx)\geq 0$ for all $\sigma$ and $\xx$, and $\Sum_{\sigma} \varphi_\sigma(\xx)=1$
for all $\xx$, and it is easy to deduce that \eqref{Csigma_mood} is fulfilled. 
The low order differential operator $\mathcal{L}^1_\sigma$ constructed this way is then of high accuracy in space and explicit in time.

\subsubsection{High Order Differential Operator $\mathcal{L}^2$}
The high order differential operator $\mathcal{L}^2_\sigma$ for a generic sub-timestep $m$ at a correction $r$ is defined as
\begin{equation}\label{eq:L2}
\mathcal{L}^2_\sigma ( \U^{n,m,r} ) =
\Sum_{K | \sigma\in K} \left( 
\Int_K \varphi_\sigma ( \U^{n,m,r} - \U^{n,0,r}) \, d\xx +  \Int_{t^{n,0}}^{t^{n,m}} \mathcal{I}_M \left( \phi_\sigma^K(\U^{(r)}) ,s\right)\,ds   \right)  
\end{equation}
Once the coefficients of the interpolating polynomial $\mathcal{I}_M$ of degree $M$ are computed, we perform
the exact integration to obtain the approximation for every row of \eqref{eq:L2} in the form
\begin{equation}
 \Int_{t^{n,0}}^{t^{n,m}} \mathcal{I}_M \left( \phi_\sigma^K(\U^{n,0,r}),\ldots , \phi_\sigma^K(\U^{n,M,r}) ,s \right)\,ds 
=  \Sum_{l=0}^M \zeta_{l,m}  \phi_\sigma^K(\U^{n,l,r}) ,
\end{equation}
where $\zeta_{l,  m}$ are approximation coefficients.
The high order differential operator $\mathcal{L}^2$ ensures a high order approximation
of the space-time term $\partial_t \U + \nabla \cdot \tens{F}(\U)$. 
The operator \eqref{eq:L2} is implicit in time.
Nevertheless, adopting the combination of the two operators $\mathcal{L}^2$ and $\mathcal{L}^1$, as in equation \eqref{eq:defcor}, we retrieve an \textit{explicit} formulation.

\begin{remark}
Having defined the low order and high order differential operators, we provide hereafter an example of the 3-step algorithm provided in Section \ref{ssec.TimeRD} for second order ($M=2$). 
\begin{enumerate}

\item Initialize for $m=1,2$: $ \U^{n,m,0} =  \U^n$
\item Apply the correction $r=0,1$ for $m=1,2$, on \eqref{eq:defcor}.
\begin{itemize}
\item For $r=0$ and $m=1,2$ we have:
\begin{equation}
\begin{split}
|C_\sigma| \left(  \U_\sigma^{n,m,1} - \U_\sigma^{n,0,1} \right) + \Delta t \, \xi_m \, \Sum_{K | \sigma\in K} \phi_\sigma^K(\U^{(0)})  
= &|C_\sigma| \left(  \U_\sigma^{n,m,0} - \U_\sigma^{n,0,0} \right) + \Delta t \, \xi_m \, \Sum_{K | \sigma\in K} \phi_\sigma^K(\U^{(0)})  
\\&
-\Sum_{K | \sigma\in K} \left( 
\Int_K \varphi_\sigma ( \U^{n,m,0} - \U^{n,0,0}) \, d\xx +  \Int_{t^{n,0}}^{t^{n,m}} \mathcal{I}_2 \left( \phi_\sigma^K(\U^{(0)}) ,s\right)\,ds   \right)  
\end{split}
\end{equation}
As $\U^{n,0,0} =\U^{n,1,0} = \U^{n,2,0}=\U^n$ and $\U^{n,0,1}=\U^n$, we can simplify it to
\begin{equation}
\begin{split}
 \U_\sigma^{n,m,1}= \U_\sigma^{n,0,1}
- \frac{1}{|C_\sigma|}  \Sum_{K | \sigma\in K} \left( 
 \Int_{t^{n,0}}^{t^{n,m}} \mathcal{I}_2 \left( \phi_\sigma^K(\U^{(0)}) ,s\right)\,ds   \right)  
\end{split}
\end{equation}
\item  For $r=1$ and $m=1,2$ we have, similarly:
\begin{equation}
\begin{split}
 \U_\sigma^{n,m,2}= \U_\sigma^{n,m,1} 
- \frac{1}{|C_\sigma|}  \Sum_{K | \sigma\in K} \left( \Int_K \varphi_\sigma ( \U^{n,m,1} - \U^{n,0,1}) \, d\xx + 
 \Int_{t^{n,0}}^{t^{n,m}} \mathcal{I}_2 \left( \phi_\sigma^K(\U^{(0)}) ,s\right)\,ds   \right)  .
\end{split}
\end{equation}
\end{itemize} 
\item Set $\U^{n+1}=\U^{n,2,2}$.
\end{enumerate}

\end{remark}

\begin{remark}
As a matter of completeness, one should note that both \eqref{eq:L1} and \eqref{eq:L2} make use of  residuals $\phi^K_\sigma( \U^{n,m,r} )$ for all $\sigma \in K$, and all $K$, which are computed via any spatial residual distribution scheme as seen in the previous sections, and where the boundary conditions are applied within the residual computation. More details on the actual spatial discretization scheme will be specified in the forthcoming sections, as it is part of the main idea behind this work.
\end{remark}

%

%
%
\section{``A Posteriori'' Multidimensional Optimal Order Detection (MOOD) Method }  \label{sec:MOOD}
\begin{figure}
  \begin{center}
    \includegraphics[width=0.7\textwidth]{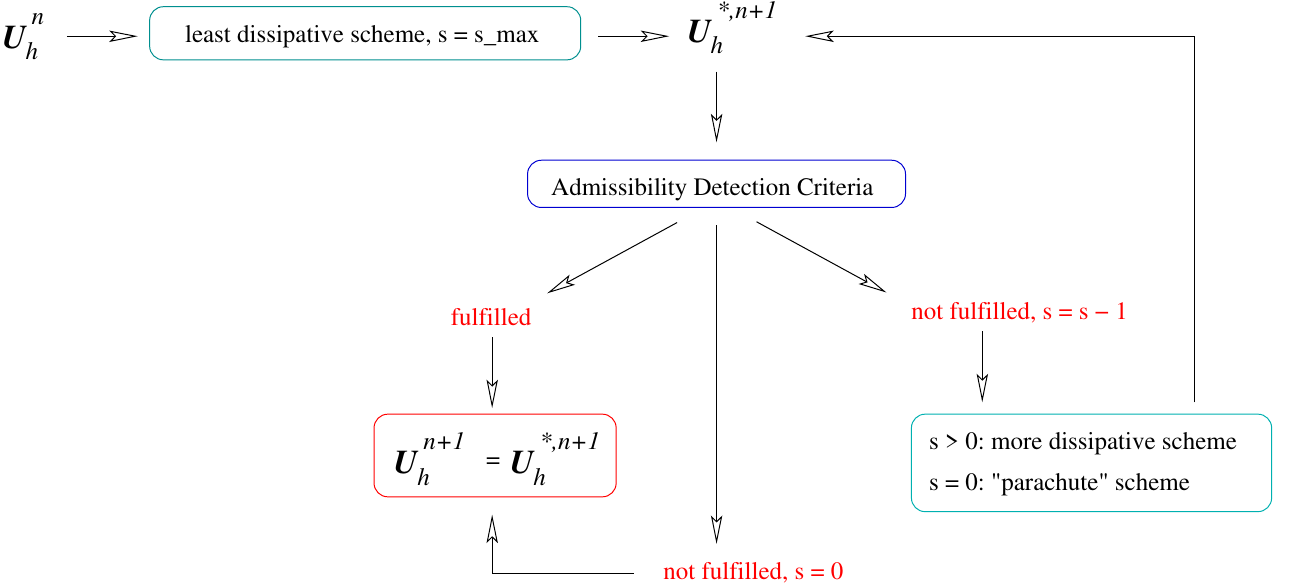}
    \caption{ \label{fig:mood}
     Loop highlighting the novel methodology with an ''a posteriori`` limiting on the RD scheme.
    }
  \end{center}
\end{figure}
\subsection{Basics: Detect, Decrement, Re-compute}\label{ssec:MOOD_basics}
The design of the residual $\phi_{\sigma}(\U^{n,m,r})$ within \eqref{eq:L2} is the main goal of this manuscript. Let us illustrate the proposed methodology via figure~\ref{fig:mood}. In a classical, ``a priori'' RD approach, as for example also done in \cite{AbgrallHO2018}, the limiting is comprised within the chosen numerical scheme. The quality of the solution entirely depends on our ability to predict how the RD scheme
behaves and when and where it fails to do so. If some inappropriate numerical data are generated in some cells, then the solution $\U_h^{n+1}$ will be imprinted without any chance to go back in time to 
possibly cure this situation. This issue is tackled within our novel approach, as follows. 

In our ``a posteriori'' technique, we start at a given time step $t_n$ with a valid solution $\U_h^n$ and compute a candidate solution for $t^{n+1}$ as $\U_h^{\star,n+1}=\textbf{RD}(\U_h^n )$ with a certain chosen spatial scheme. 
If some bad numerical data is observed in some cells, then the solution is discarded 
in those cells, and the solution is recomputed starting again from valid data at $t^n$
but using a more appropriate scheme, for instance a more robust one.
As such, some challenging situations like a lack of positivity, invalid data (NaN, Inf) 
or more classical spurious oscillations, if detected, can be handled.
\begin{figure}[t!]
  \begin{center}
 \includegraphics[width=0.7\textwidth]{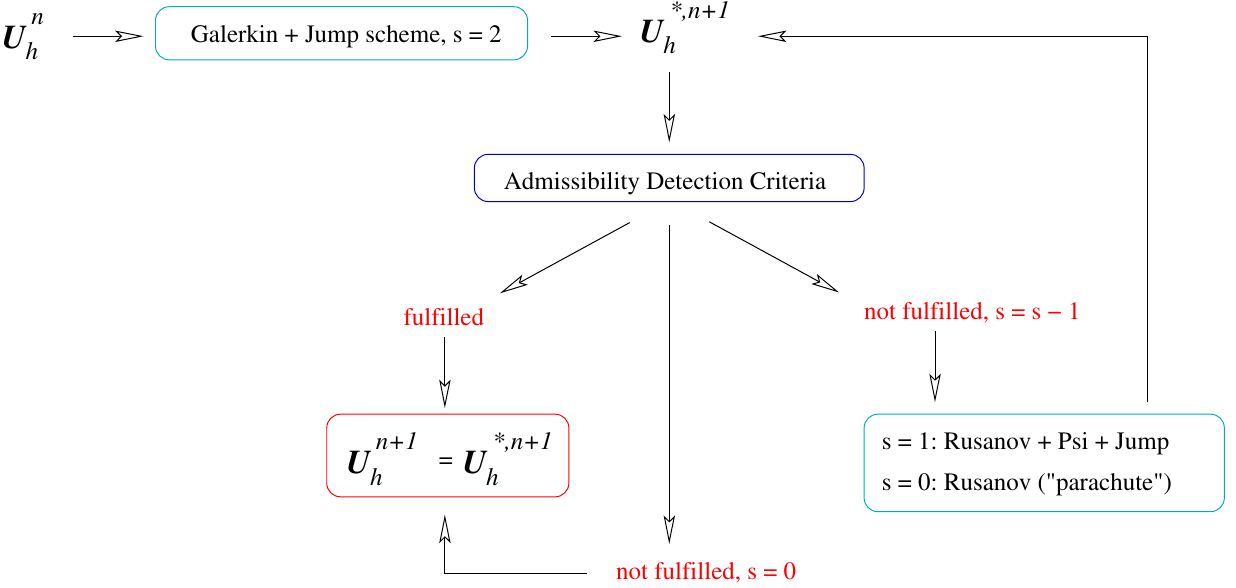}
      \caption{ \label{fig:mood_cascade}  Loop highlighting the novel methodology with the ''a posteriori`` limiting on the RD scheme, highlighting the numerical schemes used in this work.
    }
  \end{center}
\end{figure}
More specifically, the idea behind this, following figure \ref{fig:mood_cascade}, is to:
\begin{enumerate}
\item Start to compute a first candidate solution with the least dissipative and most accurate scheme possible, which we denote by $s=s_{max}$, which corresponds to a scheme providing,  for example, high-order of accuracy in case of smooth flows. In case  a cell is flagged as troubled, a more dissipative scheme $s=s_{max}-1$ is applied locally, in order to guarantee more robustness at cost of  some accuracy. 
If the forthcoming checks again outline troubled cells, the scheme is further ``decremented'' locally until, a so-called parachute scheme with $s=0$, which is first order accurate, and thus the most dissipative and more robust.
In this work, $s=2$ corresponds to the least dissipative scheme, which has been chosen to be a Galerkin method with some stabilizing terms, defined hereafter as $\phi_{\sigma,\mathbf{x}}^{K,jump}(\U_h)$, where $\mathbf{x}$ denotes the spatial part of the corresponding residual and reads:
\begin{equation}
\phi_{\sigma,\mathbf{x}}^{K,Galerkin+Jump}(\U_h)=\int_{\partial K}\varphi_\sigma \mathbf{F}(\U_h)\cdot \mathbf{n}\,d\Gamma -\int_K \nabla \varphi_\sigma\cdot \mathbf{F}(\U_h)\,d \mathbf{x} + \phi_{\sigma,\mathbf{x}}^{K,jump}(\U_h).
 \label{GalerkinJump} 
\end{equation} 
As in \cite{AbgrallHO2018}, the jump stabilization term reads
\begin{equation}
\phi_{\sigma,\mathbf{x}}^{K,jump}(\U_h)= \sum_{\text{edges of }K}
    \theta_1 h_e^2 \int_e [\nabla \U_h]\cdot [\nabla \varphi_\sigma]\,d\gamma+\sum_{\text{edges of }K}
    \theta_2 h_e^4 \int_e [\nabla^2 \U_h \mathbf{n}]\cdot [\nabla^2 \varphi_\sigma\mathbf{n}]\,d\gamma
    \label{phi_burman}
\end{equation} 
where we denote $\left[ \nabla \psi \right] = \nabla \psi|_K - \nabla \psi |_{K'}$ with $e=K \cap K'$
for any function $\psi$ and where $\mathbf{n}$ is a normal to $e$.  This scheme is referred as the Galerkin plus Jump method, or GJ.

\item In case the detection is  activated, we keep the candidate solution for $t_{n+1}$ obtained via  \eqref{GalerkinJump} for all elements $K$ except for those flagged as troubled, where we take as $s=1$ a Rusanov (i.e. a local Lax-Friedrichs) PSI scheme with some stabilizing terms, which is exactly the same scheme adopted in \cite{AbgrallHO2018} and reads
\begin{equation}\label{RusPsiJump}
 \phi_{i,\mathbf{x}}^{K,Rus+Psi+Jump}(\U_h)=\phi_{\sigma,\mathbf{x}}^{K,Rus\star}(\U_h)+ \phi_{\sigma,\mathbf{x}}^{K,jump}(\U_h) .
\end{equation}
Here, $\phi_{\sigma,\mathbf{x}}^{K,Rus\star}(\U_h)$ corresponds to a Rusanov (local Lax-Friedrichs) scheme with a psi-like blended,
first designed in \cite{DeSantis2015} and, recalling the description in \cite{AbgrallHO2018}, is written in local characteristic variables by projecting the first order residuals onto a space of left eigenvalues, as
\begin{equation}
\label{phi_char}
\hat{\phi}_{\sigma,\mathbf{x}}^{K,Rus} = \mathbf{L}\,\phi_{\sigma,\mathbf{x}}^{K,Rus}.
\end{equation}
The distribution coefficients are
\begin{equation}
\beta_\sigma^K = \dfrac{\max\left( \dfrac{\hat{\phi}_{\sigma,\mathbf{x}}^{K,Rus}}{\hat{\phi}_{\sigma,\mathbf{x}}^K},0 \right)}{\sum\limits_{\sigma'\in K} \max\left( \dfrac{\hat{\phi}_{\sigma',\mathbf{x}}j^{K,Rus}}{\hat{\phi}_{\mathbf{x}}^K},0 \right)}, \quad \hat{\phi}_{\mathbf{x}}^K = \sum_{\sigma\in K} \hat{\phi}_{\sigma,\mathbf{x}}^{K,Rus}.
\label{beta_char}
\end{equation}
and applying the blending scheme
\begin{equation}
\hat{\phi}_{\sigma,\mathbf{x}}^{K,Rus\star} = (1-\Theta)\,\beta_\sigma^K \hat{\phi}_{\mathbf{x}}^K + \Theta\,\hat{\phi}_{\sigma,\mathbf{x}}^{K,Rus},
\label{limiting}
\end{equation} 
where the blending coefficient $\Theta$ is defined by
\begin{equation}
\label{theta}
\Theta = \dfrac{\big| \hat{\phi}_{\mathbf{x}}^K \big|}{\sum\limits_{\sigma'\in K} \big| \hat{\phi}_{\sigma', \mathbf{x}}^{K,Rus} \big|}.
\end{equation}
 with $0 \leq \Theta \leq 1$, and $\Theta = O(h)$ for a smooth solution, thus ensuring accuracy and $\Theta = O(1)$ at the discontinuity, thus ensuring monotonicity \cite{Abgrall2006}.
Finally, the high-order nodal residuals are projected back to the physical space:
\begin{equation}
\label{phi_phys}
\phi_{\sigma,\mathbf{x}}^{K,Rus\star} = \mathbf{R}\,\hat{\phi}_{\sigma, \mathbf{x}}^{K,Rus\star}.
\end{equation}
This guarantees that the scheme is high-order in time and space and (formally) non-oscillatory, see \cite{Ricchiuto2010,Abgrall2006} for more details. 

The local Lax-Friedrich scheme $\phi_{i,\mathbf{x}}^{K,Rus}(\U_h)$ is defined as  
\begin{equation}
\label{Rus}
 \phi_{\sigma,\mathbf{x}}^{K,Rus}(\U_h)=\int_{\partial K}\varphi_\sigma \mathbf{F}(\U_h)\cdot \mathbf{n}\,d\Gamma -\int_K \nabla \varphi_\sigma \cdot \mathbf{F}(\U_h)\,d \mathbf{x} 
     +\alpha(\U_\sigma-\overline{\U_h}).
     \end{equation}
where $\overline\U_h$ is the arithmetic average of all degrees of freedom defining $\U_h$ in $K$.
The viscosity coefficient $\alpha_K$ is connected to the spectral radius 
\begin{equation}\label{eq:rhoS}
\rho_S \equiv \rho_S( \bm{A}(\U) ) = \max ( |\lambda_1|, \ldots , |\lambda_m| ),
\end{equation}
of the normal flux Jacobian matrix
$\bm{A}(\U)=\nabla_{\U} \F (\U)\cdot \bm{n}$  
and reads
\begin{equation}
\alpha_K = N_{DoF}\max_{\sigma \in K} \left( \rho_S(\nabla_{\U} \F (\U)\cdot \nabla \varphi_\sigma \right).
\label{alpha_Mood}
\end{equation}
Here we recall that $N_{DoF}$ corresponds to the number of DoFs in a cell $K$.

In the sequel, this scheme is refered as the Rusanov PSI plus Jump or RPJ scheme.

\item Finally, for $s=0$, i.e. our parachute first order scheme, we consider locally the classical Rusanov scheme recast as \eqref{Rus}, see \cite{RD-positivity-preprint}. 
The Rusanov scheme is also be nicknamed as Rus in the sequel.

\end{enumerate}


\begin{remark}
Having explicitly outlined that the Rusanov-Psi-Jump (RPJ) scheme allows for an $\mathcal{O}(1)$ accuracy across discontinuities might rise the question of the necessity of the parachute scheme. Indeed, the main reason behind this choice is the restrictive applicability in case of arising singularities, due, for example, to pressures close to zero. While in \cite{AbgrallHO2018} the authors have overcome this issue by applying a different formulation of the limiting, which excludes the characteristic projection, we adopt the classical Rusanov scheme locally.
\end{remark}
%

\subsubsection{Preservation of the Conservation}
The main feature of Residual Distribution scheme is that they rely on the formulation \eqref{eq:fluctuation}-\eqref{eq:phiconsBC} which has been shown in \cite{Abgrall2001Mer}  to provide an approximation that converges to the correct weak solutions.
Since we verify locally the conditions of distribution of the residual among a cell and we apply locally on each degree of freedom within a cell $K$ a numerical method, the conservation is guaranteed by construction even if the spatial discretization scheme differs between two neighbouring cells. 

\subsubsection{Positivity Preservation}
The positivity of the  internal energy and the density is guaranteed provided that the parachute scheme guaranties this property.
In the appendix \ref{positivity}, we show that this is true, under a CFL condition, for the Rusanov scheme which is our parachute scheme. More precisely, we show that we can guaranty, even with a Bernstein representation of the data that the internal energy and the density remain positive at the standard Lagrange interpolation points, provided the numerical dissipation is large enough. We also quantify this amount.

\subsection{Detection Procedure}\label{subsec:detection}
The key procedure in an ``a posteriori'' MOOD loop is the detection step.
Given the candidate solution $\U_h^{K,\star,n+1}$ in a cell $K$, the detection procedure
determines if the solution is valid and accepted to be $\U_h^{K,n+1}=\U_h^{K,\star,n+1}$, or unvalid and the solution in 
this bad cell $K$ must be recomputed, i.e. we start again from $\U_h^n$ and apply a more dissipative scheme.
The detection criteria in this work are based on physical/modeling and numerical considerations.
The underlying physics based on the system of PDEs solved must be ensured. For instance
in the case of the hydrodynamics system of equations,
we check for the positivity of the density at each degree of freedom in the 
cell, that is, if the cell $K$ fulfils the 
\begin{itemize}
\item Physical Admissibility Detection criteria
\begin{equation}
 \label{eq:PAD}
\text{PAD}_K= \left\{ 
    \begin{array}{lll}
      1 & \text{if} & \forall \sigma \in K, \; \; \rho_\sigma^{\star,n+1}<0,\\
      0 & \text{else} & 
    \end{array}\right.
\end{equation}
\end{itemize}

The numerical solution in $K$ can not be any undefined or unrepresentable data value
such as Not-A-Number (NaN)\footnote{We define an approximation $x$ to be NaN if $x\text{.EQ.}x$ is .FALSE.} or Infinity (Inf).
In other words we test for 
\begin{itemize}
\item Computational Admissibility Detection criteria
\begin{equation}
 \label{eq:CAD}
 \text{CAD}_K = \left\{ 
    \begin{array}{lll}
      1 & \text{if} & \exists \sigma \in K, \; \; {\U}_\sigma^{\star,n+1}=\text{NaN}  \; \text{ or } \; {\U}_\sigma^{\star,n+1}=\text{Inf},\\
      0 & \text{else} & 
    \end{array}\right.
\end{equation}
\end{itemize}
In case we are within a plateau area, we make sure to not break that area by applying a 
\begin{itemize}
\item Plateau Detection criteria
\begin{equation}
\label{eq:Plateau}
\text{P}_K = \left\{ 
    \begin{array}{lll}
      0 & \text{if} & \exists \sigma \in K,\quad |M_K^n-m_K^n| \geq \mu^3 ,\\
      1 & \text{else} & 
    \end{array}\right.
\end{equation} 
    \end{itemize}
where the neighbourhood $\mathcal{V}(K)$ is the set of cells surrounding $K$, and the relaxed parameter is
given by $\mu=|K|^{1/d}$, with $d$ the size of the considered dimensions within our problem set. The bounds are defined by
\begin{equation}
M_K^n = \max_{K' \in \mathcal{V}(K),  \sigma \in K' } \left( \U_\sigma^{n} \right),  \qquad
m_K^n = \min_{K' \in \mathcal{V}(K),  \sigma \in K' } \left( \U_\sigma^{n} \right). 
\end{equation}
We then test the solution against oscillatory behaviour via a 
\begin{itemize}
\item Numerical Admissibility Detection criteria
\begin{equation} \label{eq:NAD}
\text{NAD}_K = \left\{ 
    \begin{array}{lll}
      1 & \text{if} & \text{  DMP}_K=1 \text{ and SE}_K=1\\
      0 & \text{else.} & 
    \end{array}\right.
\end{equation}
\end{itemize}
This criteria basically allows our solution for an essentially non-oscillatory behaviour and constitutes of two criteria, where, in case  the first one is activated, only the second one will allow to state whether the cell is troubled at all or not. The first one is a so-called

\begin{itemize}

\item[A.] \emph{Relaxed Discrete Maximum Principle (DMP) criteria}
\begin{equation}
\text{DMP}_K = \left\{ 
    \begin{array}{lll}
      0 & \text{if} &  m^n-\epsilon <\U_\sigma^{\star,n+1} < M^n +\epsilon\\
      1 & \text{else} & 
    \end{array}\right.
\end{equation}
\end{itemize}
with \begin{equation}\label{relaxation}
\epsilon=\max\left(\epsilon_1 \left( |M^n - m^n| \right), \epsilon_2\right).
\end{equation}
The value of $\epsilon_1$ is chosen such that the candidate solution could possibly exceed the extrema
but only by a small fraction of the total jump, and will also be subject of analysis and discussion in the numerical section.

In case the cell is marked with DMP$_K=1$, we perform a further check through a 
\begin{itemize}
\item[B.] \emph{Smoothness Extrema Criteria (SE)}, 
in order to exclude the possibility of a mistakenly flagged cell, as, for example in case of natural oscillations with a coarse mesh, this might occur.
\end{itemize}
To this end, we compare in this manuscript two different approaches  for the SE criteria. The first one can be found in the classical MOOD approaches, as for example in \cite{Clain2011,Diot2012}. 
The second considered approach is based on a recent work \cite{Vilar2018}, which has introduced the limiting of \cite{Wang2009,Kuzmin2010} in the MOOD context. This criteria are summarized hereafter.

\subsubsection{A Classical Smoothness Criteria}
Following the idea of \cite{Clain2011,Diot2012,Diot2013}, a more classical Smoothness Extrema criteria (CSE), is generally recast as
\begin{equation}\label{eq:CSE}
\text{CSE}_K = \left\{ 
    \begin{array}{lll}
      1 & \text{if} & \chi_\sigma^{max} \cdot \chi_\sigma^{min}\geq -\mu \quad\text{AND}\quad  \left(\, \big| \frac{\chi_\sigma^{min}}{\chi_\sigma^{max}}\big|< \frac{1}{2}, \quad\text{OR} \quad \max\left(|\chi_\sigma^{max}|,|\chi_\sigma^{min}|\right) \geq\mu\, \right)\\
      0 & \text{else} & 
    \end{array}\right.
\end{equation}  
where we define
$$\chi_\sigma^{min}=\min_{K' \in \mathcal{V}(K),  \sigma \in K' } (\chi_\sigma),\quad \chi_\sigma^{max}=\max_{K' \in \mathcal{V}(K),  \sigma \in K' } (\chi_\sigma)$$ with $\chi_\sigma$  the second order derivative $\chi_\sigma=\nabla^2 \varphi_\sigma \U_\sigma^{\star,n+1}$.
  
\subsubsection{A Linearised Smoothness Criteria}
The second considered approach,  the Linearized Smoothness Criteria (LSE), has been inspired by \cite{Vilar2018}, where the detection criterion has been formulated as a  generalized moment limiter as in \cite{Wang2009} along the hierarchical slope limiting of \cite{Kuzmin2010}.  In particular, our criterion is based on the check of the possible presence of either a discontinuity in the derivative of the approximation between neighbouring elements or the appearance of a change in sign of the derivative between neighbouring cells.  In both cases, the verification would imply the presence of a numerical oscillation, requiring the recomputation with a more dissipative schema. 
This technique can be recast in 1D as  
\begin{equation}\label{eq:LSE}
\text{LSE}_K = \left\{ 
    \begin{array}{lll}
      0 & \text{if} & \min_{K' \in \mathcal{V}(K),  \sigma \in K' } (\hat{\alpha_\sigma}) =1\\
      1 & \text{else} & 
    \end{array}\right.
\end{equation}  
where 

\begin{equation}
\hat{\alpha}_{\text{e}_j}=\left \{ 
    \begin{array}{lll}
    \min\left(1, \dfrac{\partial U_{max,e}^{\star,n+1}- \overline{\partial U_K}^{\star,n+1}}{\widetilde{\partial U_{e_j}}^{\star,n+1}-\overline{\partial U_K}^{\star,n+1}} \right), & \text{if} & \widetilde{\partial U_{e_j}}^{\star,n+1}> \overline{\partial U_K}^{\star,n+1}\\
    & & \\
     \min\left(1, \dfrac{\partial U_{min,e}^{\star,n+1}- \overline{\partial U_K}^{\star,n+1}}{\widetilde{\partial U_{e_j}}^{\star,n+1}-\overline{\partial U_K}^{\star,n+1}}\right), & \text{if} & \widetilde{\partial U_{e_j}}^{\star,n+1}< \overline{\partial U_K}^{\star,n+1}\\
      & & \\
    1 & else &
      \end{array}\right.  
\end{equation}

On $K$ we do have $\partial U_{\sigma}= \int_{e} \nabla \varphi_{\sigma} U_h d\gamma$ and  $\overline{\partial U_K}= \frac{1}{N_{DoF}}\sum_{\sigma \in K} \partial U_{\sigma} $. The same definitions apply for each of the neighbour element $K'$, such that we can define $\partial U_{\sigma'}$ and $\overline{\partial U_{K'}}$.\\
Further, we also define on each element $\partial^2 U_{\sigma}=\int_e \nabla^2 \varphi_{\sigma} U_h \mathbf{n}   d\gamma$, and accordingly $\overline{\partial^2 U_K}= \frac{1}{N_{DoF}}\sum_{\sigma \in K} \partial^2 U_{\sigma} $. Again, the same definition applies to cell $K'$.\\
Then one can define the smoothness across an edge $e_j$ of cell $K$ via $\widetilde{\partial U_{e_j}}= \overline{\partial U_K}+\overline{\partial^2 U_K} \cdot (\mathbf{x}-\mathbf{x}_c)$. 
Here $\mathbf{x}$ denoted the barycentre of the cell $K$, while $\mathbf{x}_c$ to the physical coordinate of the midpoint of the edge being considered.\\
Finally, we define the maximum (minimum, respectively) on an edge $e$ between cell $K$ and cell $K'$ as $\partial U_{max,e}= \max \left( \overline{\partial U_K},  \overline{\partial U_{K'}} \right)$, so that one has the maximum contribution between the current and neighbouring cell.

In 2D, the check between neighboring cells can be done either by comparing the gradients across edges or across vertices. In order to be more consistent with the 1D case, we considered the vertex scenario. In particular, we look for the gradients generated from each neighbour cell sharing the targeted vertex of the element.  Following this, the convex hull of this gradients is compared with the gradient generated by the element under inspection and in case the latter does not result to be belonging within the convex hull,  this would imply that a numerical oscillation is occurring, instead of a true extrema.
More specifically, for the convex hull formulation we have adopted the classical Quickhull algorithm where the points, i.e. the neighboutring gradients were first reordered to create a sort of boundary with a clockwise orientation and the point consisting in the gradient of the element under inspection is verified by taking the normalized cross product to all the other gradients of the neighbouring cells.   The element's gradient  belongs to the convex hull if  the product of the minimum cross-product value with the maximum one are below a defined threshold value, here considered as $10^{-4}$, a value that has shown an expected behaviour, providing robustness on the considered benchmarks problems.

\subsubsection{To Summarize the Detection Criteria}
The cell will be thus flagged as 'good' if its candidate solution fulfils all detection criteria
and it is not a direct neighbour  of a bad cell, otherwise it is flagged as accepted.
Therefore we have \textit{de facto} a set  $\mathcal{B}$ of cells 
to be sent back in time to $t^n$ for re-computation. It is defined by  
\begin{equation}
\label{eq:bad}
\mathcal{B} = \left\{  K\in \Omega_h, \; \text{s.t.} \;
  \left( \text{PAD}_K \times \text{CAD}_K\times\text{P}_K \times\text{NAD}_K = 1 \right) \; \text{or} \; 
  \left( \exists K' \in \mathcal{V}(K), \; K' \in \mathcal{B} \right)  \right\}.
\end{equation} 
\begin{remark}
The iterative MOOD loop is driven by the 'detection procedure' to pull apart
good cells from bad ones. This loop always converges because there is a finite number 
of cells and schemes in the cascade, and the solution provided by the parachute scheme
is always accepted. 
One may observe that the neighbouring cells of a bad one could be destabilized in the next iterate 
as the fluxes in the bad cell are recomputed with a different scheme and contribute 
to the new candidate solution in the previously detected good cell. 
As such one may fear that the correction of one bad cell may spread far away. 
Although there is no mechanism to prevent this behaviour, 
we have not experimented such dramatic phenomena and this is subject to future research. 
In general few percentages of the total number of cells demand a recomputation (see cf. section \ref{sec:numerics}).
Therefore the extra-cost of recomputing the same cells several times is acceptable
\end{remark}
\begin{remark}
Remark that within this relatively non-intrusive 
``a posteriori'' MOOD paradigm, the need of using specific and complex limiting procedure does 
vanish. The stabilization and robustness is gained by the use 
of a preferred low-order scheme ('parachute')
where and when a ('detection') calls it appropriate, 
while the high accuracy is reached on
smooth parts of the flow by the use of one high-order scheme (from the 'cascade').
\end{remark}

\begin{remark}\label{Remark_vars}
In the context of the Euler equation, we have observed that 
whether the detection criteria are applied to all physical variables, or only one, as for example the density, does not play any role at all (see cf. Figure \ref{Fig:SO_WC_vars} of the numerical experiments section), and thus, in order to spare useless computational costs, we advice to adopt the proposed detecting strategy on a single quantity. Of course, it seems almost pointless to remark, that any local treatment of a cell has to be applied to all variables within at the marked area.
\end{remark}



%
%

\section{Numerical Experiments} \label{sec:numerics}

This section introduces and describes a list of representative test cases for the system of PDEs given by the Euler's equations.
Numerical solutions given by the novel strategy are proposed to assess the gain brought by the use of the ``a posteriori'' blending strategy.\\
We shall refer to the second order scheme obtained by using linear shape functions on each element as $\mathcal{B}^1$. Higher order approximations are obtained by choosing quadratic ($\mathcal{B}^2$) or cubic ($\mathcal{B}^3$) Bernstein polynomials as shape functions.
The numerical benchmark problems have been run with the same parameters as in \cite{AbgrallHO2018}, i.e. for the $\mathcal{B}^1$ approximation we consider $M=2$ and $R=2$, for $\mathcal{B}^2$ we input $M=3$ and $R=3$ and, finally, for $\mathcal{B}^3$ we take $M=4$ and $R=4$ in the algorithm presented in Section \ref{ssec.TimeRD}.
All test cases are advanced in time using the Courant-Friedrichs-Lewy condition $\Delta t= \text{CFL}\cdot \min_{\sigma} \big{(} \frac{\Delta x_{\sigma}}{|u_{\sigma}+c_{\sigma}|}\big{)} $, where $\Delta x_{\sigma}$ represents the volume of the cell corresponding to the considered degree of freedom $\sigma$ and $|u_{\sigma}+c_{\sigma}|$ the spectral radius of the solution in $\sigma$.
We have set for all the considered tests a fixed $CFL=0.1$.
The parameters of \eqref{phi_burman} $\theta_1$ and $\theta_2$ depend on the order of accuracy and on the typology of considered system, i.e. they change from 1D to 2D and from the wave equation to the Euler system. In the following considered benchmark problems, we set empirically the values of $\theta_1$ and $\theta_2$ that show a robust stabilization capability. 
In the one-dimensional case, we have set the parameters in \eqref{phi_burman} as follows: $\mathcal{B}^1$ $\theta_1=1$ and $\theta_2=0$; $\mathcal{B}^2$ $\theta_1=1$ and $\theta_2=0$; $\mathcal{B}^3$ $\theta_1=3$ and $\theta_2=10$.
Furthermore,  we shall refer via ``MOOD'' to the novel proposed strategy of this manuscript, which includes the three cycles for $s=0,1,2$ chosen via the ``a posteriori'' detection criteria.
The label ``no MOOD'' will refer to those results obtained with a stabilized blended Rusanov scheme, i.e. we apply the sole scheme given by equation \eqref{RusPsiJump}-\eqref{alpha_Mood} and it will not comprehend the subcell arrangement proposed in \cite{AbgrallHO2018}.  

In the following plots, all the degrees of freedom are being displayed by e.g. solid/dashed lines. Nevertheless, for readability purposes, symbols such as circles are not displayed at each degree of freedom and are meant to help distinguish each approximation. 
\subsection{Numerical Results for 1D Test Cases}
\subsubsection{Convergence Study: Smooth Isentropic Flow}\label{Isoflow_1D}
The first considered test case  is performed to assess the accuracy of our scheme on a smooth isentropic flow problem on an Euler system of equations in one dimension introduced in \cite{ChengShu2014}. 
The initial data for this test problem is the following:
\begin{equation*}
\rho_0(x) = 1 + 0.9999995\sin(\pi x), \quad u_0(x) = 0, \quad p_0(x) = \rho^{\gamma}(x,0),
\end{equation*}
with $x \in [-1,1]$, $\gamma=3$ and periodic boundary conditions. 
The exact density and velocity in this case can be obtained by the method of characteristics and is explicitly given by
\begin{equation*}
\rho(x,t) = \dfrac12\big( \rho_0(x_1) + \rho_0(x_2)\big), \quad u(x,t) = \sqrt{3}\big(\rho(x,t)-\rho_0(x_1) \big),
\end{equation*}
where for each coordinate $x$ and time $t$ the values $x_1$ and $x_2$ are solutions of the non-linear equations
\begin{align*}
& x + \sqrt{3}\rho_0(x_1) t - x_1 = 0, \\
& x - \sqrt{3}\rho_0(x_2) t - x_2 = 0.
\end{align*}

\begin{figure}[H]
\centering
\includegraphics[width=0.42\textwidth]{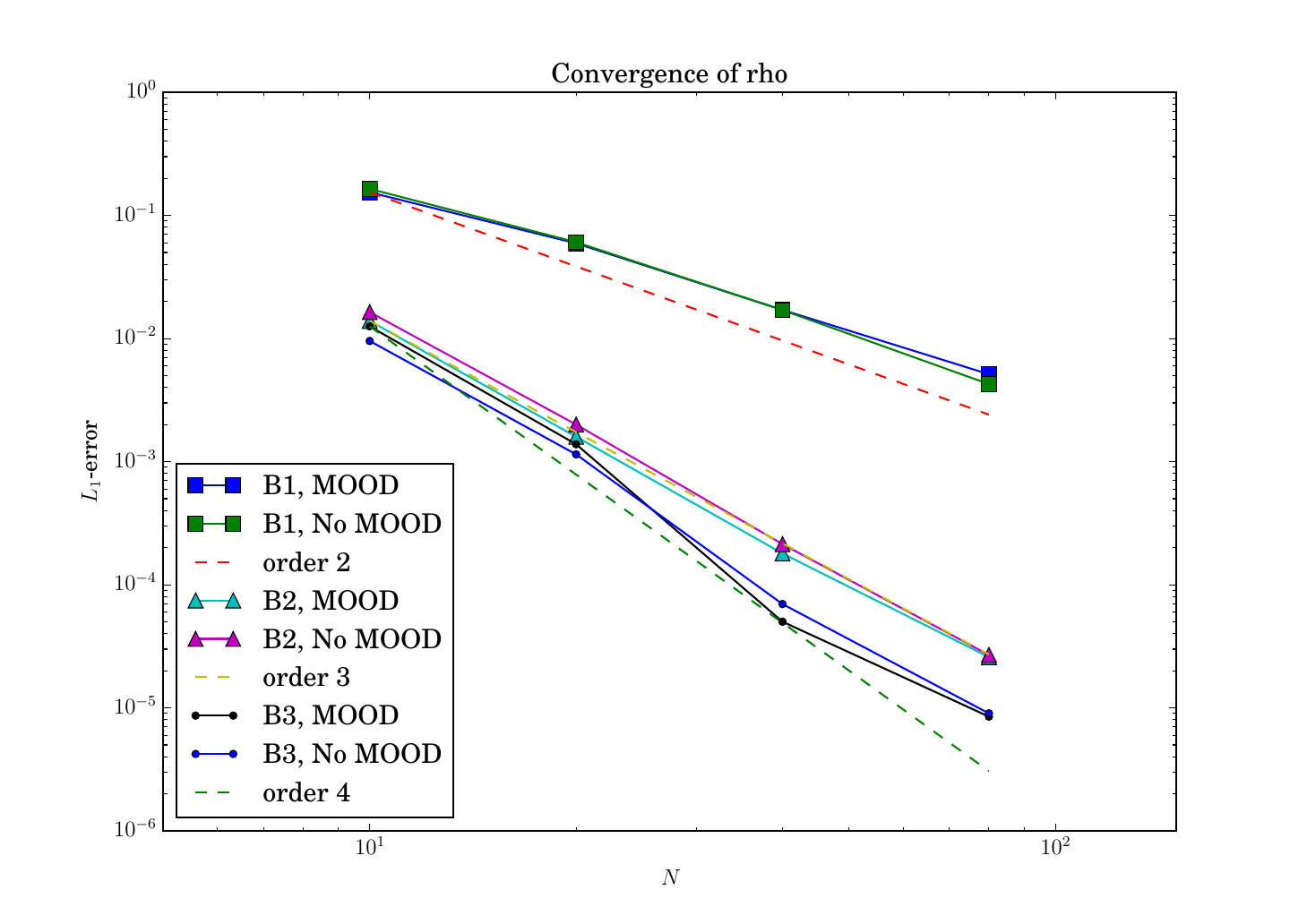}
\includegraphics[width=0.39\textwidth]{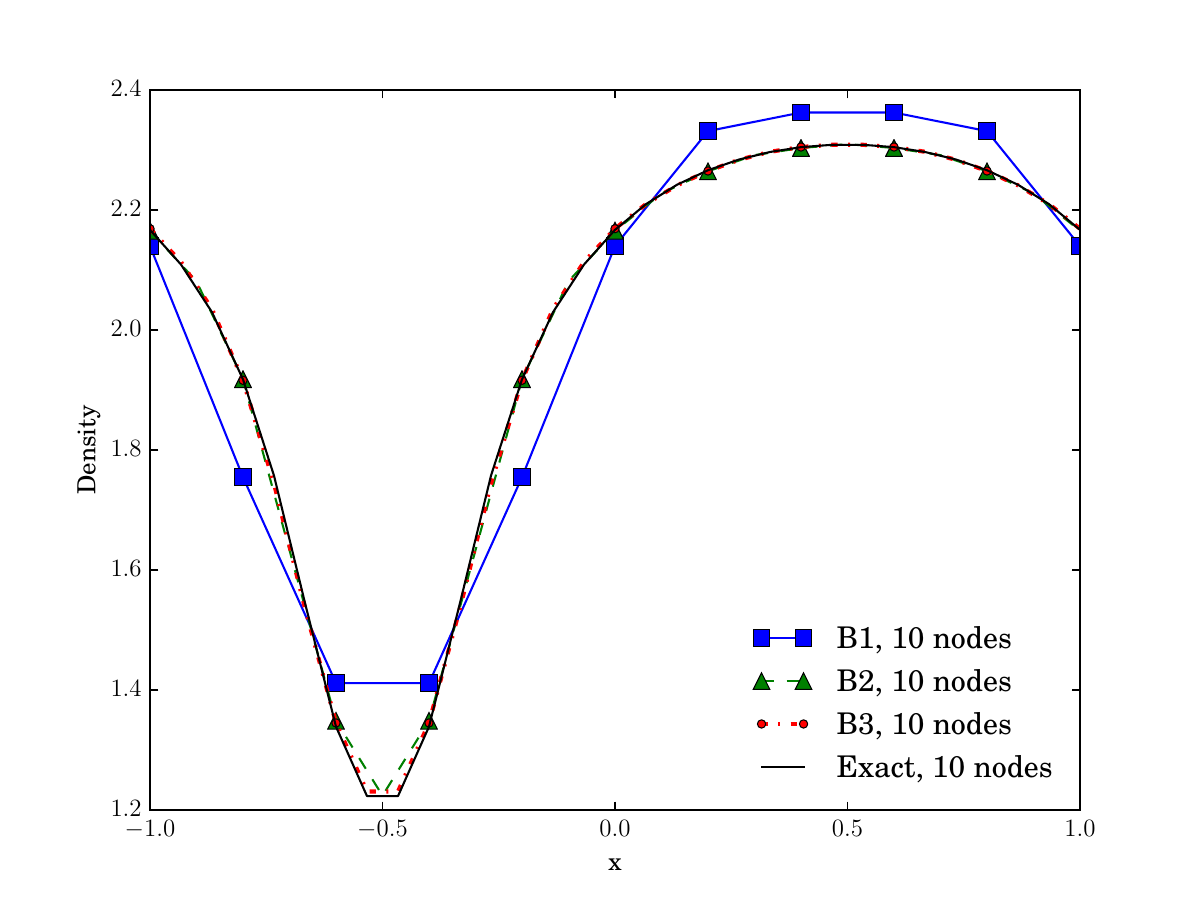}
\caption{Smooth isentropic flow in 1D at $T=0.1$. Left: Convergence plot for the MOOD method (LSE approach) and the one of \cite{AbgrallHO2018} (denoted as 'no MOOD') on the density. Right: Comparison of the density on 10 nodes obtained by the exact solution w.r.t. the approximations obtained by $\mathcal{B}^1, \,\mathcal{B}^2$ and $\mathcal{B}^3$.}
\label{Fig:convergenceorigV} 
\end{figure}

The convergence of the second ($\mathcal{B}^1$), third ($\mathcal{B}^2$) and fourth ($\mathcal{B}^3$) order RD schemes is demonstrated in Figure~\ref{Fig:convergenceorigV},  with the left figure showing the error in the $L_1$ norm for the density vs. the number of cells.
One can observe, that comparing the novel strategy (``MOOD''), against the one obtained via the ``no MOOD'' approach, the sought order of accuracy on a smooth problem is kept.
The accuracy order is throughout guaranteed with the novel strategy: on a smooth problem, none of the detection criteria is activated and the solution is approximated by the $s=s_{max}$ scheme. 
\begin{figure}[H]
\centering
\subfigure[Convergence of the density]{\includegraphics[width=0.45\textwidth]{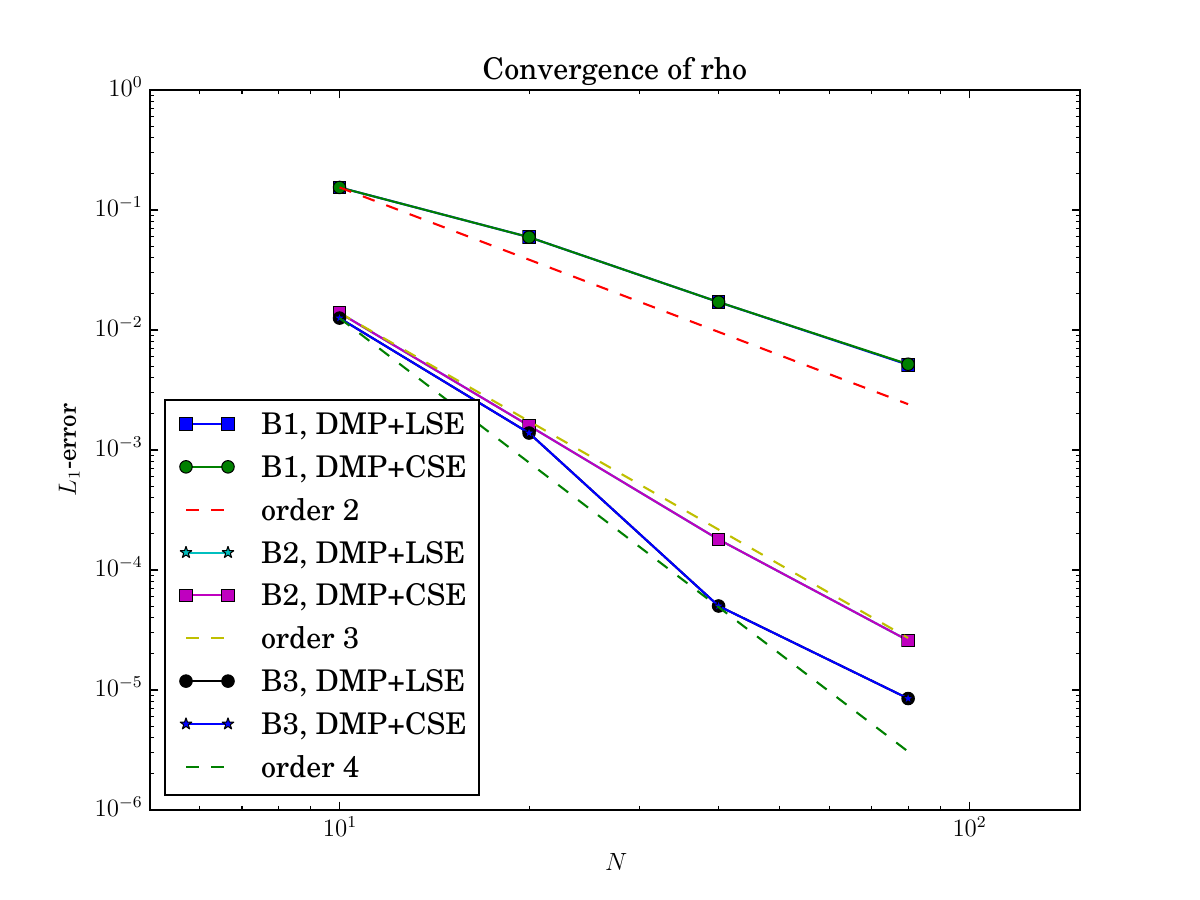}}
\subfigure[Convergence of the pressure]{\includegraphics[width=0.45\textwidth]{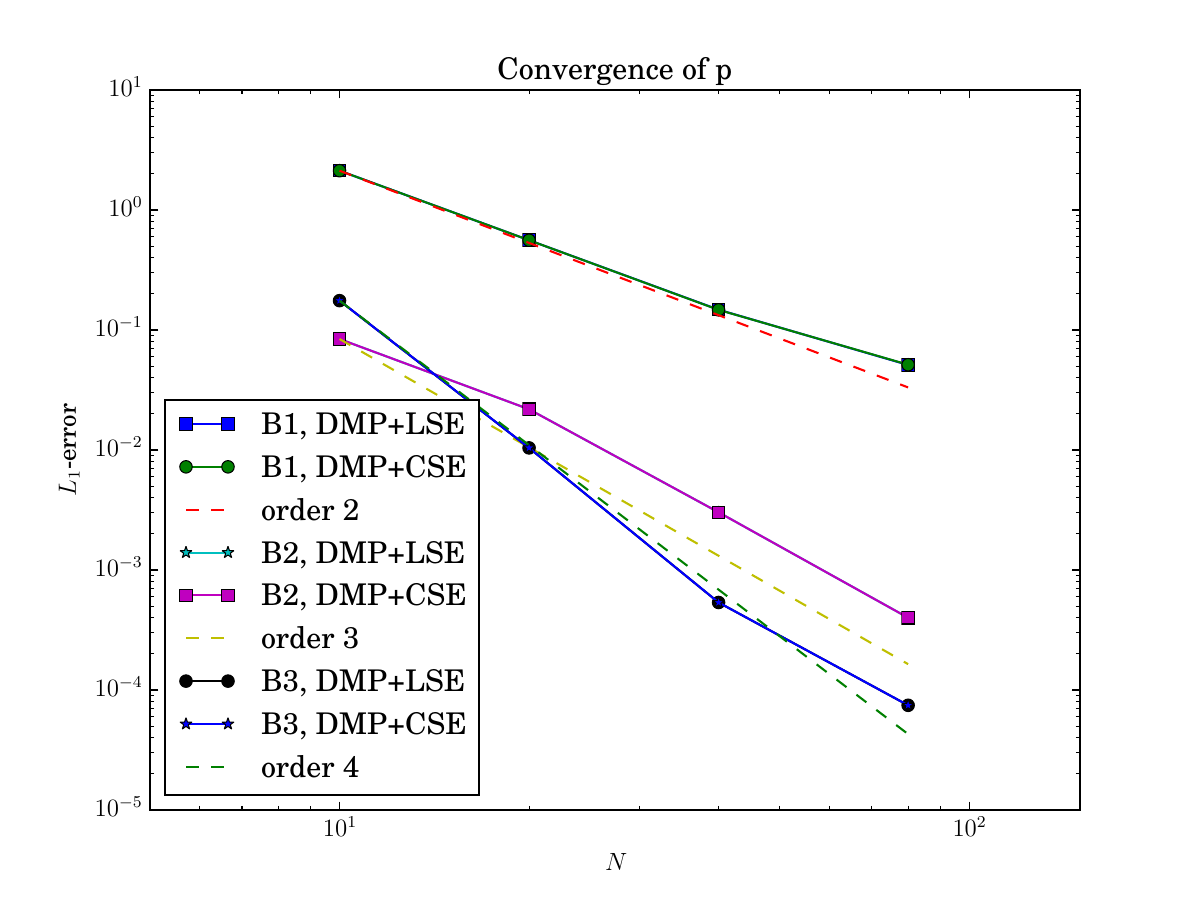}}
\caption{Convergence plot for the MOOD method with two different NAD criteria on a smooth isentropic flow in 1D at $T=0.1$.}
\label{Fig:convergenceV}
\end{figure}
Comparing in Figure \ref{Fig:convergenceV} the convergence for an approximation given by the novel methodology, with the NAD criteria composed by the DMP with $\epsilon_1=\epsilon_2=0$ and the linearised smoothness extrema (LSE) criteria,  together with the one given by the MOOD methodology with the NAD given by a relaxed DMP with $\epsilon_1=10^{-3}$ and $\epsilon_2=0$ and the classical smoothness extrema (CSE) criteria.
We observe an overall good convergence rate for all the variables, which perfectly agrees with the results obtained in \cite{AbgrallHO2018}. This observation is extremly relevant, as it assesses the quality of the introduced ``a posteriori'' limiting within the high order residual distribution method of \cite{AbgrallHO2018}.


\subsubsection{Sod's Shock Tube Problem}

The Sod shock tube is a classical test problem for the assessment of numerical methods for solving the Euler equations. Its solution consists of a left rarefaction, a contact and a right shock wave. The initial data for this problem is given as follows:
\begin{equation*}
(\rho_0, u_0, p_0) = 
\begin{cases}
 (1, 0, 1), \quad &x < 0, \\
 (0.125, 0, 0.1), \quad &x > 0.
\end{cases}
\end{equation*}

We present in the following first two different tests, where we show the solution for different polynomial-orders once keeping the same amount of cells, i.e. 100, and once keeping the same amount of DoFs. Further, we show the how the detection technique works in practice and how the a posteriori technique compares to an a priori one on same and different amount of cells. We then show the convergence towards the exact solution by increasing the amount of cells. Finally, we compare the reduced (RPJ as parachute scheme) and complete cascade (with Rus as parachute) and compare also the two different numerical admissibility detection techniques, i.e. pure DMP, DMP+CSE and DMP+LSE (c.f. Section \ref{subsec:detection}). 
 In case not specified otherwise, the ``a posteriori'' limiting is designed with the detection criteria described in the previous section, where, in particular, the numerical admissible detection (NAD) criteria constitutes of the relaxed discrete maximum principle (DMP) and the linearized smoothness extrema criteria (LSE).\\

\subsubsection*{Comparison between different polynomial-orders with the same amount of cells}
\begin{figure}[H]
\centering
\subfigure[Density]{\includegraphics[width=0.45\textwidth]{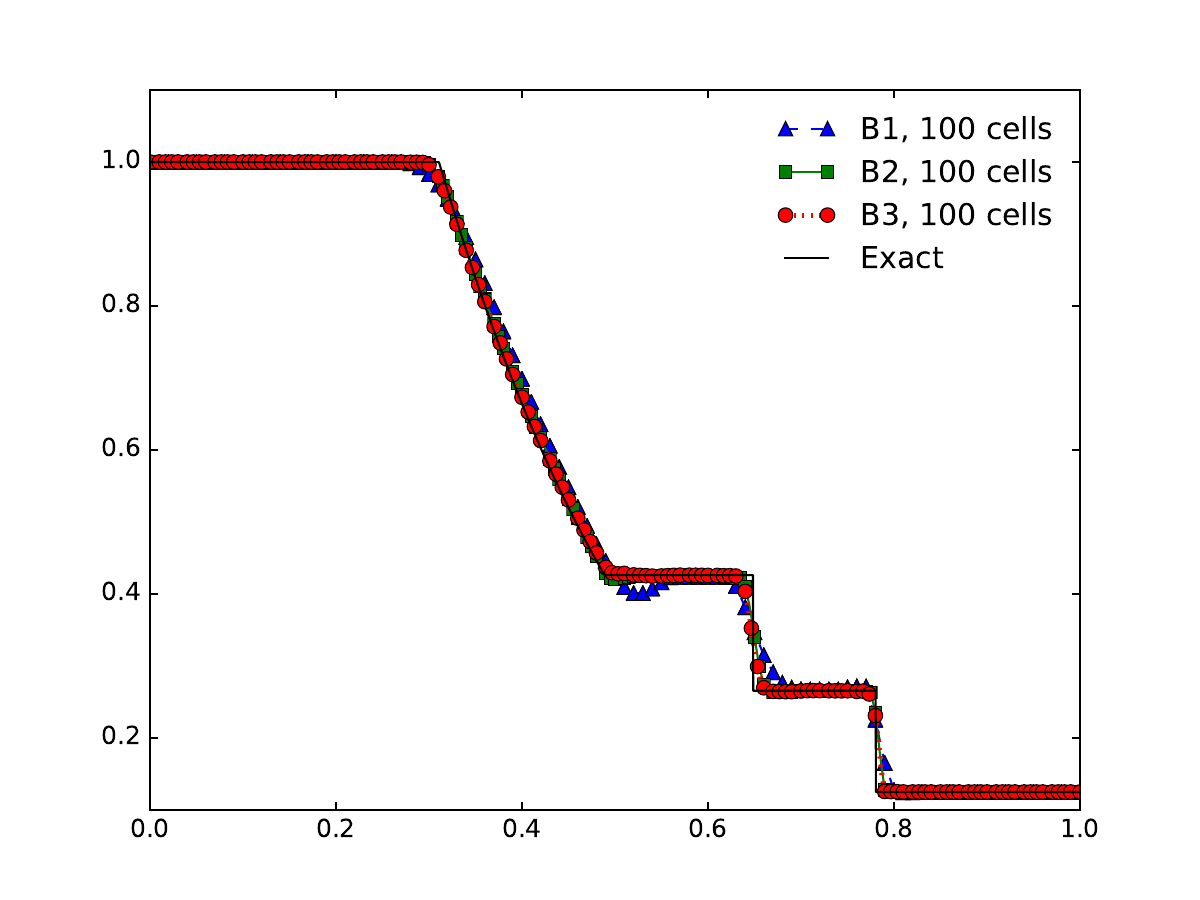}\label{Fig:Sod_density_Bs_100}}
\subfigure[Density, zoom]{\includegraphics[width=0.45\textwidth]{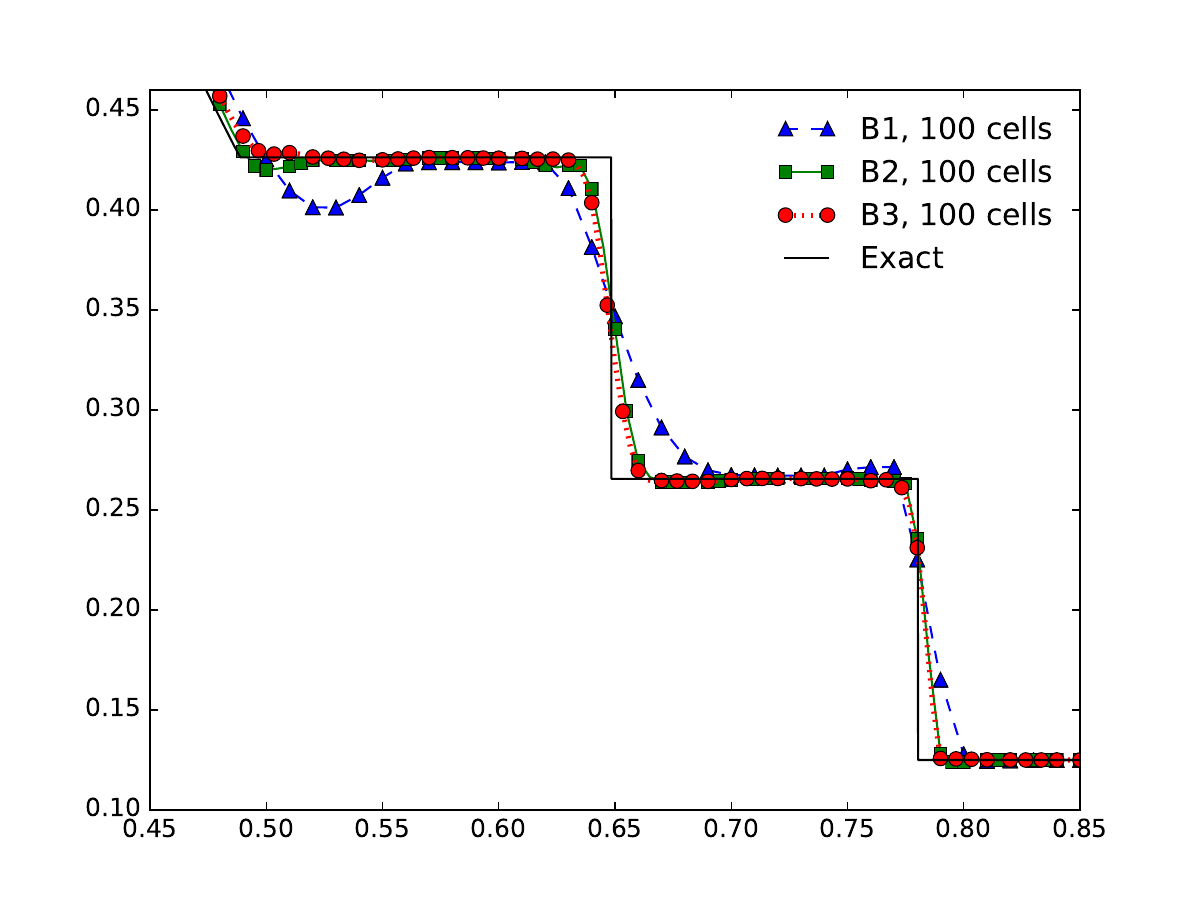}\label{Fig:Sod_density_Bs_100_zoom}}\\
\subfigure[Velocity]{\includegraphics[width=0.45\textwidth]{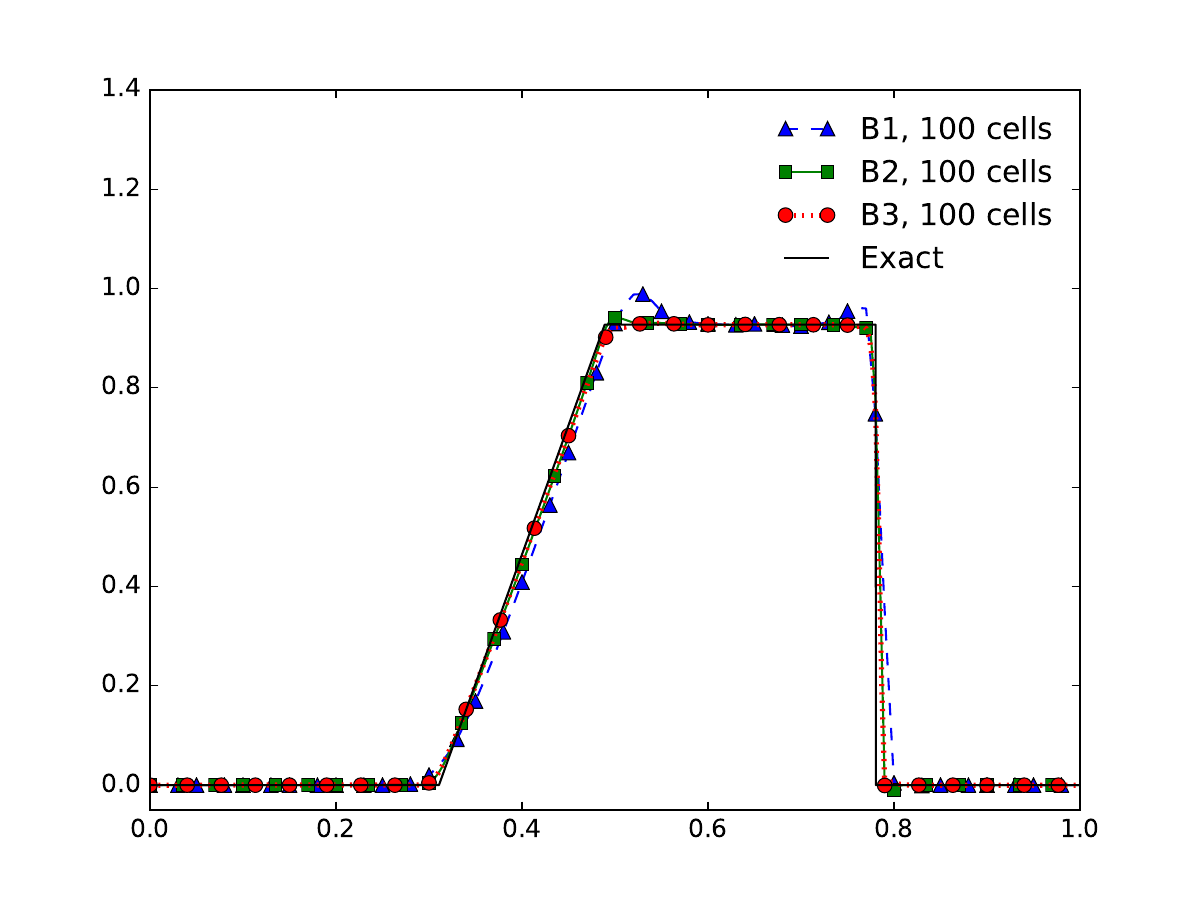}\label{Fig:Sod_velocity_Bs_100}}
\subfigure[Pressure]{\includegraphics[width=0.45\textwidth]{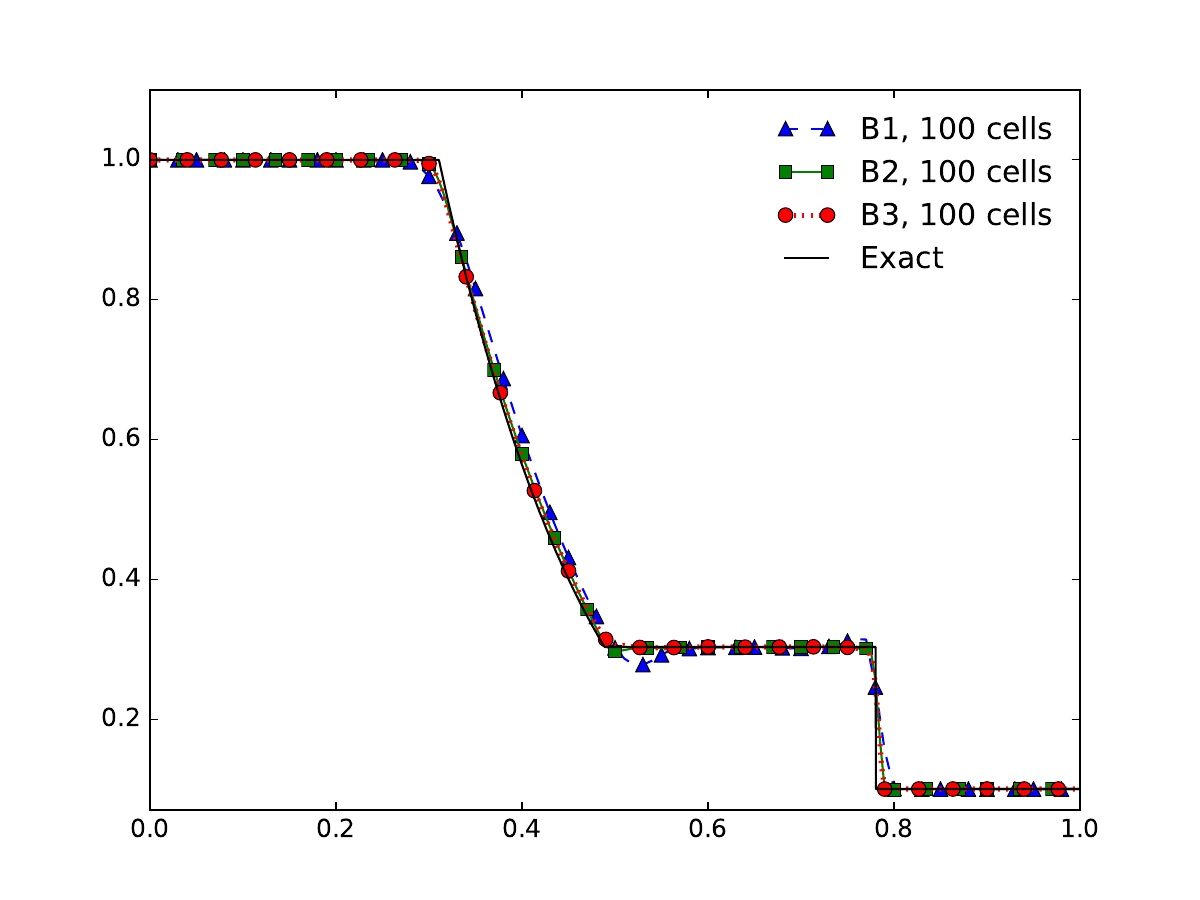}\label{Fig:Sod_pressure_Bs_100}}
\caption{Sod 1D. Density, velocity and pressure for $\mathcal{B}^1$, $\mathcal{B}^2$ and $\mathcal{B}^3$ at $T=0.16$. Note: not all cells are marked for increased readability.}
\label{Fig:Sod_Bs_comparison_100} 
\end{figure}

A first test considers the comparison of $\mathcal{B}^1$, $\mathcal{B}^2$ and $\mathcal{B}^3$ on a mesh of 100 cells with respect to the exact solution. In Figure \ref{Fig:Sod_Bs_comparison_100}, the density of the complete solution (cf. Figure \ref{Fig:Sod_density_Bs_100}) and a zoom between $x=[0.45,0.83]$ (cf. Figure \ref{Fig:Sod_density_Bs_100_zoom}) show how increasing the order, i.e. going from $\mathcal{B}^1$ to $\mathcal{B}^3$, the solution gains in accuracy. 
 
These pictures, together with the one of the velocity (cf. Figure \ref{Fig:Sod_velocity_Bs_100}) and the pressure (cf. Figure \ref{Fig:Sod_pressure_Bs_100}) display under and overshoots of the solution in $\mathcal{B}^1$, as the second order does not allow to provide an accurate approximation with such little nodes. Generally, this is due to the Galerkin scheme with the jump stabilization, as this scheme, which we refer to also with $s=2$ in reference to Figure \ref{fig:mood_cascade}, gives an approximation of the solution across the jump with high numerical oscillations. These are mostly damped by the activated detection criteria which locally treat those oscillations with more dissipating numerical schemes smoothing  the solution. 
One can see that, increasing the order,  i.e. testing from $\mathcal{B}^1$ to $\mathcal{B}^3$, these oscillations are less evident and, thus, the approximation gains in quality.

\subsubsection*{Comparison between different polynomial-orders with the same amount of DoFs}
\begin{figure}[H]
\centering
\subfigure[Density]{\includegraphics[width=0.45\textwidth]{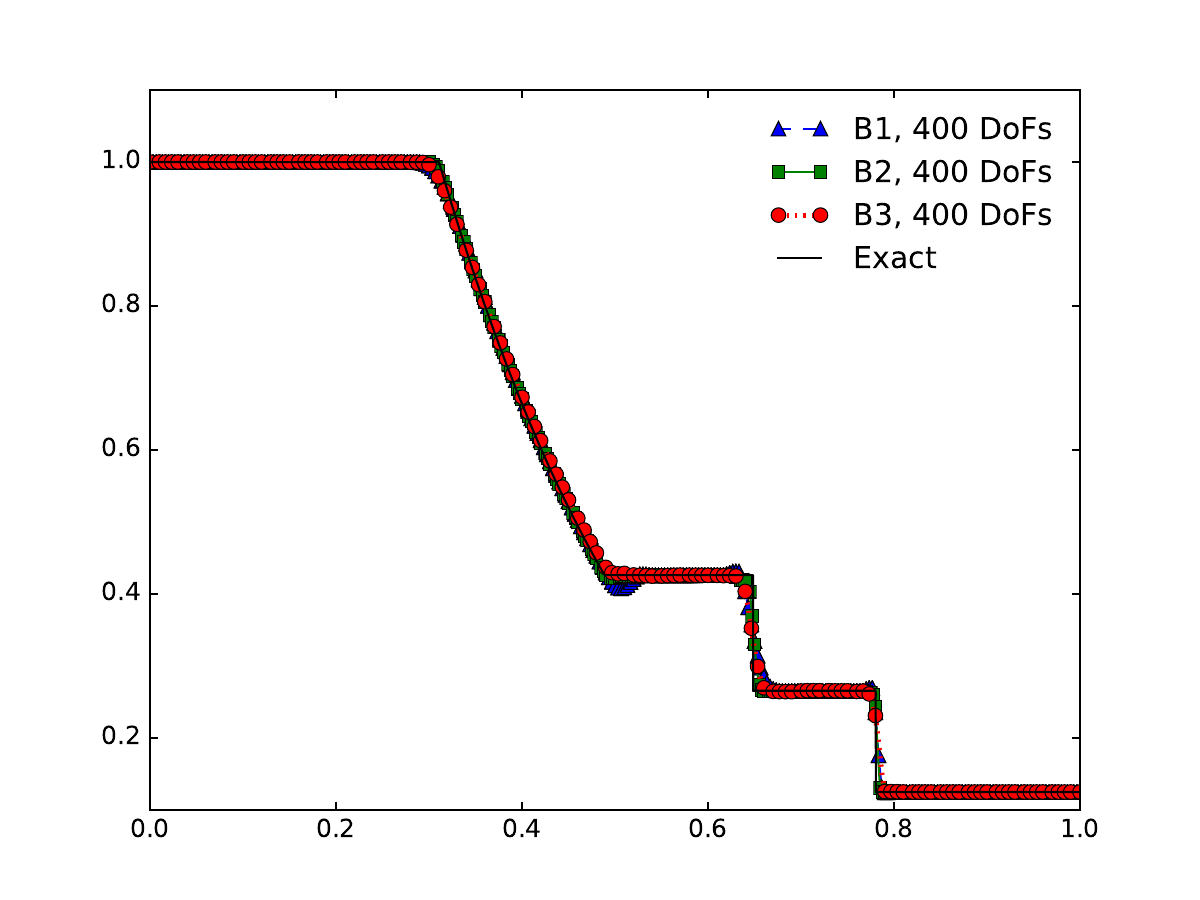}}
\subfigure[Density, zoom]{\includegraphics[width=0.45\textwidth]{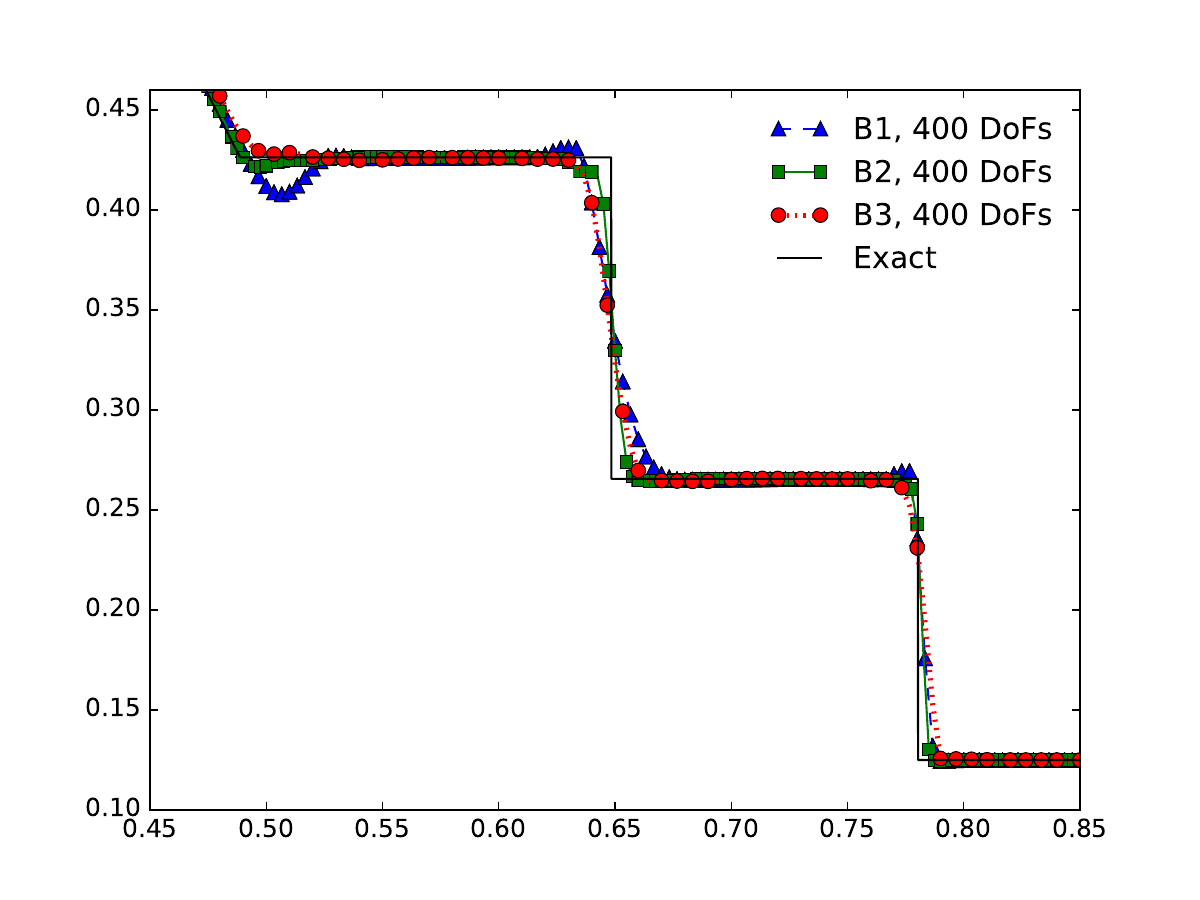}}\\
\caption{Sod 1D. Density for $\mathcal{B}^1$, $\mathcal{B}^2$ and $\mathcal{B}^3$ with the same amount of DoFs at $T=0.16$. Note: not all cells are marked for increased readability.}
\label{Fig:Sod_Bs_comparison_200_dof} 
\end{figure}

In the following study, we have considered $\mathcal{B}^1$, $\mathcal{B}^2$ and $\mathcal{B}^3$, each with a total of 400 degrees of freedom on the whole domain. The NAD criteria considers here the DMP+LSE check. Figure \ref{Fig:Sod_Bs_comparison_200_dof} shows once  more, after the convergence study of a previous section, that there is no gain at all as in $\mathcal{B}^3$ enhancing the amount of cells for $\mathcal{B}^1$ to get the same amount of DoFs. 

\subsubsection*{Detecting Technique in Practice}
\begin{figure}[H]
\centering
\subfigure[$1$st iteration, spatial scheme indicator]{\includegraphics[width=0.45\textwidth]{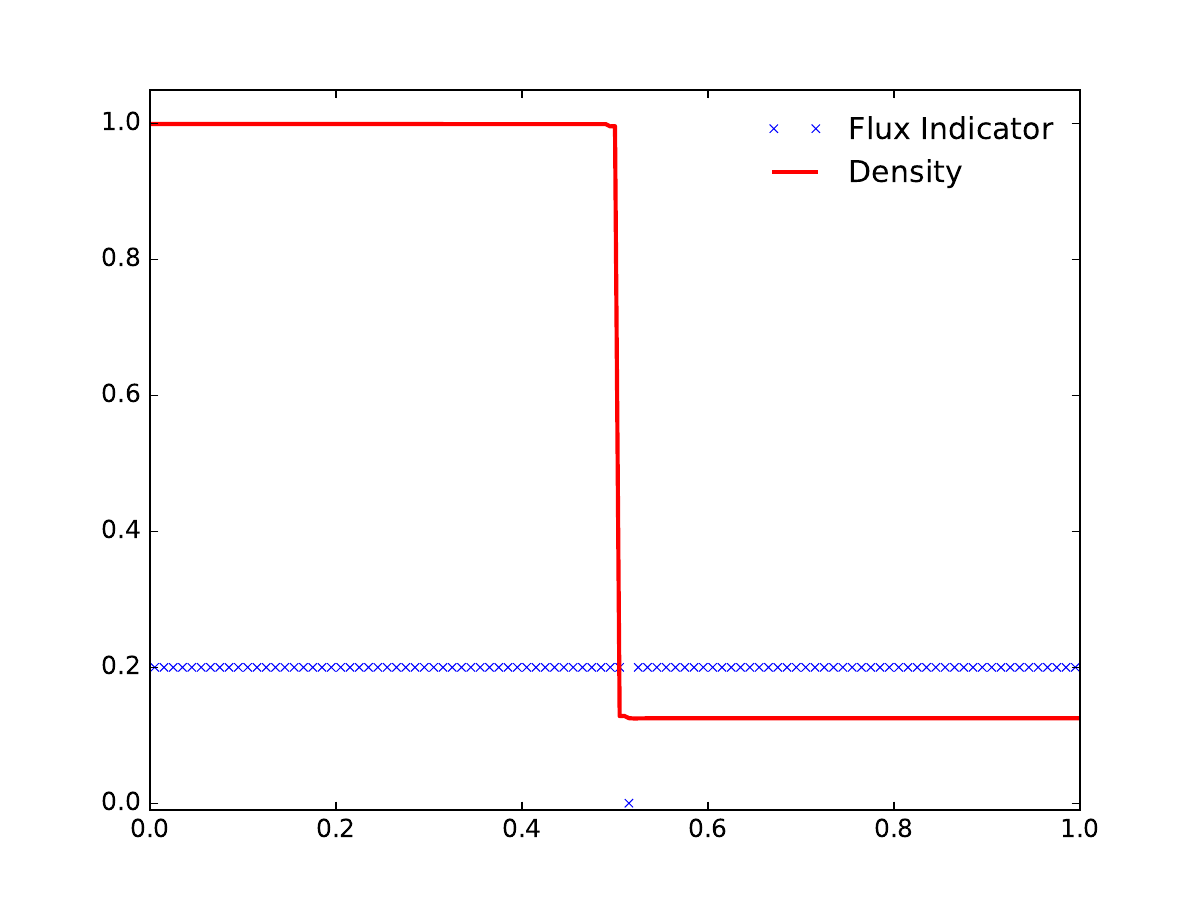}\label{Fig:Sod_scheme_1}}
\subfigure[$50$th iteration, spatial scheme indicator]{\includegraphics[width=0.45\textwidth]{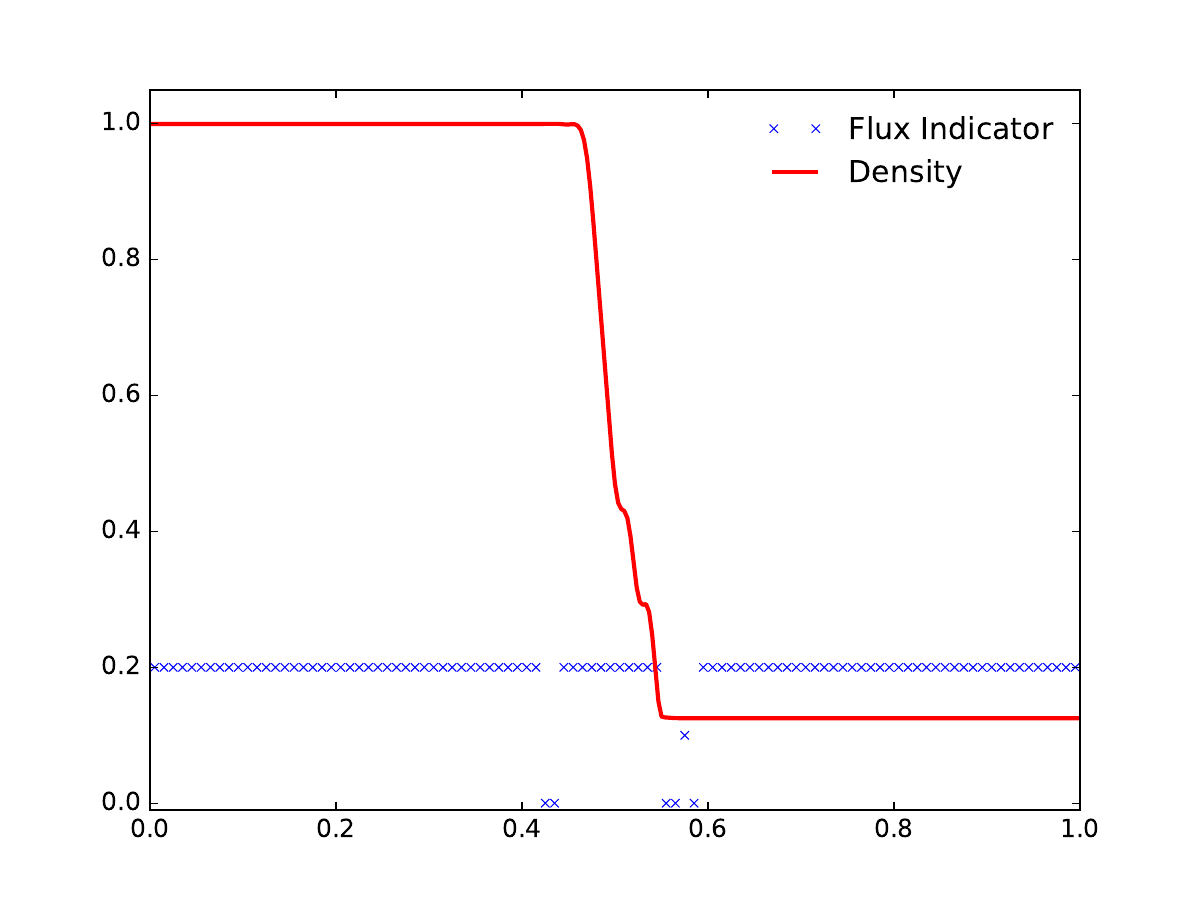}\label{Fig:Sod_scheme_10}}\\
\subfigure[$1$st iteration, detection criteria indicator, $s=1$]{\includegraphics[width=0.45\textwidth]{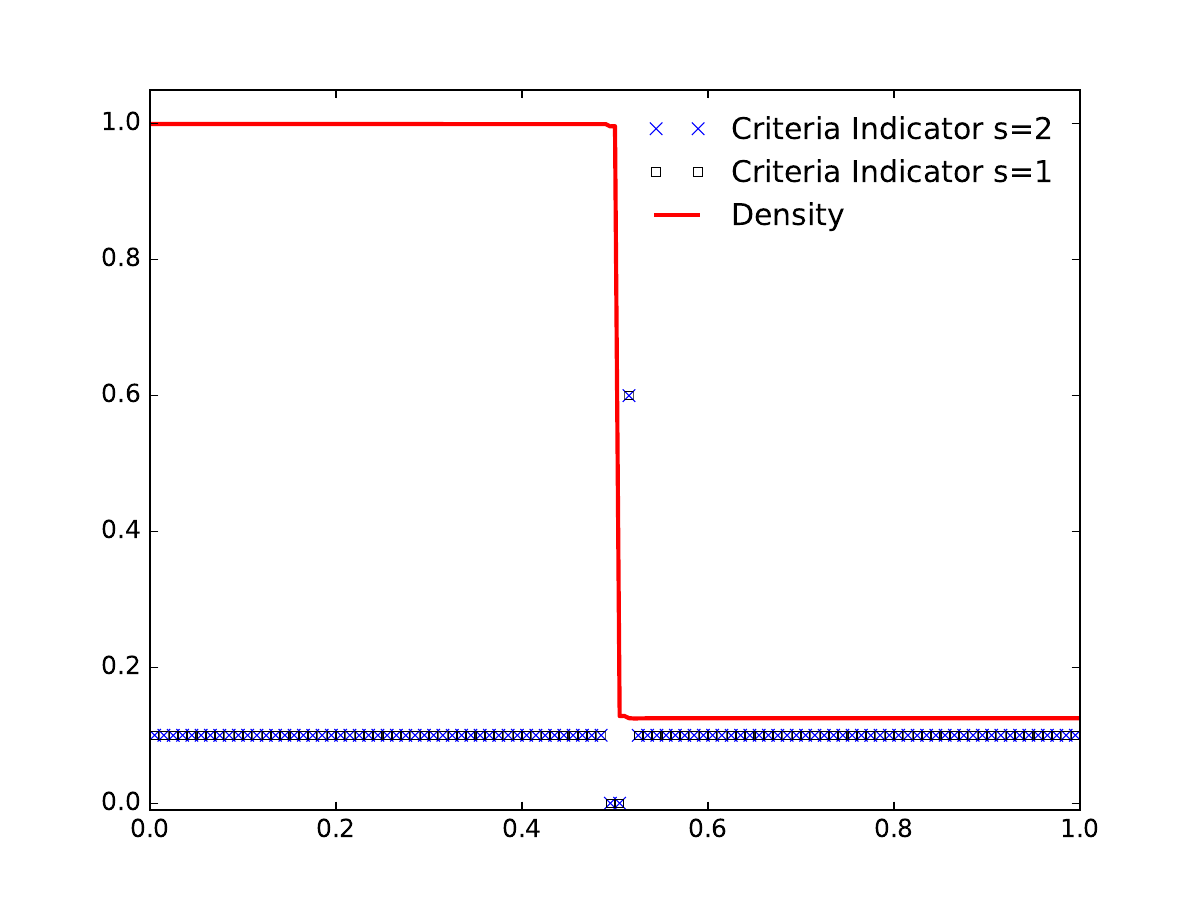}\label{Fig:Sod_detect_1_2}}
\subfigure[$50$th iteration, detection criteria indicator, $s=1$]{\includegraphics[width=0.45\textwidth]{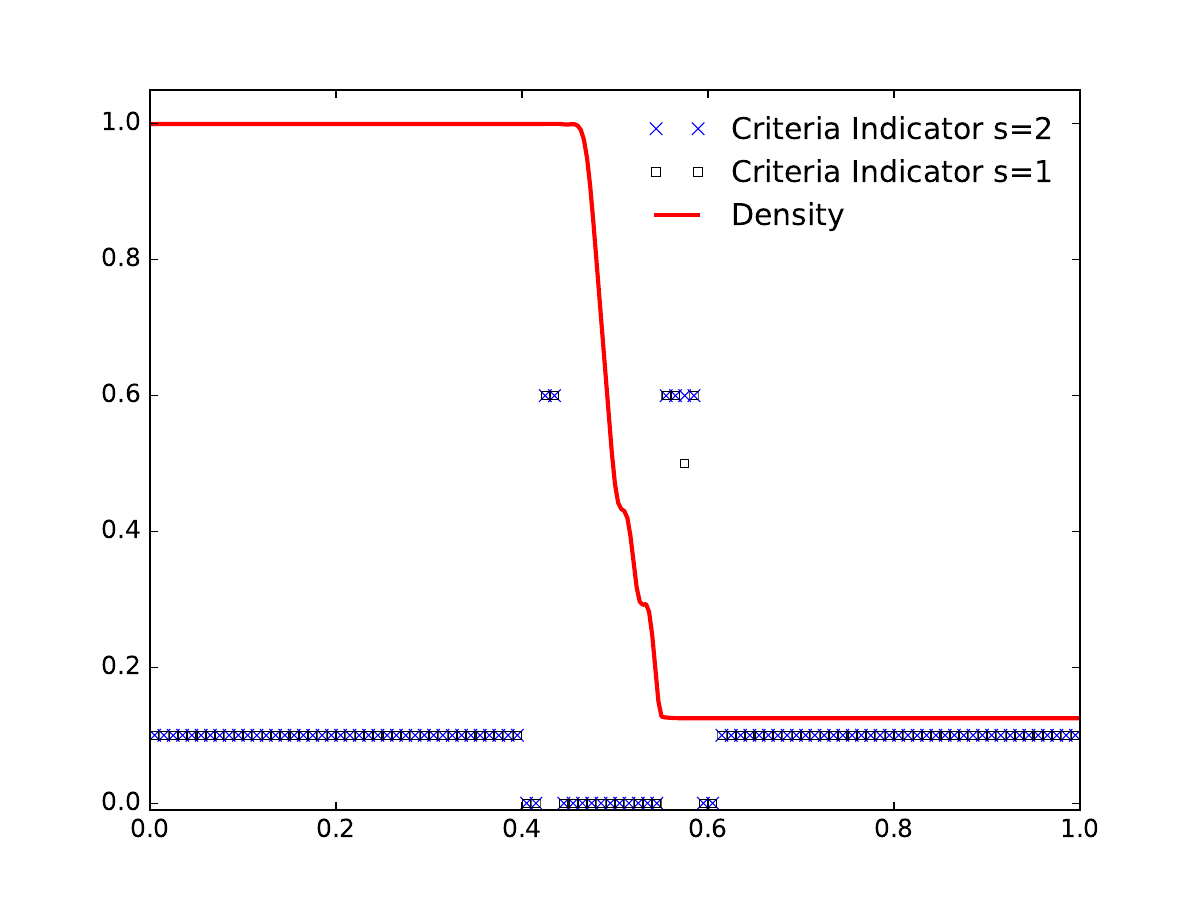}\label{Fig:Sod_detect_10_2}}
\caption{Sod 1D. Detection Criteria activation and considered spatial scheme for the $1$st (left) and $50$th (right) iteration in time on $\mathcal{B}^3$.}
\label{Fig:Sod_detection} 
\end{figure}

To better understand how the detection technique works, we show in Figure \ref{Fig:Sod_detection} the very first iteration in time, i.e. $n=1$ and the $n=50$-th iteration in time. Figures \ref{Fig:Sod_scheme_1} and \ref{Fig:Sod_scheme_10} show the approximated solution for each of the considered iterations along with a flag that indicates which typology of scheme is applied. Here, when the flag is at $0.2$ on the vertical axis, we adopt the standard Galerkin with stabilizing jump terms, i.e. $s=2$ or GPJ, as spatial approximation.
Further, this indicator is at $0.6$ for the Rusanov with PSI limiting and stabilizing jump terms ($s=1$ or RPJ) and at $0$ for the Rusanov  scheme. One can note how in the very first iteration only nodes at the shock interface are being treated by the most dissipative scheme, while, at the $50$th iteration, the Rusanov scheme is applied left and right of the interface, where one would normally presume oscillations developing for non damped schemes. Only one cell is locally flagged with the RPJ scheme.
Moreover, to see what is actually happening behind this decision of which flux is activated, we show in Figures \ref{Fig:Sod_detect_1_2} and \ref{Fig:Sod_detect_10_2} the detection criteria for the first cycle via crosses and second cycle via squares, i.e. respectively when passing from GPJ (s=2) to RPJ (s=1) and from RPJ (s=1) to Rusanov (s=0) in the cascade (see Figure \ref{fig:mood_cascade}). Here the $0.1$ vertical axis values of the indicator corresponds to the plateau detection, whereas $0.6$ corresponds to NAD detection, $0.5$ would correspond to the sole relaxed DMP activation without the LSE, allowing thus the flux to be kept as from the candidate solution. The values of $0.2$ would correspond to the CAD criteria, $0.3$ to the PAD and $0$ is generally signalling that none of the detections has been considered. 

\subsubsection*{The Proposed Method Compared to an ``A Priori Technique''}

\begin{figure}[H]
\centering
\subfigure[$\mathcal{B}^1$]{\includegraphics[width=0.45\textwidth]{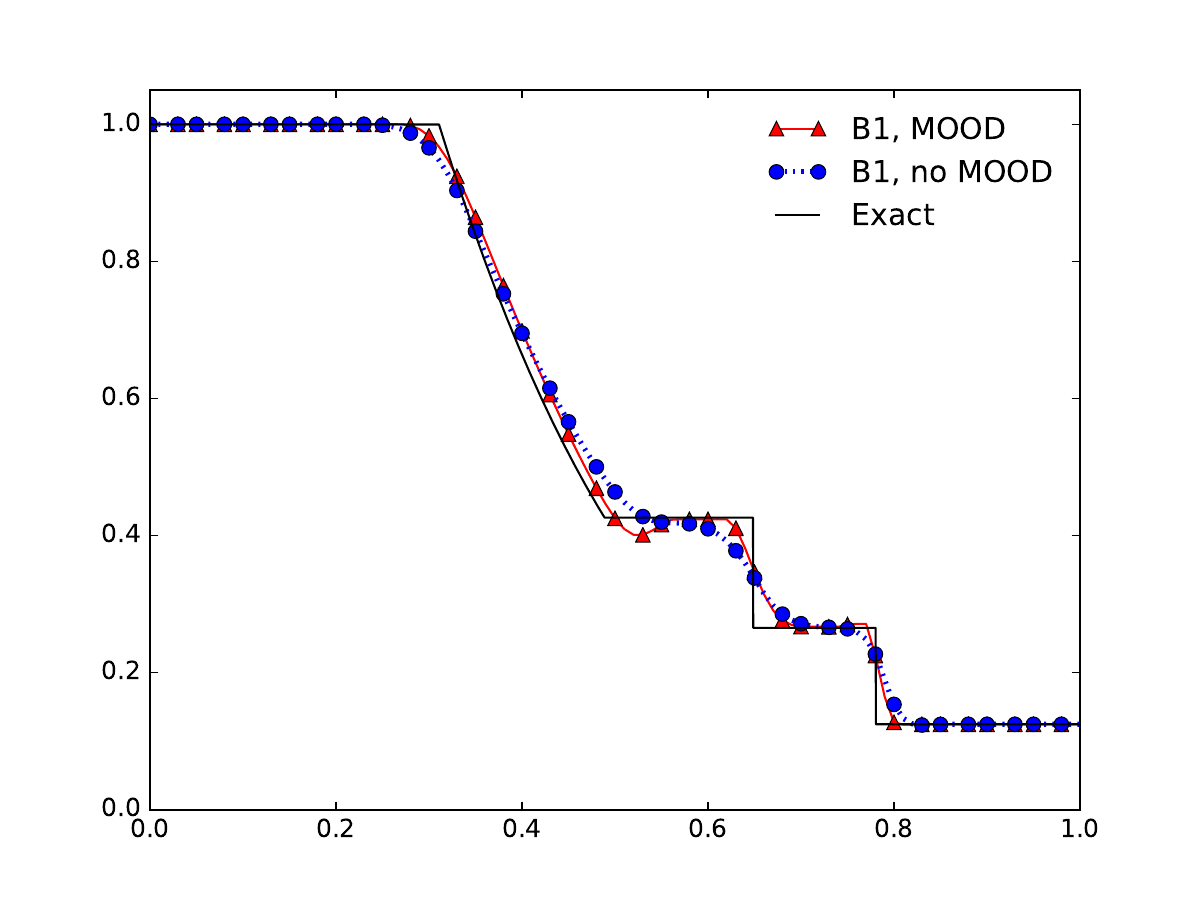}}
\subfigure[$\mathcal{B}^2$]{\includegraphics[width=0.45\textwidth]{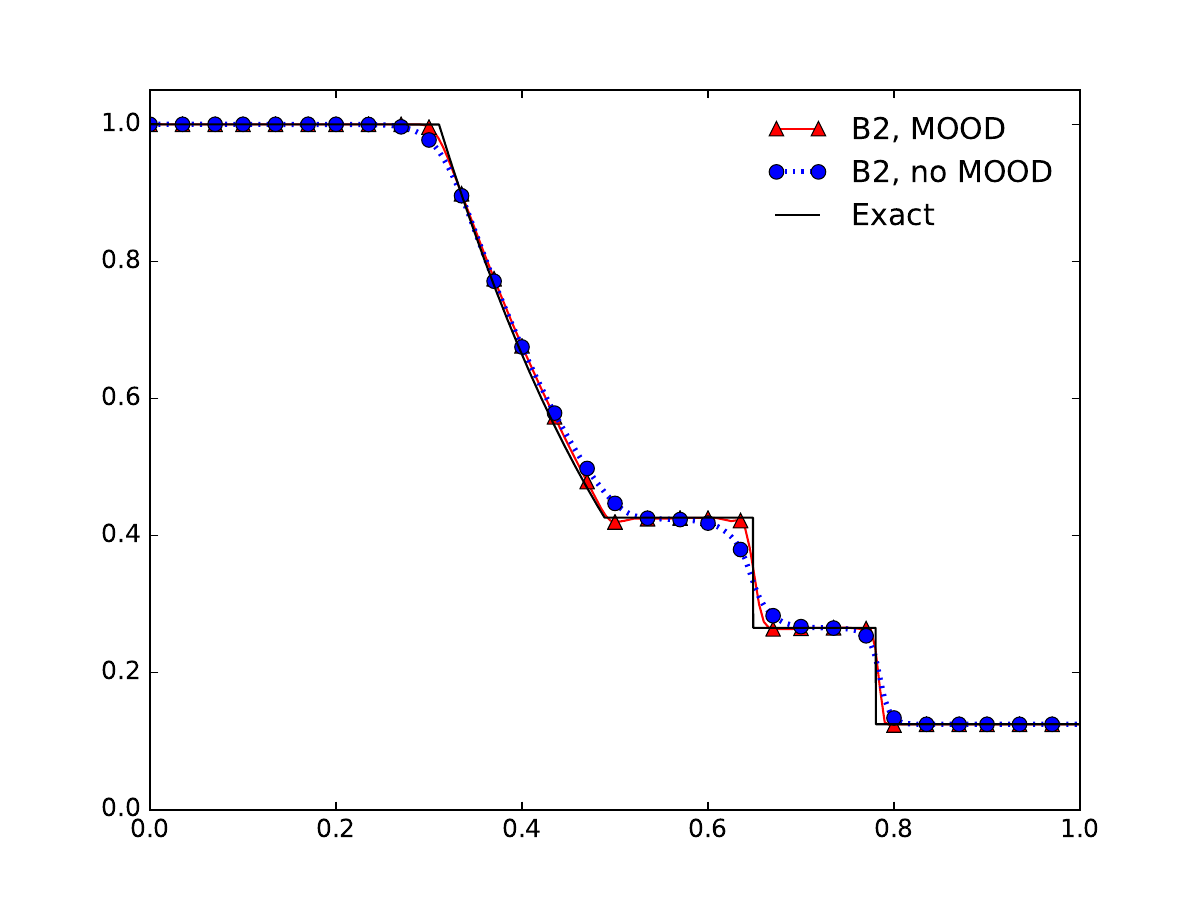}} 
\caption{Sod 1D. Comparison of the densities between the MOOD and no-MOOD schemes on $100$ cells at $T=0.16$.}
\label{Fig:Sod_MoodVSnoMood_Bs} 
\end{figure}
To validate the presented methodology, we consider in Figure \ref{Fig:Sod_MoodVSnoMood_Bs} the comparison between our novel approach and the one of \cite{AbgrallHO2018} where we have a pure RPJ scheme. 
In the following, by referring to the pure RPJ scheme we are basically considering the Rusanov PSI jump scheme without any MOOD and without the sub-cell implementation (cf. section 3.5 in \cite{AbgrallHO2018}). This juxtaposition allows to see the advantage which might bring the ``a posteriori limiting'', taking the same quantity of cells, as here for example $100$.
\subsubsection*{Less DoFs for a Comparable Quality}
\begin{figure}[H]
\centering
\subfigure[Density, $\mathcal{B}^3$]{\includegraphics[width=0.45\textwidth]{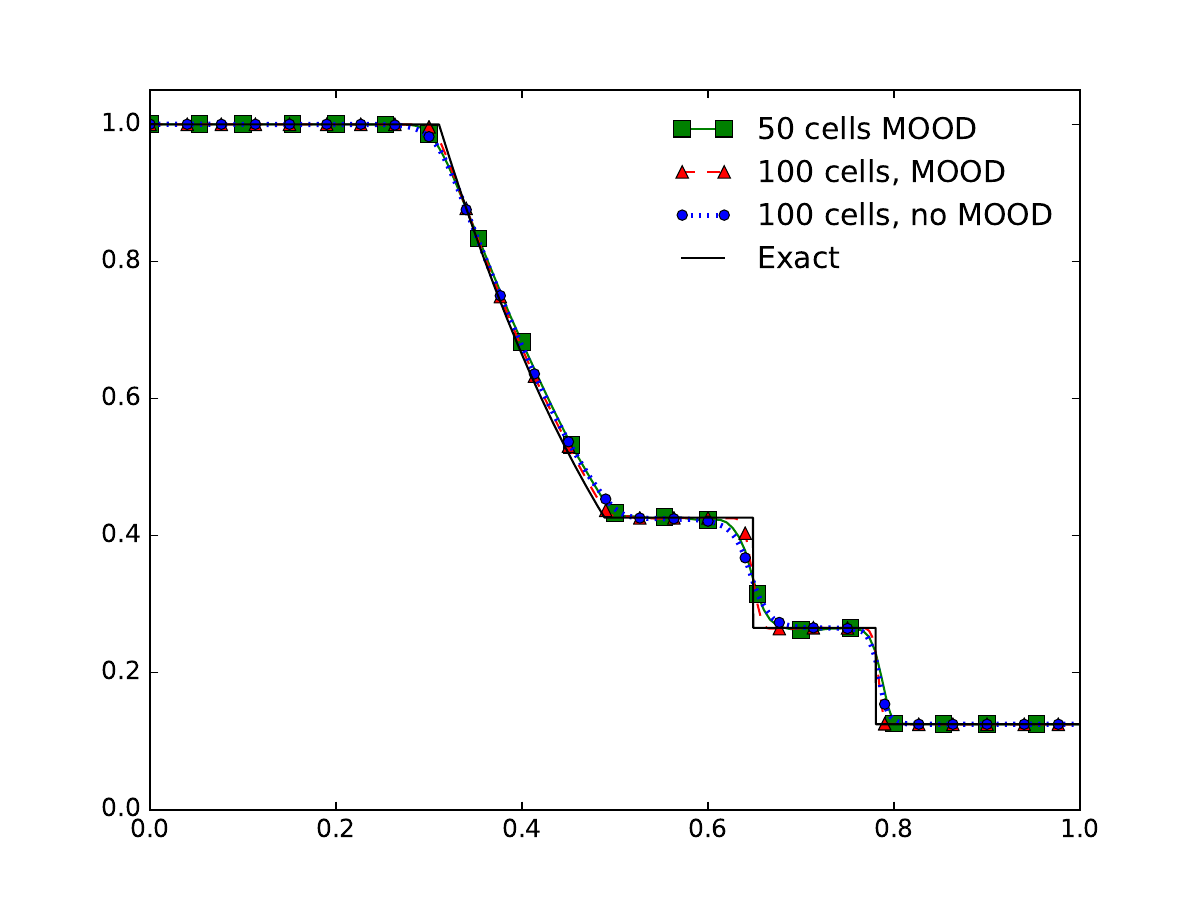}}
\subfigure[Density, $\mathcal{B}^3$, zoom]{\includegraphics[width=0.45\textwidth]{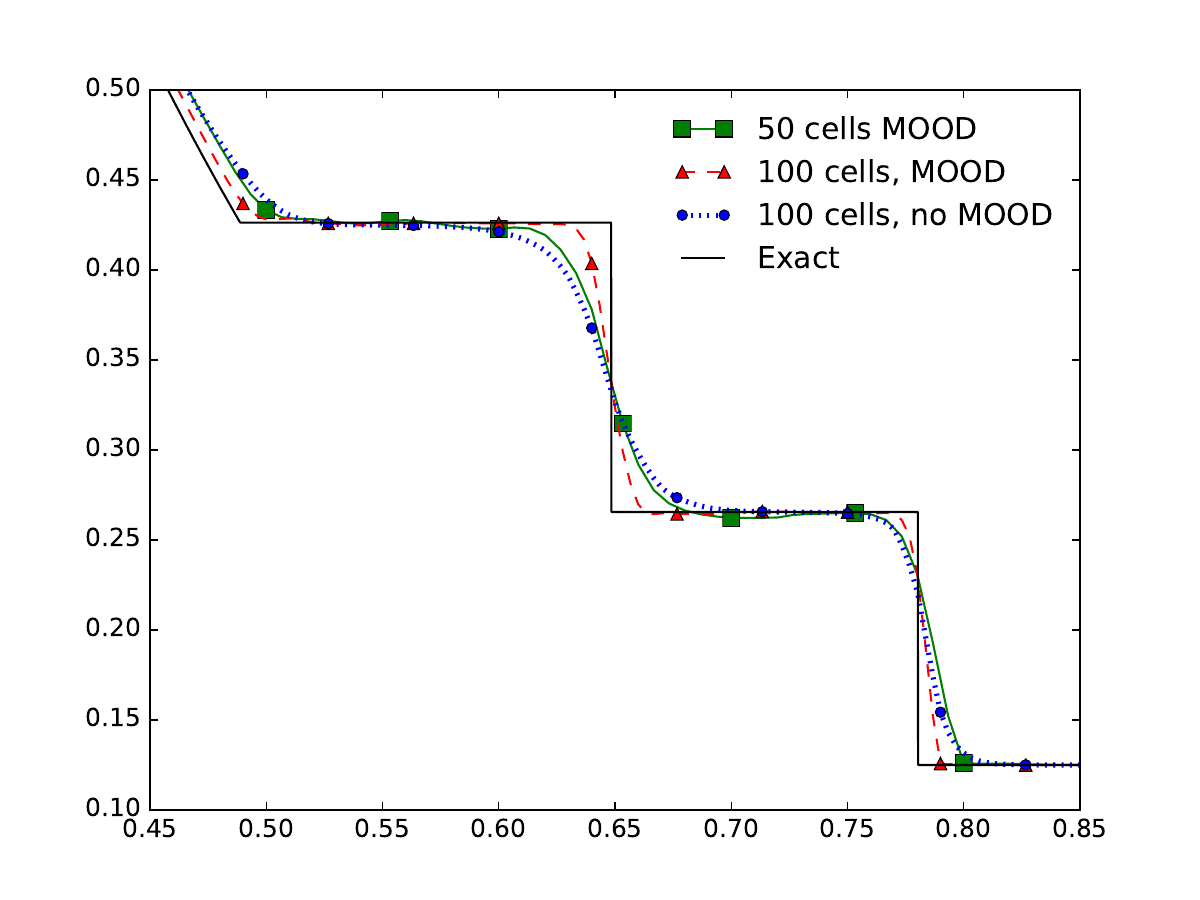}}
\caption{Sod 1D. Comparison of the densities with zoom for different mesh sizes for the MOOD and non-MOOD schemes at $T=0.16$}
\label{Fig:Sod_MoodVSnoMood_meshes} 
\end{figure}

To disclaim eventual remarks that it might be true that the ``a posteriori'' technique is more accurate w.r.t. an ``a priori'', but at some efficiency cost in terms of velocitywe have compared in Figure \ref{Fig:Sod_MoodVSnoMood_meshes}
 the approximations obtained by the ``a posteriori'' technique on $50$ and $100$ cells and compared it to the pure RPJ scheme without MOOD on $100$ cells. The resulting observation is extremely interesting as we can obtain a more detailed result with the coarser mesh and the proposed method compared to the extremely dissipative solution of the pure GPJ without the novel limiting strategy.

\subsubsection*{Convergence Towards the Exact Solution}
\begin{figure}[H]
\centering
\subfigure[Density, $\mathcal{B}^3$]{\includegraphics[width=0.45\textwidth]{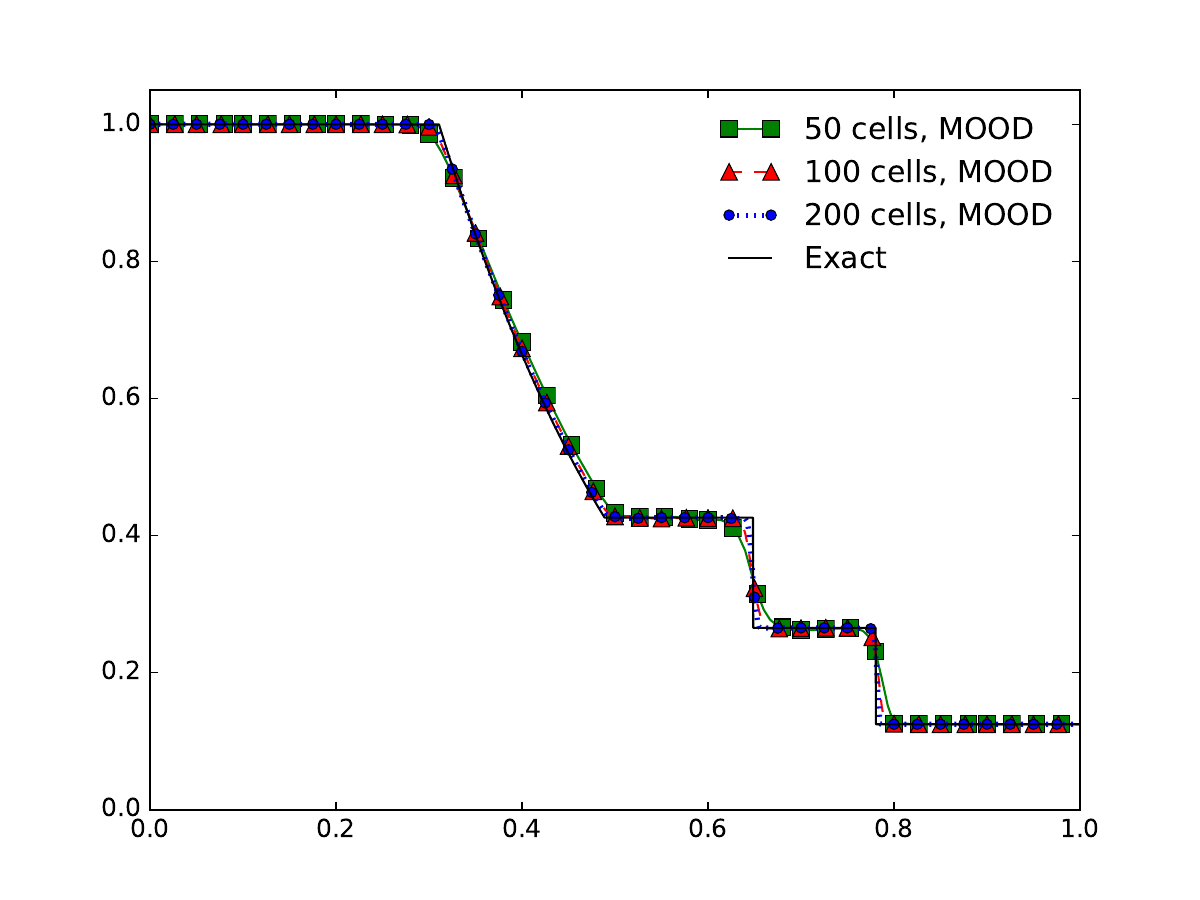}}
\subfigure[Density, $\mathcal{B}^3$, zoom]{\includegraphics[width=0.45\textwidth]{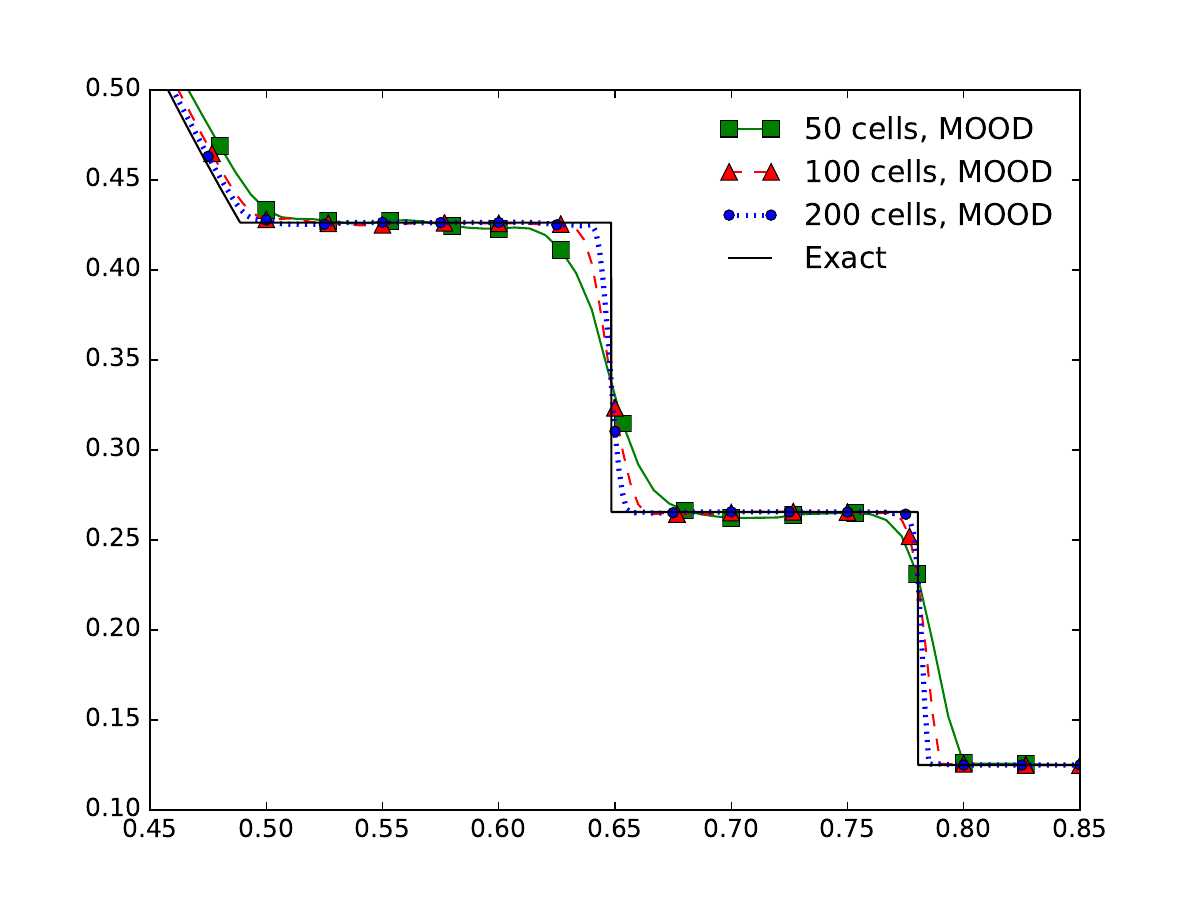}}
\caption{Sod 1D. Mesh convergence for $\mathcal{B}^3$ at $T=0.16$.}
\label{Fig:Sod_convergence} 
\end{figure}

To be able to see whether we do converge to the exact solution on shocking flows, Figure \ref{Fig:Sod_convergence} compares the approximation of our new method for $\mathcal{B}^3$ on $50$, $100$ and $200$ cells. The results show perfect agreement with the expected behaviour.

\subsubsection*{Differences for a Complete vs. Reduced Cascade}	\label{ComplRed_Sod}
 One natural question might arise, whether there is a difference in the solution by applying just GJ and RPJ or the three schemes including the Rusanov, in terms of a cascade.
\begin{figure}[H]
\centering
\subfigure[Density, complete cascade]{\includegraphics[width=0.45\textwidth]{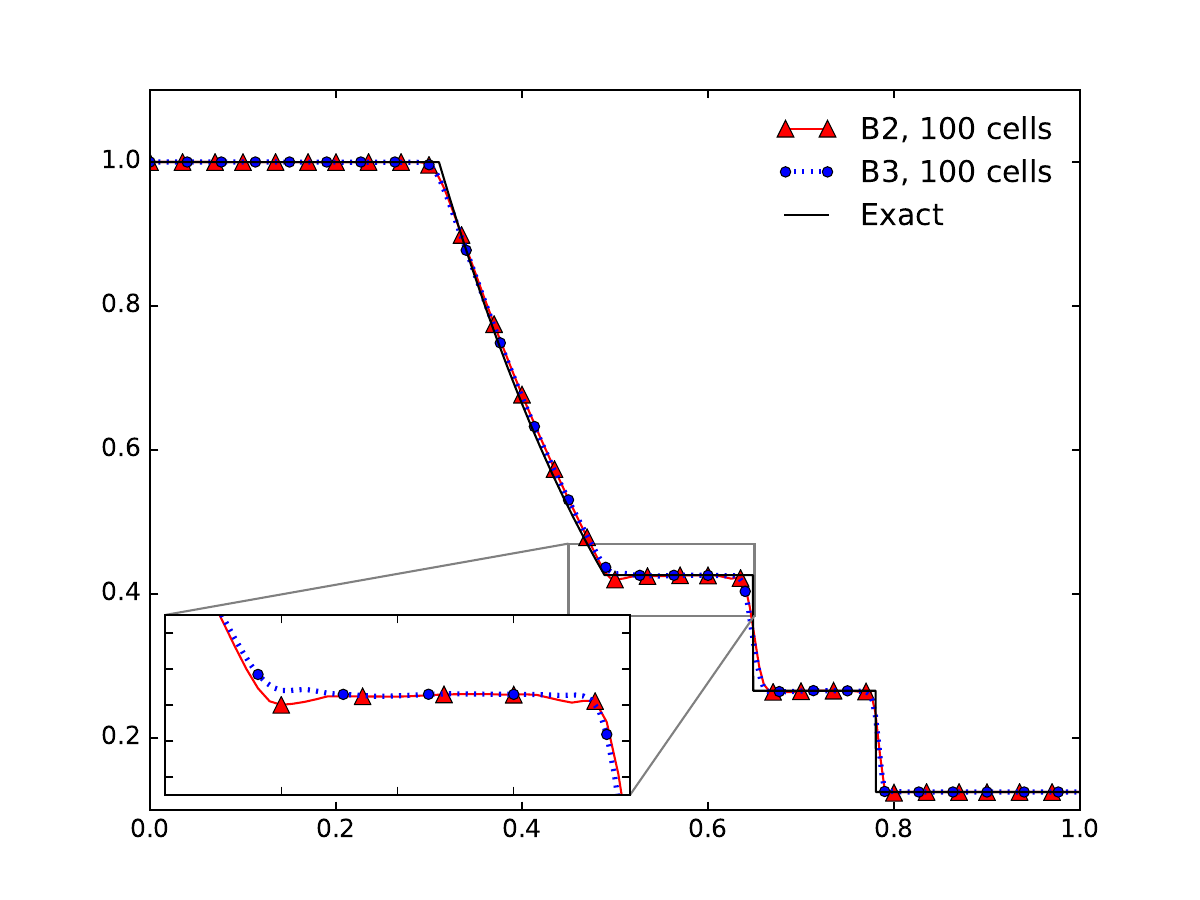}}
\subfigure[Density, reduced cascade]{\includegraphics[width=0.45\textwidth]{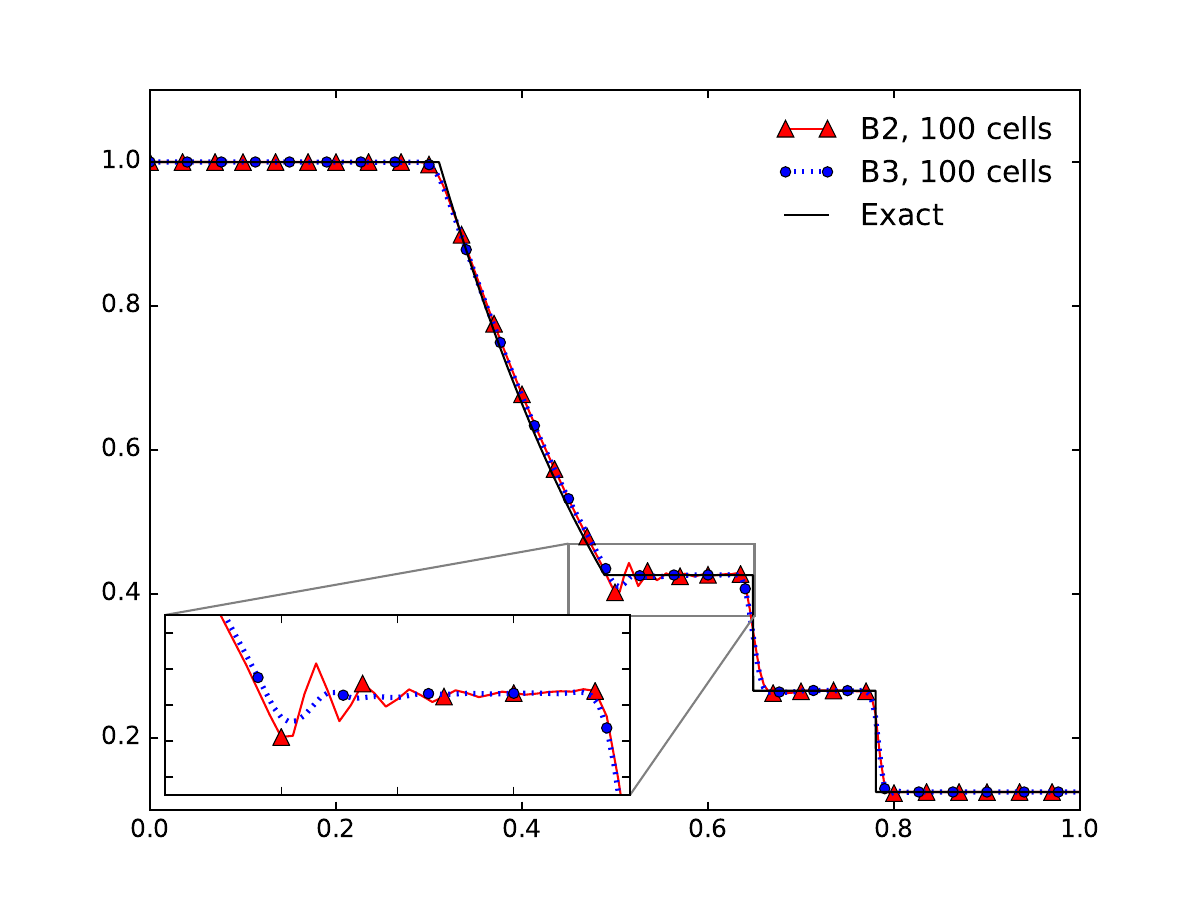}}\\
\subfigure[Density,complete cascade]{\includegraphics[width=0.45\textwidth]{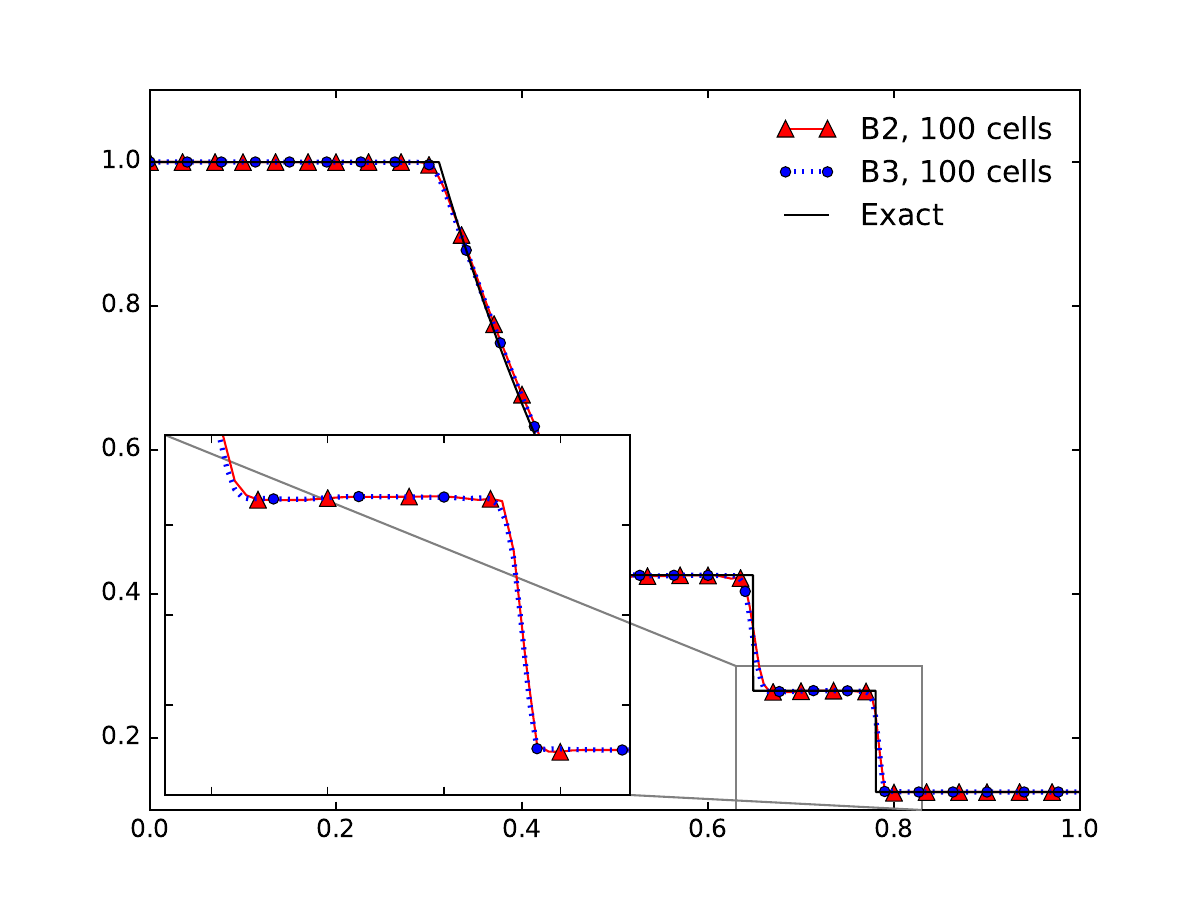}}
\subfigure[Density, reduced cascade]{\includegraphics[width=0.45\textwidth]{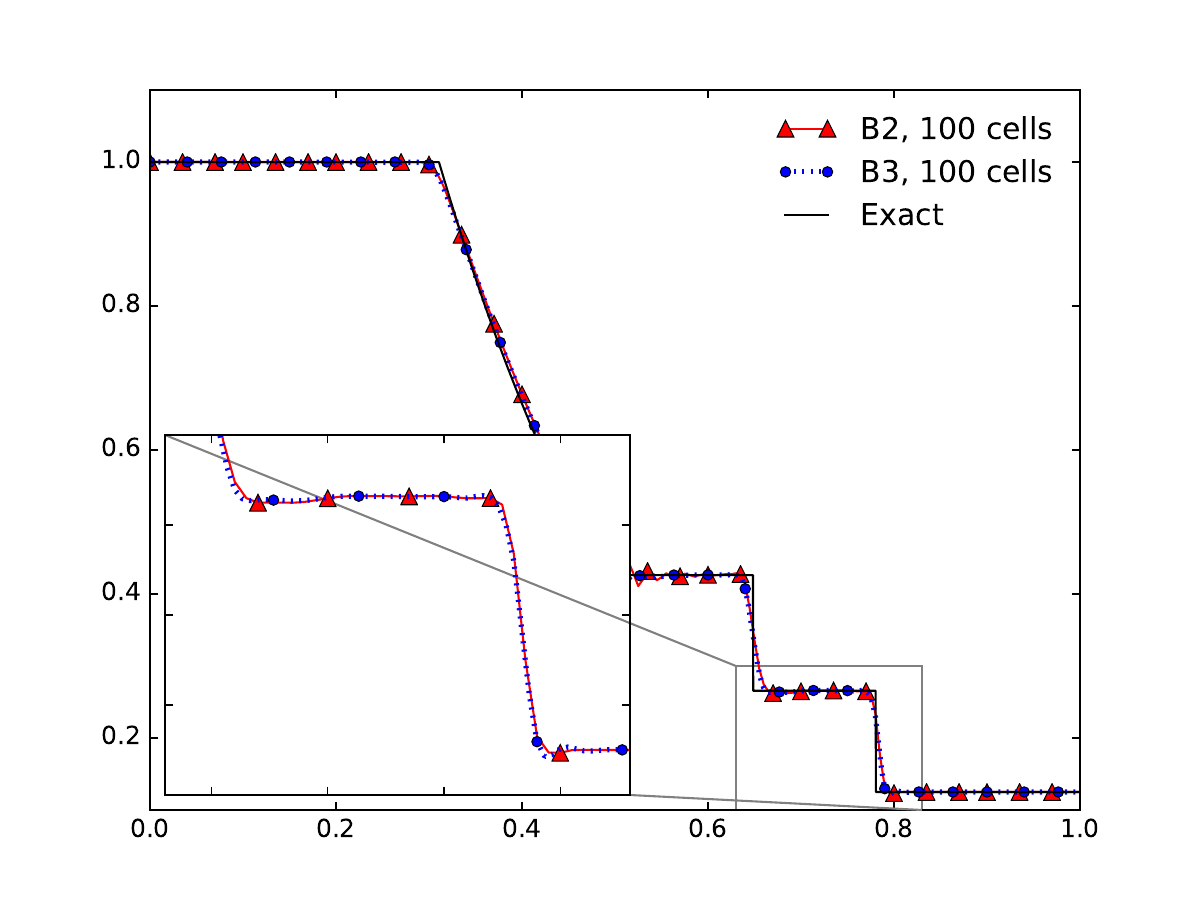}}
\caption{Sod 1D. Comparison of densities with zoom between a full (left) and reduced (right) cascade for $100$ cells on $\mathcal{B}^2$ and $\mathcal{B}^3$ at $T=0.16$}
\label{Fig:Sod_zoom_cascades} 
\end{figure}
 To compare the results obtained by the complete cascade, i.e. going from $s=2$ to $s=0$, i.e. from the GPJ to the parachute scheme, and a reduced version from $s=2$ to $s=1$, i.e. stopping before the parachute scheme. In Figure \ref{Fig:Sod_zoom_cascades} we show the results of this juxtaposition with some appropriate zooming of the areas of interest. We clearly see, as expected, how the solution results smoother for the complete cascade, and how the parachute scheme, thus, allows for a more truth-some approximation. 
 
 \subsubsection*{Comparing Different Numerical Admissibility Detection Techniques}	\label{NAD_Sod}
\begin{figure}[H] 
\centering
\subfigure[Density]{\includegraphics[width=0.45\textwidth]{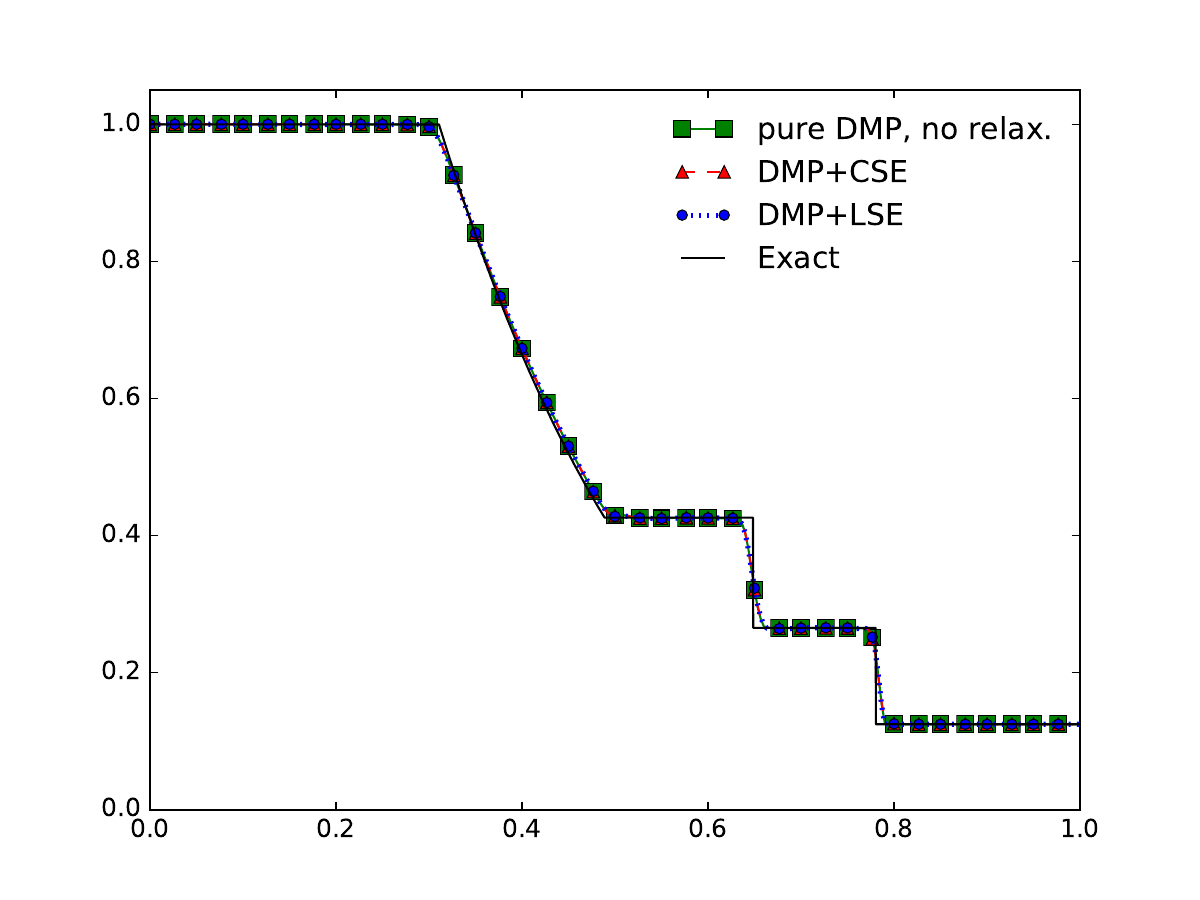}}
\subfigure[Density, zoom]{\includegraphics[width=0.45\textwidth]{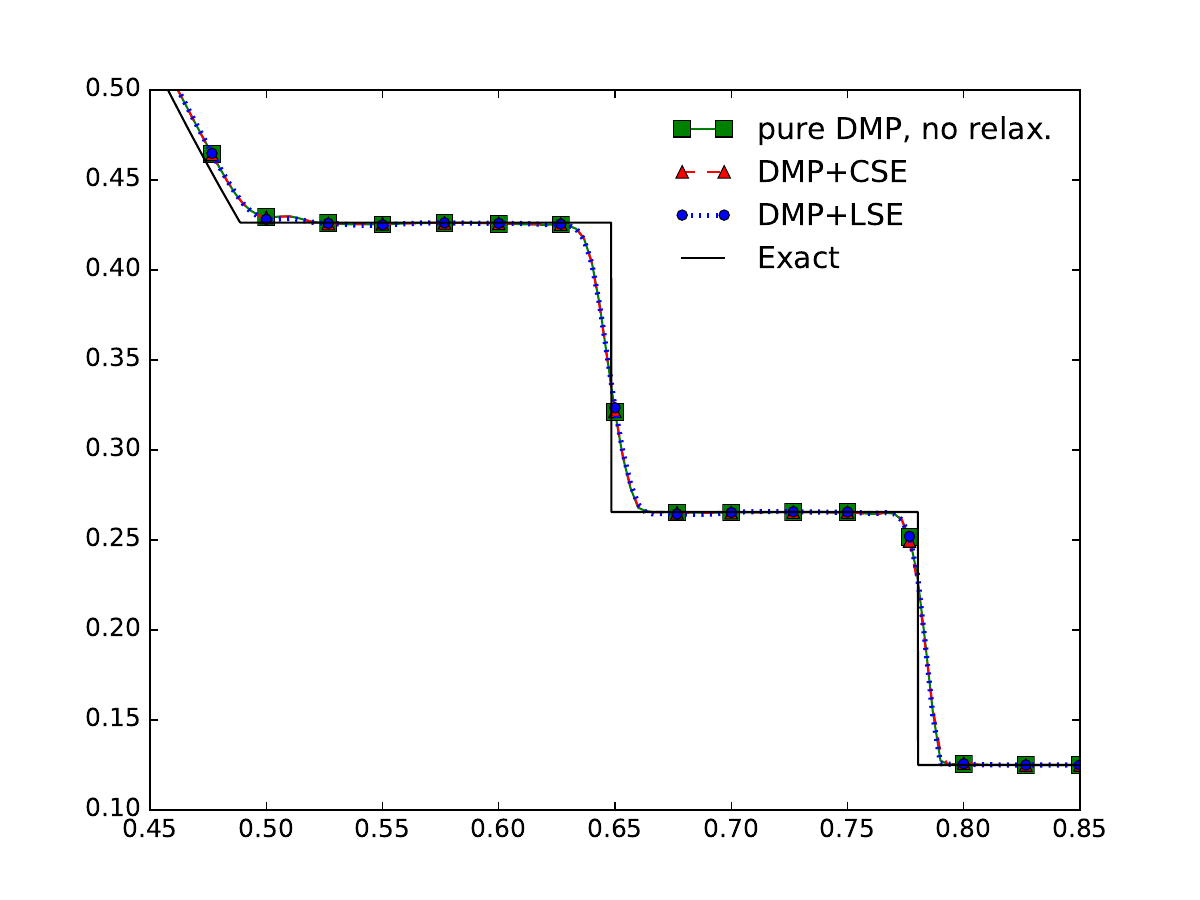}}
\caption{Sod 1D. Comparison of densities with zoom between different NAD criteria for $100$ cells on $\mathcal{B}^3$ at $T=0.16$}
\label{Fig:Sod_zoom_cascadesD} 
\end{figure}

Finally, hereafter is a remark on how different NAD criteria affect the solution, as shown in Figure \ref{Fig:Sod_zoom_cascadesD}. We consider here the pure DMP, with $\epsilon_1=\epsilon_2=0$ in \eqref{relaxation} without any smoothness extrema criteria. Along this approximation we compare the relaxed DMP with $\epsilon_1=10^{-3}$ and $\epsilon_2=0$ with the CSE together with the DMP with $\epsilon_1=\epsilon_2=0$ with the LSE criteria. What we deduce is that, apparently, in case of strong interacting discontinuities, there is no remarkable difference between these strategies. This does not hold nevertheless in case we are considering naturally oscillating solutions, as we shall see in the forthcoming section with the Shu-Osher problem.

\subsubsection{Shu-Osher Problem}

This test case, introduced in \cite{shuOsher1989}, is intended to demonstrate the advantages of high order schemes for problems involving some structure in smooth regions. In this test, we solve the Euler equations with initial conditions containing a moving Mach 3 shock wave which later interacts with periodic perturbations in density. The initial data for this problem is defined as follows:
\begin{equation*}
W=[\rho,u,p]=
\begin{cases}
[3.857143, 2.629369, 10.333333], \; &-5 \leq x \leq -4, \\
[1 + 0.2\sin(5x), 0, 1], \; &-4 < x \leq 5.
\end{cases}
\end{equation*}

\subsection*{Comparison of Different Numerical Admissibility Detection Techniques}

\begin{figure}[H]
\centering
\subfigure[Density]{\includegraphics[width=0.45\textwidth]{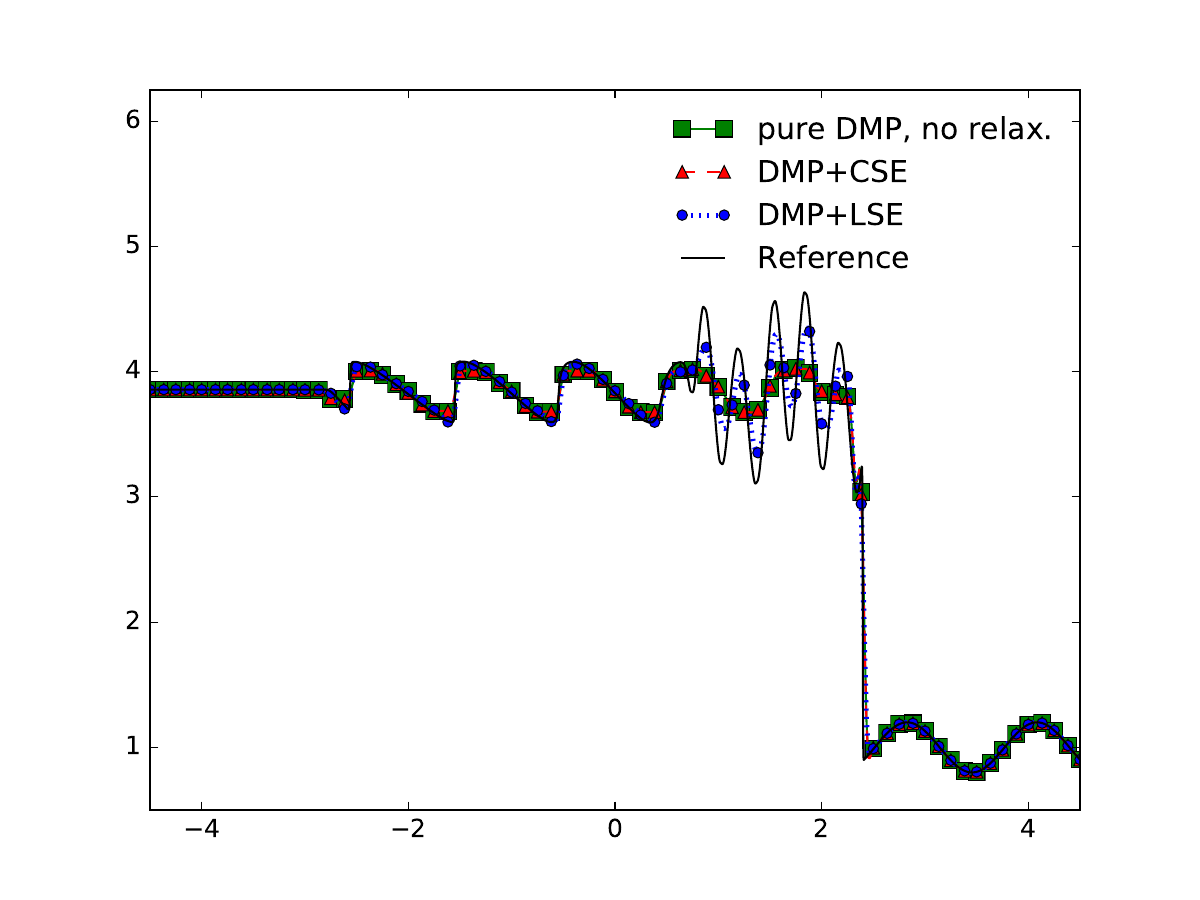}\label{Fig:SO_CompMethod_zoom}}
\subfigure[Density, zoom]{\includegraphics[width=0.45\textwidth]{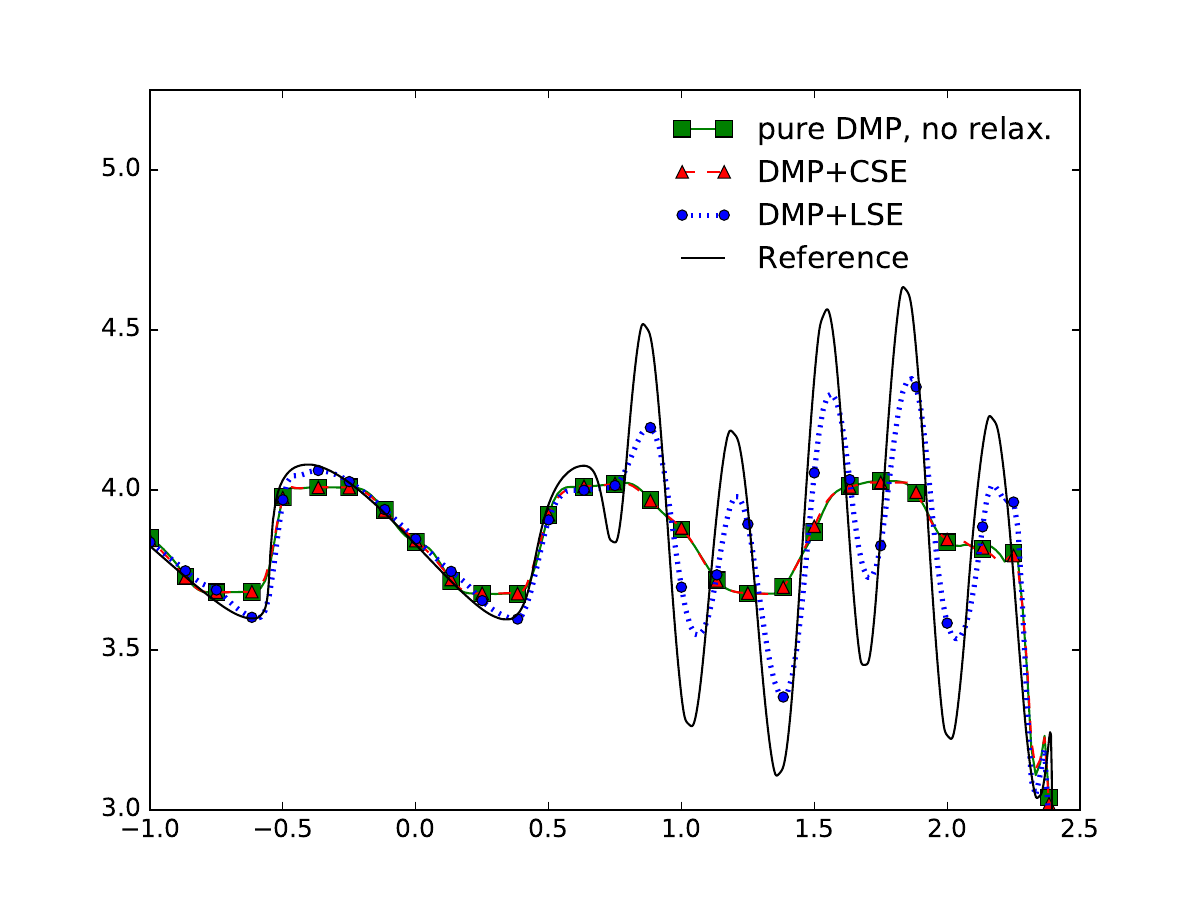}\label{Fig:SO_CompMethod_nozoom}}
\caption{Shu-Osher. Comparison of densities with zoom between different NAD criteria for $200$ cells on $\mathcal{B}^3$ at $T=1.8$.}
\label{Fig:SO_Detection_1}
\end{figure}

\begin{figure}[H]
\centering
\subfigure[Density]{\includegraphics[width=0.45\textwidth]{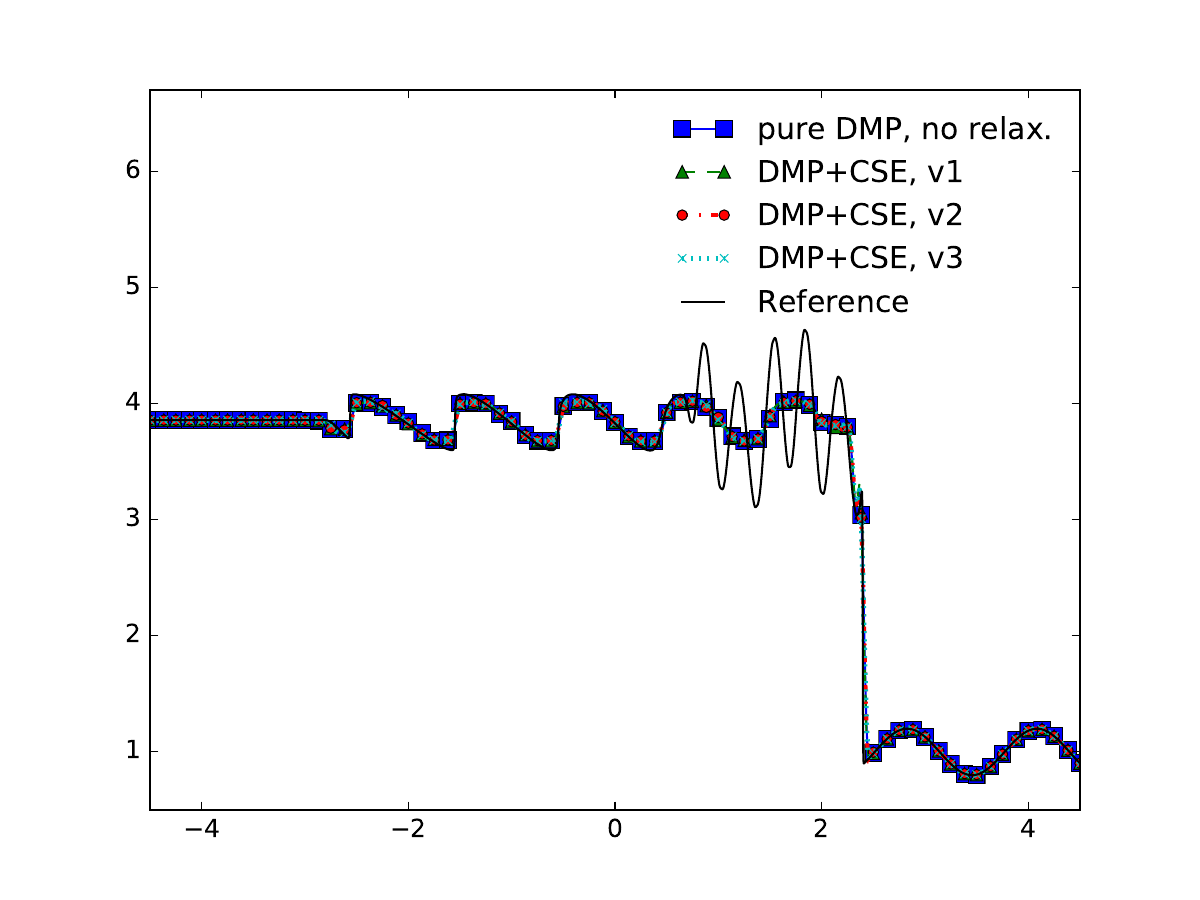}\label{Fig:SO_Raphcomp}}
\subfigure[Density, zoom]{\includegraphics[width=0.45\textwidth]{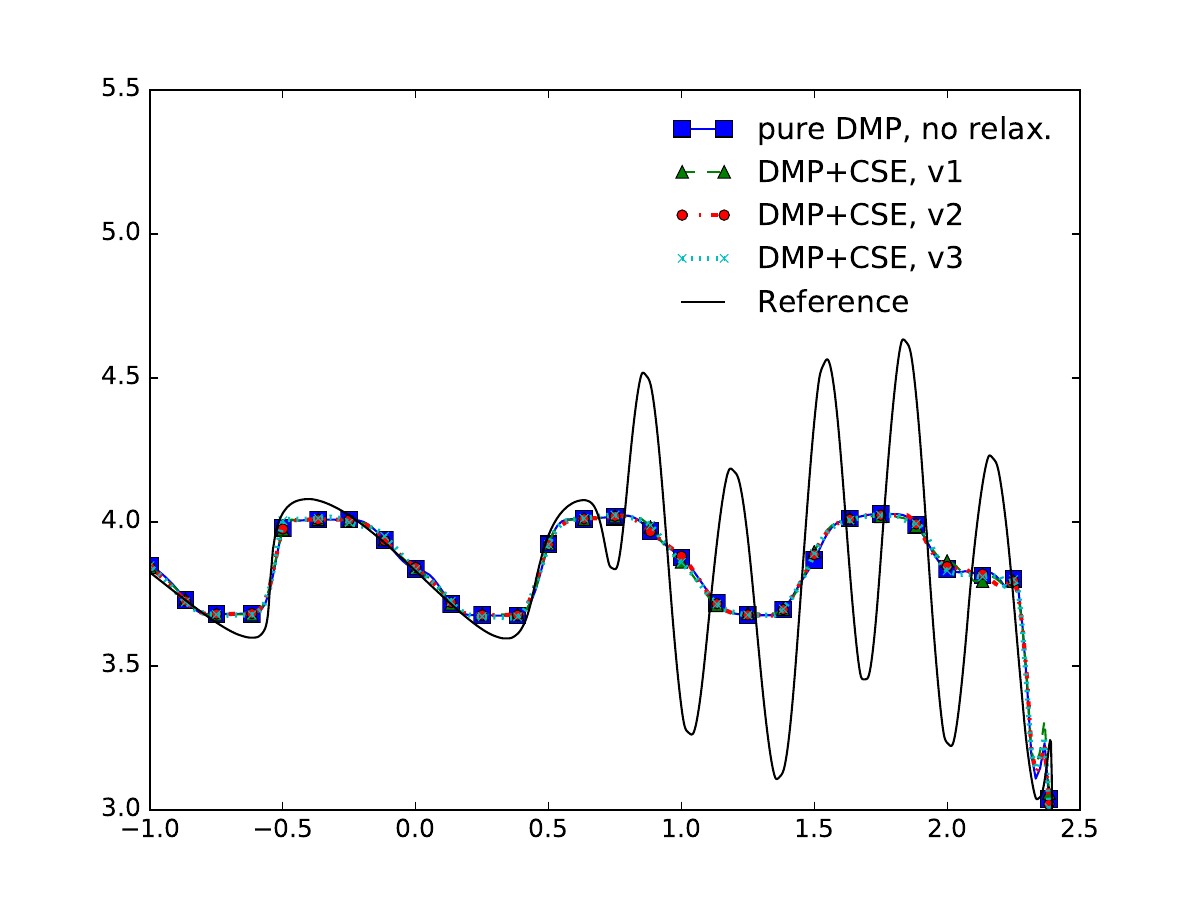}\label{Fig:SO_Raphcomp_zoom}}
\caption{Shu-Osher. Comparison of different classical smoothness detection criteria on $\mathcal{B}^3$ at $T=1.8$.}
\label{Fig:SO_Detection_2}
\end{figure}

Following the observations of Figure \ref{Fig:Sod_zoom_cascades}, we have compared in Figure \ref{Fig:SO_CompMethod_zoom} for the Shu-Osher problem the same detection criteria, i.e. the pure DMP without the smoothness extrema detection and $\epsilon_1=\epsilon_2=0$, the relaxed DMP with $\epsilon_1=10^{-3}$ and $\epsilon_2=0$ with the CSE, and moreover, the DMP with $\epsilon_1=\epsilon_2=0$ with the LSE criteria.
While not huge differences appear to be between the pure DMP and the DMP with CSE, these two approaches are remarkably giving less accurate approximations compared to the one provided by the DMP with LSE.
Furthermore, we have investigated if the dissipative behaviour, caused by an exceeding detection of areas which should indeed not be dissipated at all, might improve, in case of different relaxation parameters $\epsilon_1$ and $\epsilon_2$ of \eqref{relaxation}. Figure \ref{Fig:SO_Raphcomp_zoom} summarizes this study for the chosen parameters: 'DMP+CSE, v1' corresponds to $\epsilon_1=\epsilon_2=0$, 'DMP+CSE, v2' to $\epsilon_1=10^{-3}$ and $\epsilon_2=0$; 'DMP+CSE, v3' to $\epsilon_1=10^{-3}$ and $\epsilon_2=10^{-4}$.
The resulting differences allow us to remark that there is apparently no choice that allows for a sharper solution within this test case and for our novel approximating technique.

Indeed, changing the values of $\epsilon_1$ or $\epsilon_2$ has no impact. What allows the LSE to improve the quality of the solution is the actual meaning behind the LSE: the numerical oscillation detection is based on a change of sign or jump of the approximation's derivatives. This is an extremely fine way to catch the oscillation, whereas the CSE depends on a second derivative. The result is that the LSE is able to capture the natural peaks, while the CSE cuts the peaks as it perceives them as numerical non-physical oscillations.
This motivates us to consider in the forthcoming tests the NAD criteria composed by the DMP with the LSE criteria.

\subsubsection*{Detecting Technique in practice}
\begin{figure}[H]
\centering
\subfigure[$1$st iteration, Spatial Scheme Indicator]{\includegraphics[width=0.45\textwidth]{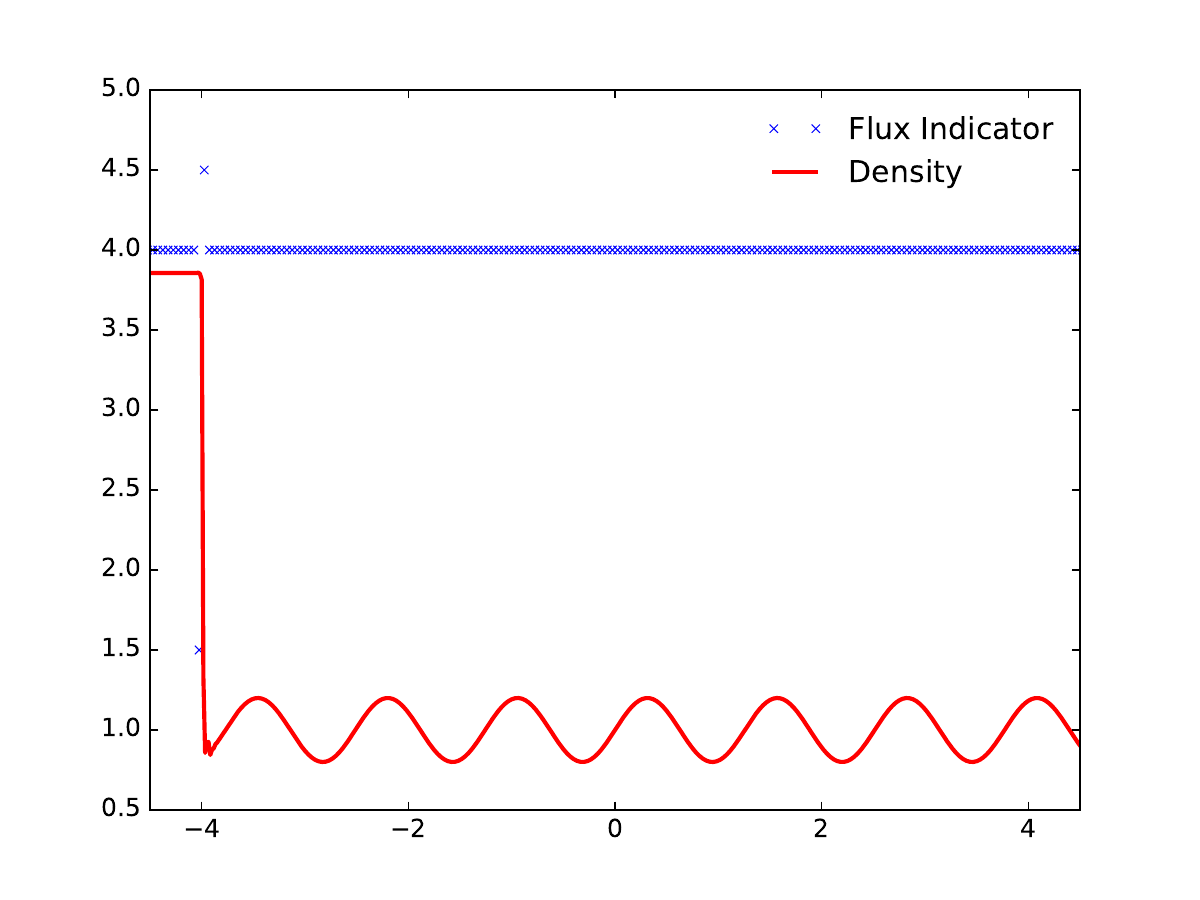}\label{Fig:SO_detection_1}}
\subfigure[$10$th iteration, Spatial Scheme Indicator]{\includegraphics[width=0.45\textwidth]{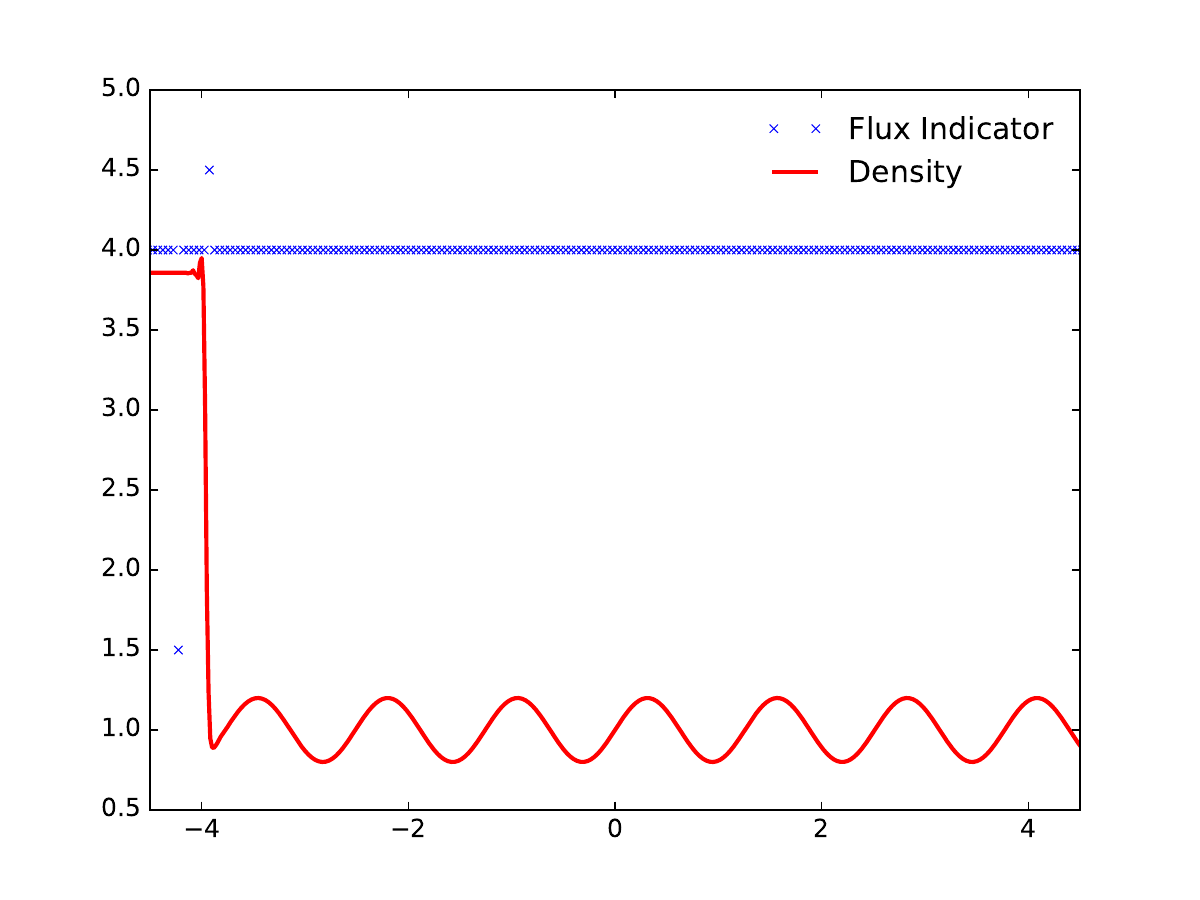}\label{Fig:SO_detection_10}}
\caption{Shu-Osher. Considered spatial scheme for the $1$st and $10$th iteration in time on  $\mathcal{B}^3$ with $200$ cells.}
\label{Fig:SO_Detection}
\end{figure}

To be able to show the detection capabilities of the considered NAD, from now on DMP+LSE, and show the overall MOOD procedure for this class of benchmark problems, we show the very first iteration in time, i.e. $n=1$ and the $n=10$-th iteration in time. In Figure \ref{Fig:SO_Detection}, we display the approximated density w.r.t the spatial scheme indicator, where this indicator is at $4.5$ for the Rusanov with PSI limiting and stabilizing jump terms ($s=1$ i.e. RPJ) and at $1.5$ for the Rusanov (Rus) scheme and else $4$ for the Galerkin with jump. One can note how in the very first iteration only nodes at the very shock interface are being treated by a more dissipative scheme. Here we can also appreciate the capability of the detection to ignore, indeed, natural oscillations within the solution.

\subsubsection*{Comparison between different orders and mesh convergence}
\begin{figure}[H]
\centering
\subfigure[$200$ cells]{\includegraphics[width=0.45\textwidth]{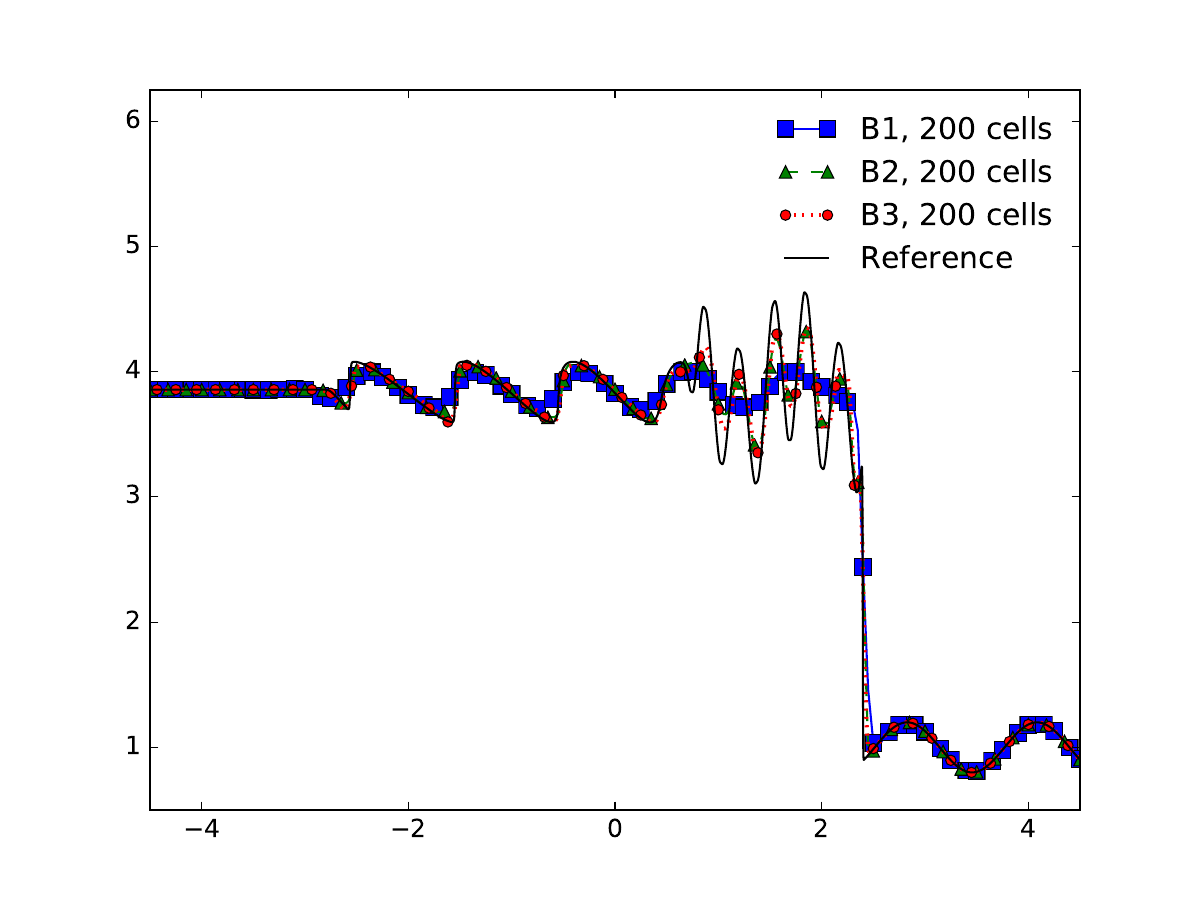}\label{Fig:SO_400_Bs}}
\subfigure[$200$ cells, zoom]{\includegraphics[width=0.45\textwidth]{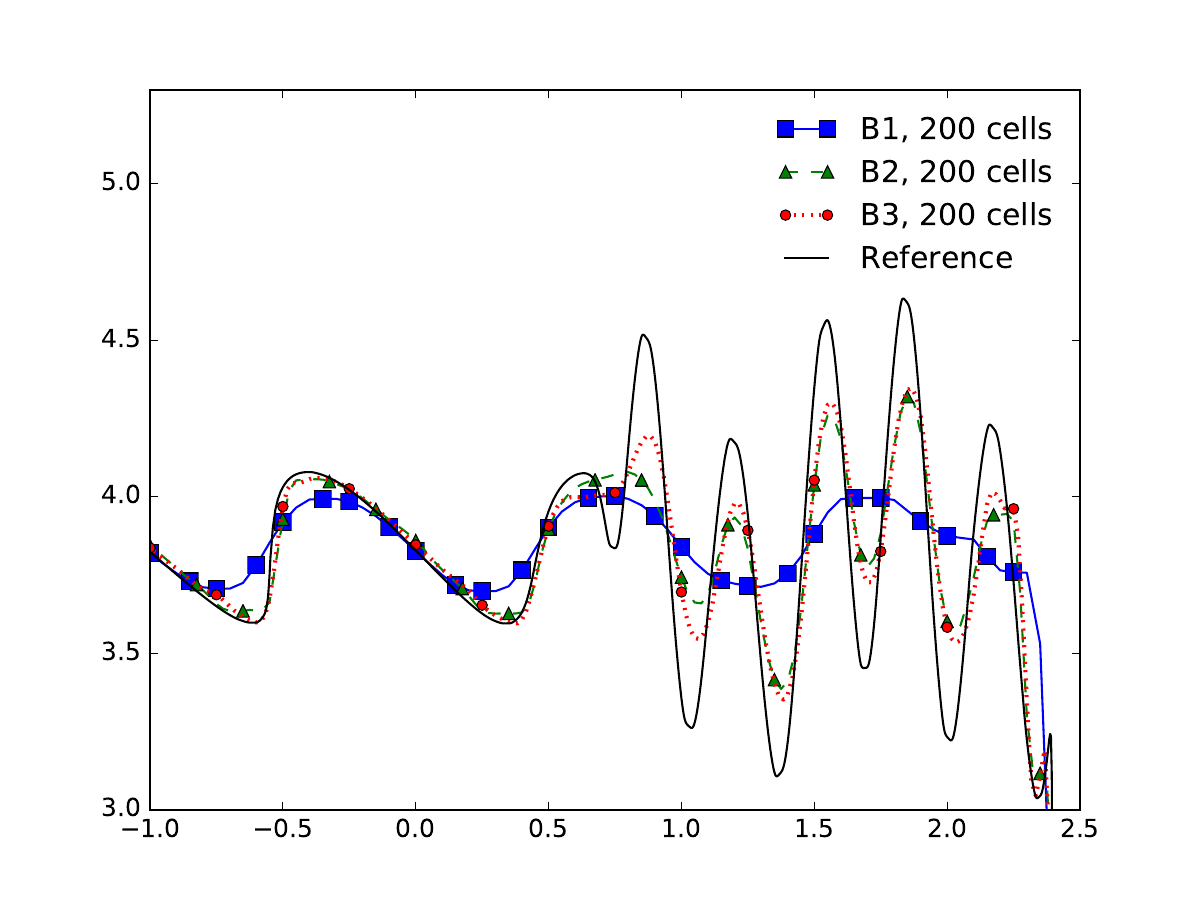}\label{Fig:SO_400_Bs_zoom}}\\
\subfigure[$800$ cells]{\includegraphics[width=0.45\textwidth]{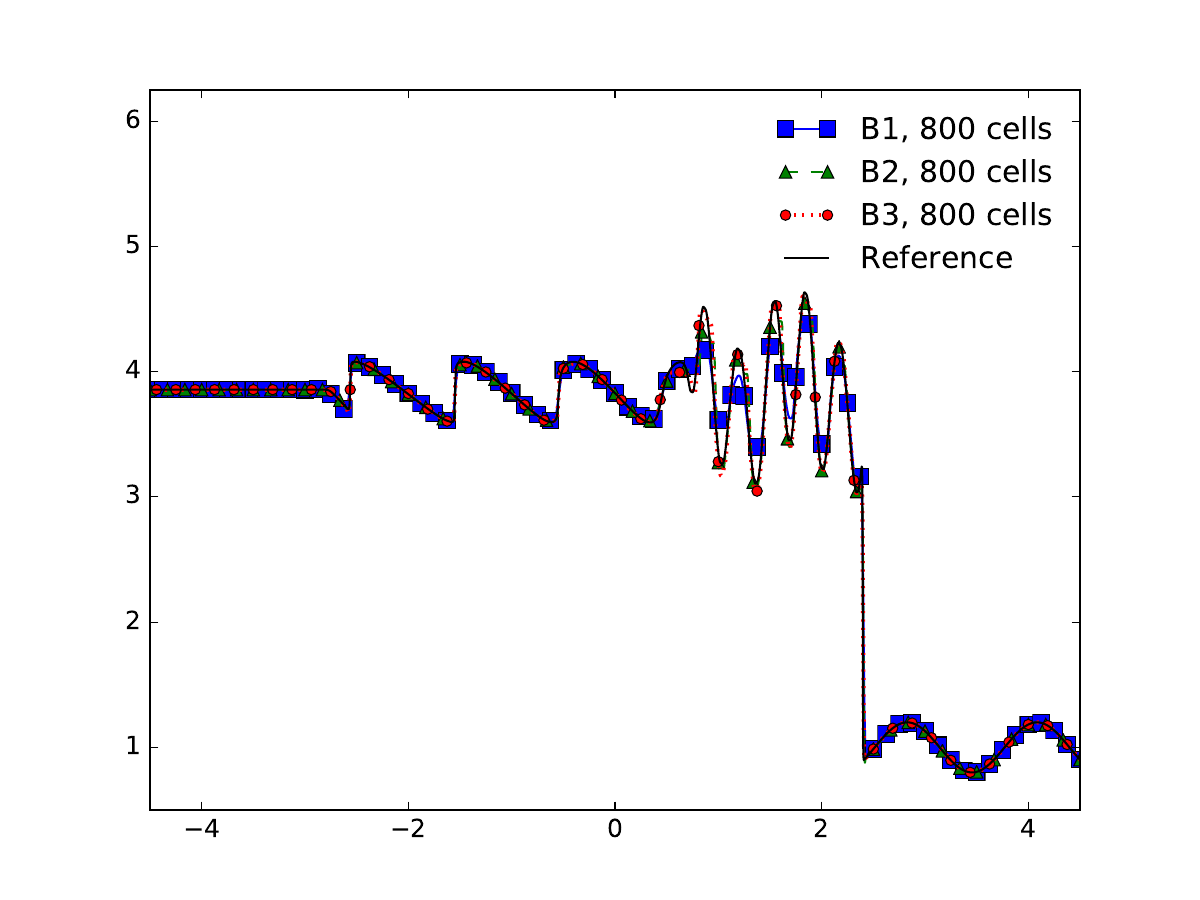}\label{Fig:SO_800_Bs}}
\subfigure[$800$ cells, zoom]{\includegraphics[width=0.45\textwidth]{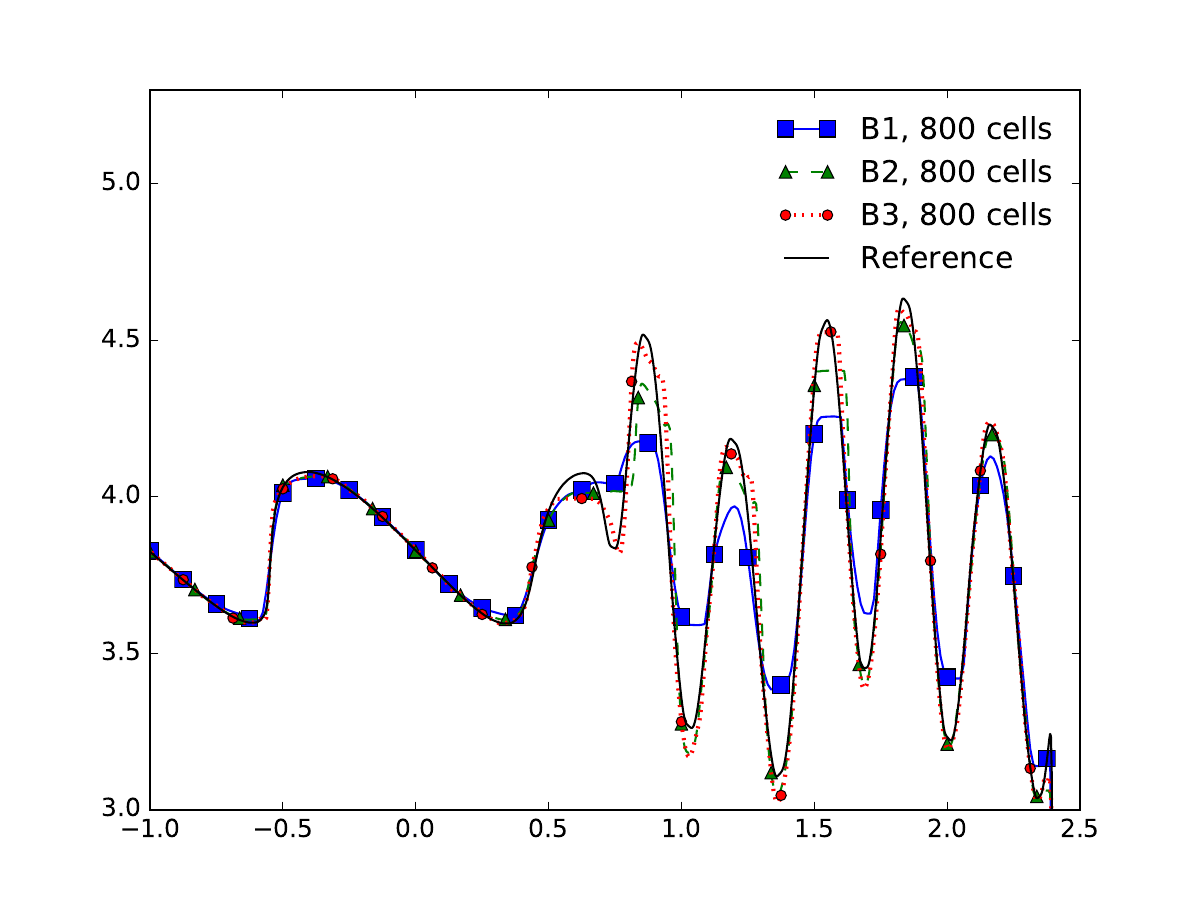}\label{Fig:SO_800_Bs_zoom}}\\   
\caption{Shu-Osher. Comparison for the density between $\mathcal{B}^1,\;\mathcal{B}^2,\,\mathcal{B}^3$ at $T=1.8$.}
\label{Fig:SO_BS}
\end{figure}

The more accurate approximation obtained by the fourth order scheme in comparison to the second and third order is clearly visible in this benchmark problem in Figure \ref{Fig:SO_BS}, and increasing the number of mesh elements within the domain strongly increases the quality of the solution, too.

\subsubsection*{The proposed method with respect to an ``a priori technique''}

\begin{figure}[H]
\centering
\subfigure[$200$ cells]{\includegraphics[width=0.45\textwidth]{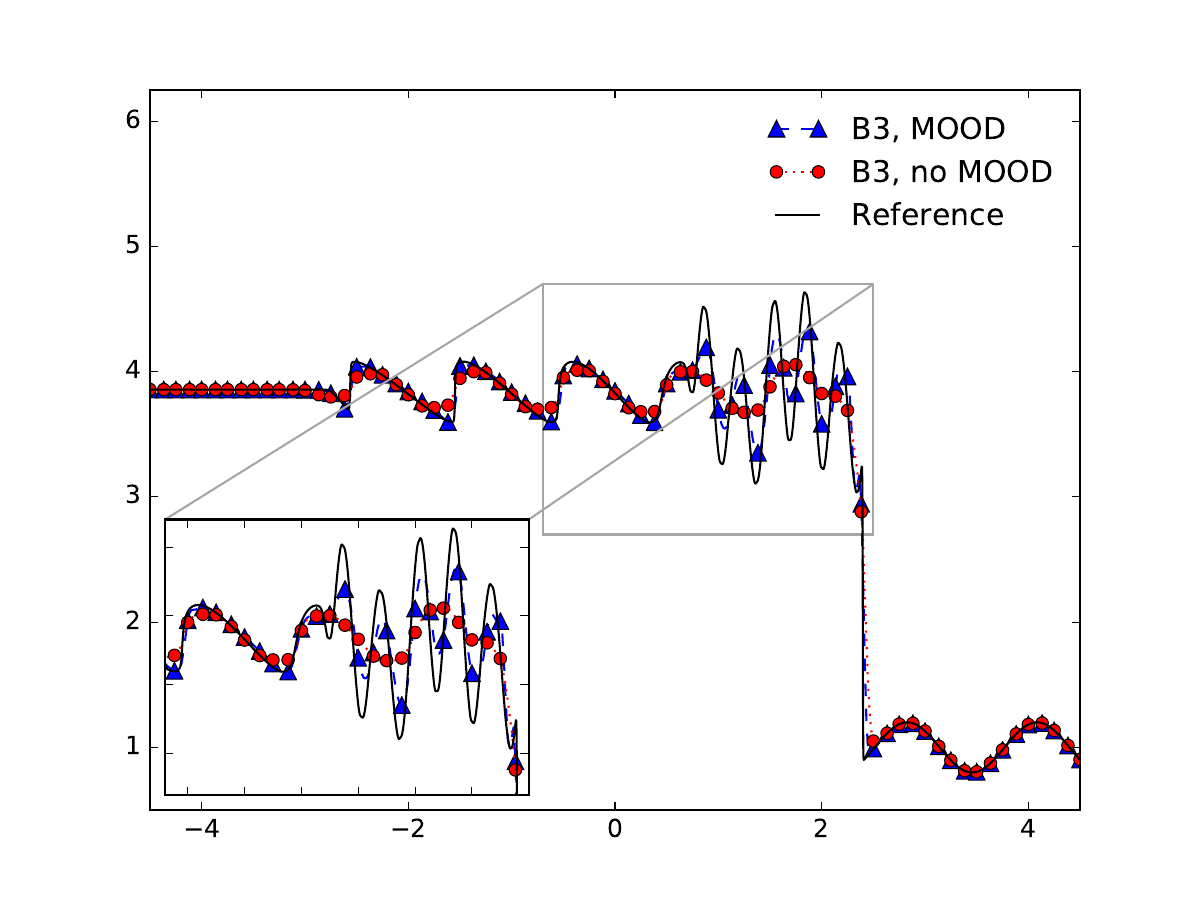}\label{Fig:SO_200_OldNew}}
\subfigure[$800$ cells]{\includegraphics[width=0.45\textwidth]{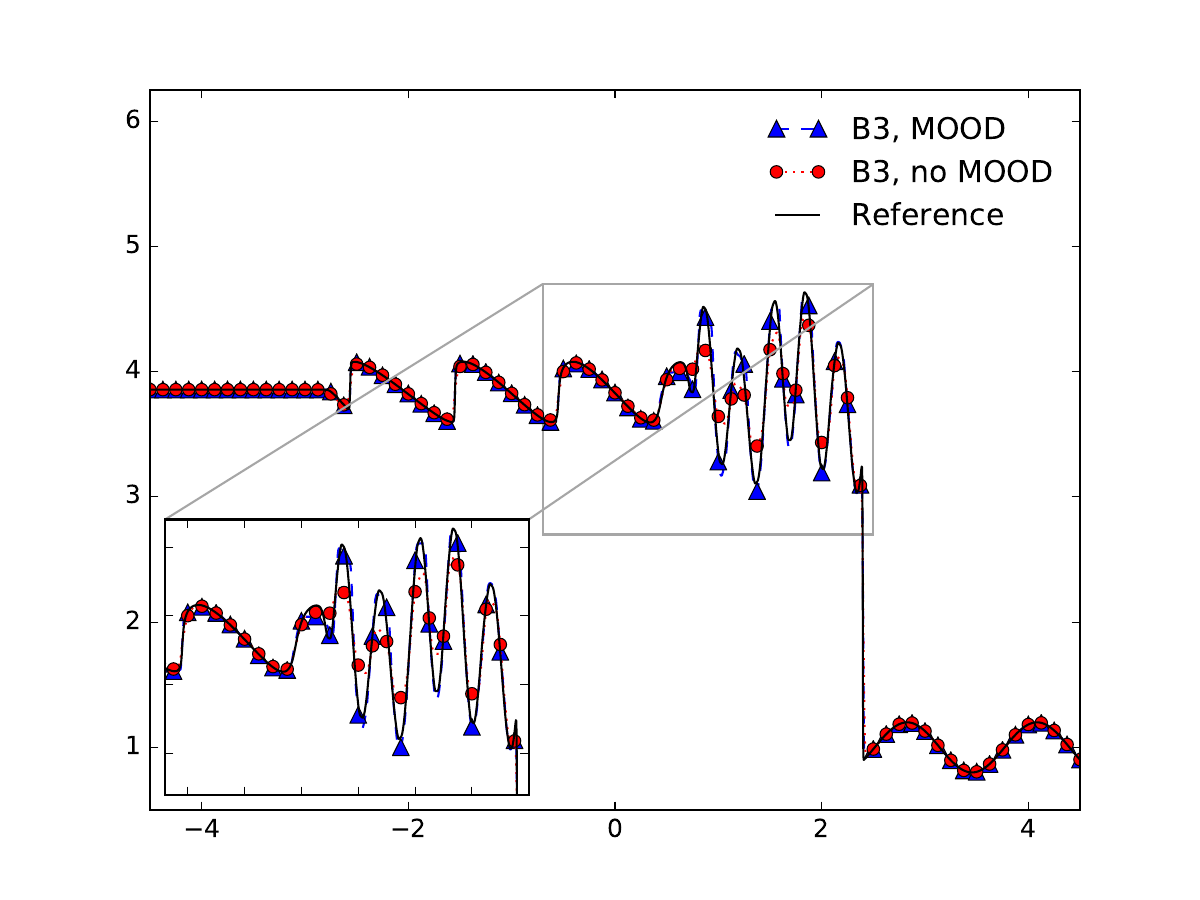}\label{Fig:SO_800_OldNew}}
\caption{Shu-Osher. Comparison of the densities between the MOOD and non-MOOD schemes on $\mathcal{B}^3$ at $T=1.8$.}
\label{Fig:SO_MeshComp_OldNew}
\end{figure}
Finally, comparing in Figure \ref{Fig:SO_MeshComp_OldNew} the novel methodology for $\mathcal{B}^3$ with the pure RPJ without any ``a posteriori limiting'', provides us with the acknowledgement that, again, also in the case of problems with both natural oscillations and shocks, we are now able to provide a more detailed approximation. This holds even when we increment considerably our mesh size, as in Figure \ref{Fig:SO_800_OldNew}.


\subsubsection{Woodward-Colella Problem}

The interaction of blast waves is a standard low energy benchmark problem involving strong shocks reflecting from the walls of the tube with further mutual interactions. The initial data is the following:
\begin{equation*}
(\rho_0,u_0,p_0) =
\begin{cases}
[1, 0, 10^3], \; &0 \leq x \leq 0.1, \\
[1, 0, 10^{-2}], \; &0.1 < x < 0.9, \\
[1, 0, 10^2], \; &0.9 \leq x \leq 1.
\end{cases}
\end{equation*}

\begin{figure}[H]
\centering
\subfigure[$400$ cells]{\includegraphics[width=0.45\textwidth]{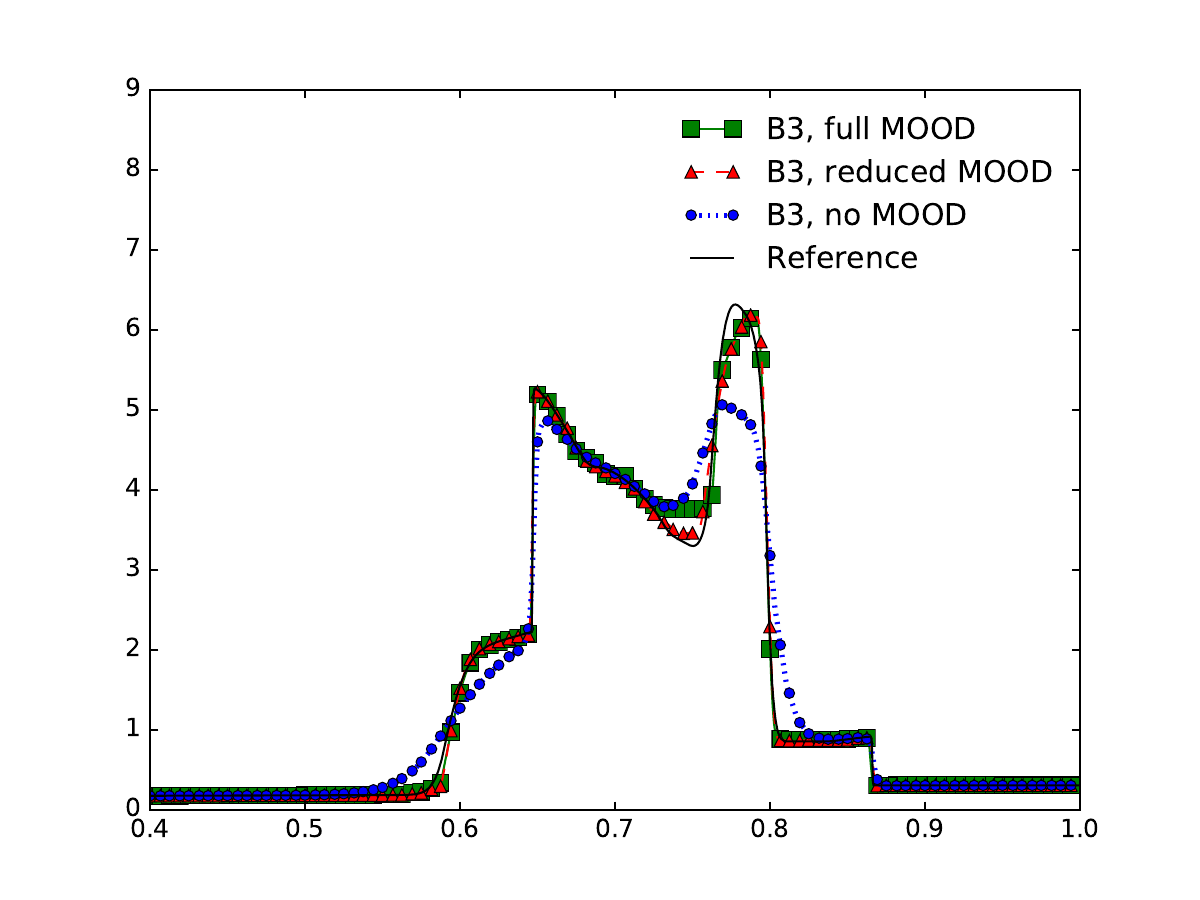}}
\subfigure[$800$ cells]{\includegraphics[width=0.45\textwidth]{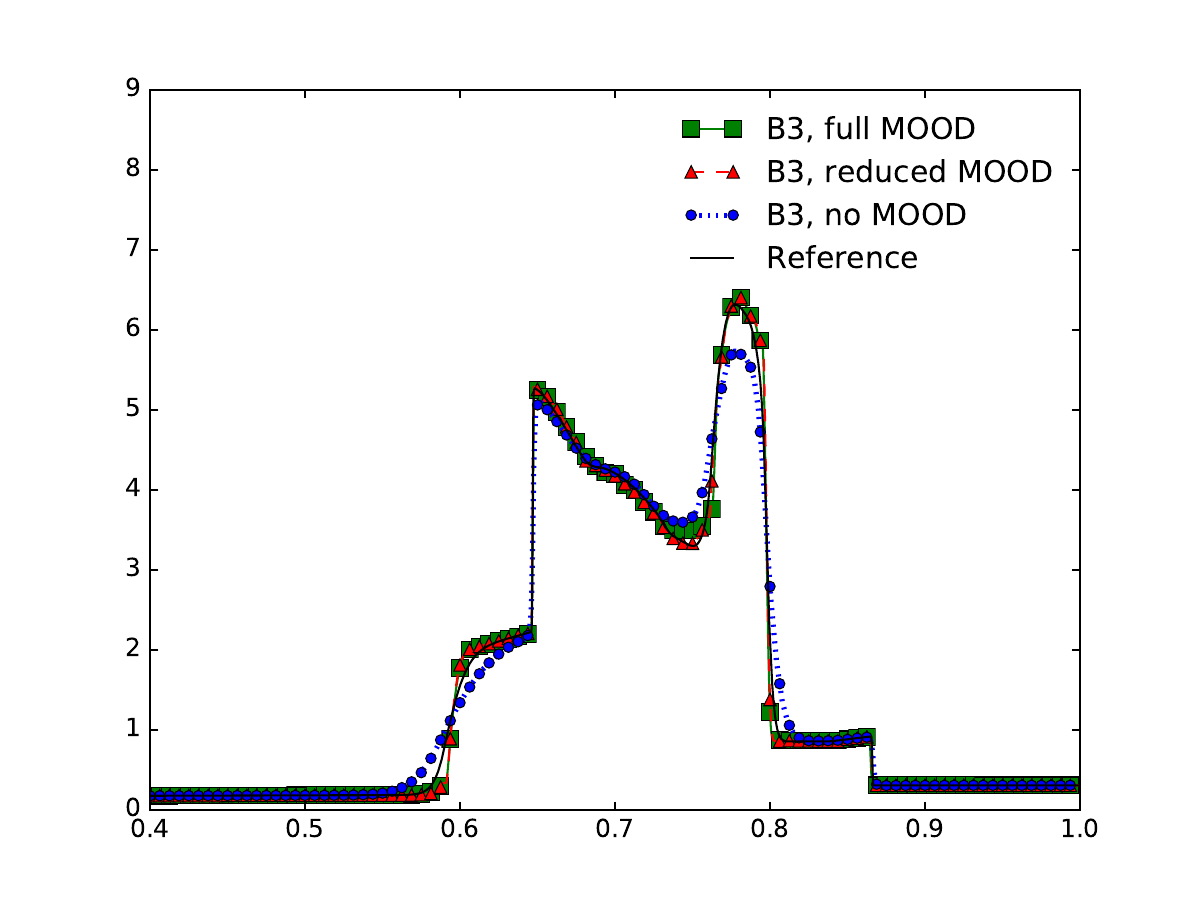}} 
\caption{Woodward-Colella. Comparison between the MOOD and non-MOOD schemes on $\mathcal{B}^3$ at $T=0.038$.}
\label{Fig:WC}
\end{figure}

Comparing the results of our novel method for both the reduced, i.e. parachute at $s=1$, and the full cascade, against the one of a pure RPJ without the MOOD strategy, we can observe in Figure \ref{Fig:WC} how the solution is well approximated already on a 400 cell mesh with B3, and further mesh refinement shows the expected convergence to the exact solution. 
The plots show a very good overall behaviour, of the numerical scheme even for this extremely demanding test case. Some wiggles are due probably to an extreme flagging activity that locally activates the parachute scheme and are one topic of future research. 
We remark that the expected solution has been computed via the pure RPJ on a mesh of $3200$ cells, reasoning why the novel method reaches and surpasses the so-called expected solution.

\begin{figure}[H]
 \centering
\subfigure[Shu-Osher, density, $200$ cells, $\mathcal{B}^3$]{\includegraphics[width=0.45\textwidth]{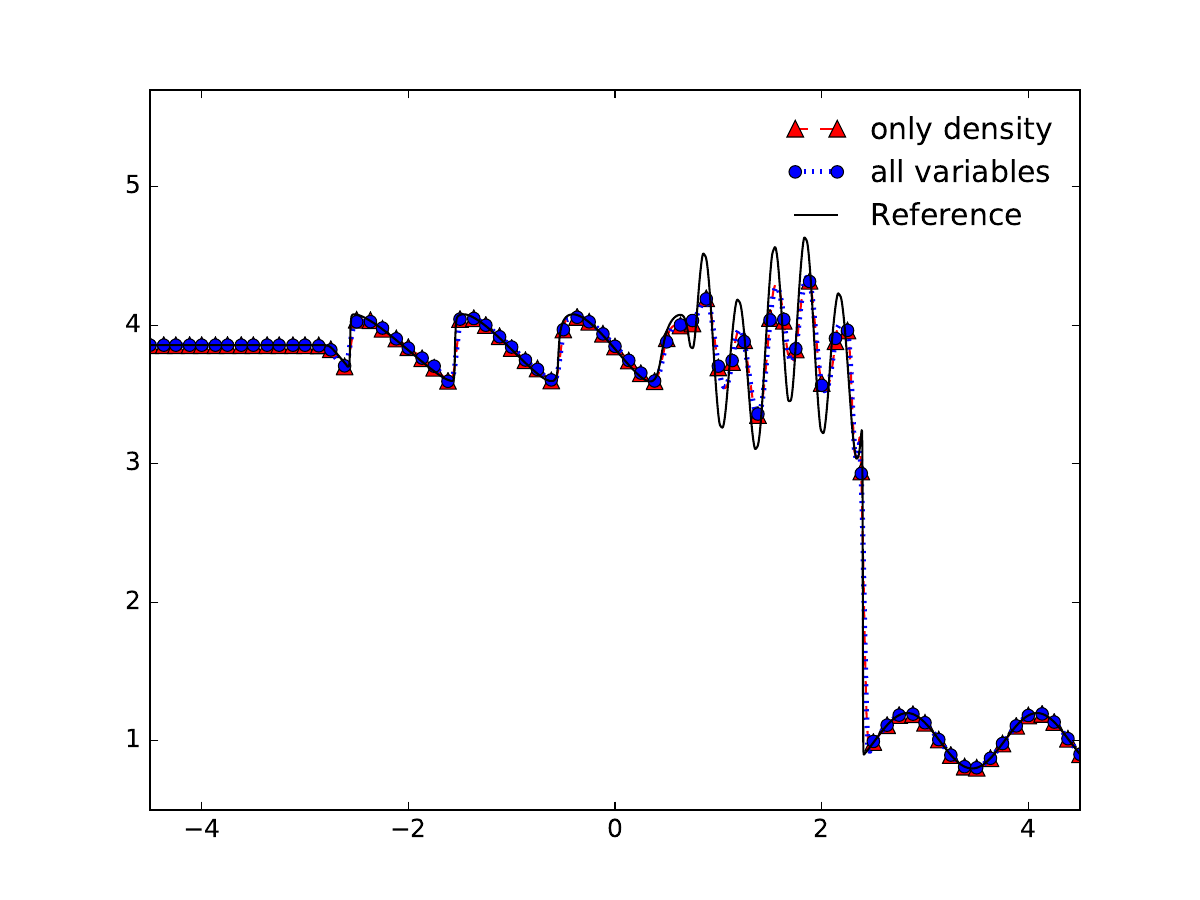}}
\subfigure[Woodward-Colella, density, $400$ cells, $\mathcal{B}^3$]{\includegraphics[width=0.45\textwidth]{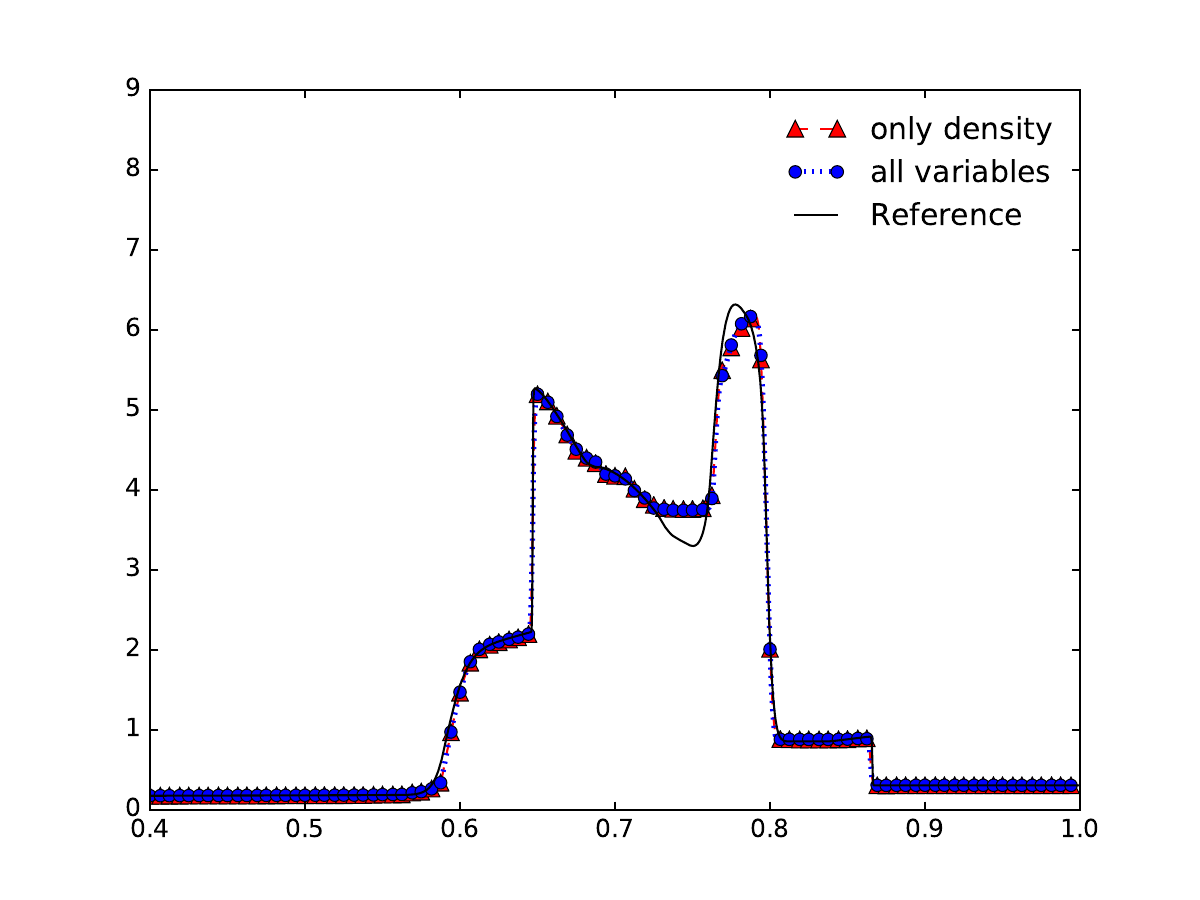}} 
\caption{Shu-Osher \& Woodward-Colella. Detection applied uniquely to the density variable vs. to all conserved primitive variables on $\mathcal{B}^3$.}
\label{Fig:SO_WC_vars} 
\end{figure}
As a final comment to the numerical section for the one-dimensional test cases, we compare in Figure \ref{Fig:SO_WC_vars} the results obtained on two different test cases for two different detection ideas.
One idea carries out all detection criteria, i.e. the CAD, plateau and NAD ones on the density only, while the other considers the CAD, plateau and NAD for all primitive variables. Following the Remark  \ref{Remark_vars}, there is no evident difference, and as such, we have considered in all the previous and forthcoming computations, the case of the sole density as parameter to be detected. 
\subsection{Numerical Results for 2D Test Cases}
\subsubsection{2D Shock Vortex Problem}
\begin{figure}[H]
\begin{center}
\subfigure[Density]{\includegraphics[width=0.49\textwidth]{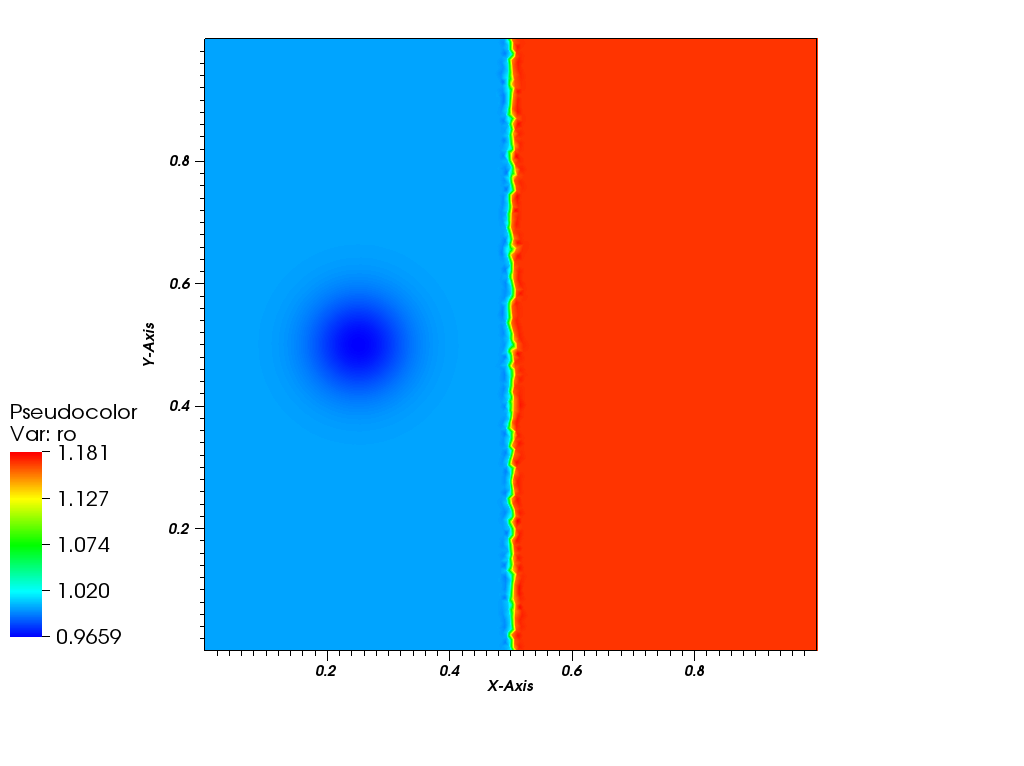}\label{shockVortex_B1_It1}}
\hspace{-1.5cm} \subfigure[MOOD flag]{\includegraphics[width=0.49\textwidth]{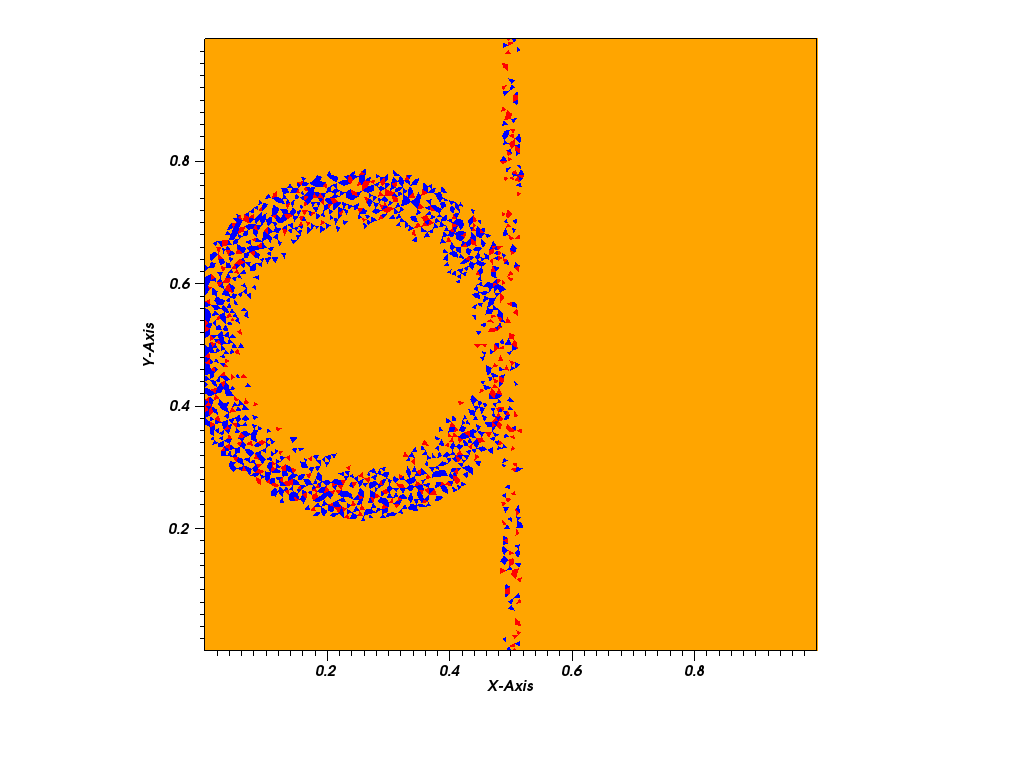}\label{shockVortex_B1_detection}}
\caption{Shock Vortex Interaction. Results for $\mathcal{B}^1$ at the first iteration with the density on the left and the flag activation on the right, with orange s=2, red s=1 and blue s=0.}\label{shockVortex_B1}
\end{center}
\end{figure}

In the following, we have considered a classical benchmark that encloses both a verification on the high order preservation by considering a smooth area given by a vortex,  and a robustness check across a shock, which has been studied in several works, such as \cite{Ellzey1995}.
In particular, we consider a domain in x-y direction given by a $1 \times 1$ m lenght area, where at $0.5$m on the x-axis a diaphragm divides a right hand area with the conditions set to  
  $\rho_R=\rho_0 \left(\frac{\gamma-1}{\gamma+1}+\frac{2}{(\gamma+1)M^2}\right)^{-1},$
    $u_R=\rho_0 u_0\rho_R^{-1},$ $v_R=0$ and $P_R=\rho_R^{-1} \rho_0^{-1}$, where $\rho_0=1$,  $u_0=M\sqrt{\gamma}$,  $M=1.1$ and $\gamma=1.4$.
On the left hand side of the diaphragm, we define 
         $\rho_L= \rho_0 \left( \frac{T}{T_0}\right)^{1/\gamma}$,
         $u_L=u_0+\epsilon  y\kappa$, $v_L=v_0-\epsilon x\kappa$ and 
         $ P_L=\rho T$.
Here $  T=1-\frac{\gamma-1}{4\alpha\gamma} \epsilon^2  \kappa^2$ and the parameters $T_0=1$,   $v_0=0$,  $\epsilon=0.3$ and $\alpha=0.204$.  Further, $\kappa=\exp\left(\alpha\left(1-\left[z_x^2+z_y^2\right]\right)\right)$,  with
 $z_x=20 \left(x-0.25\right)$ and
 $z_y=20 \left(y-0.5 \right)$.
 Considering an unstructured mesh composed by $N=26918$ cells, with a CFL of $0.125$, we test the shock vortex benchmark on $\mathcal{B}^1$, $\mathcal{B}^2$ and $\mathcal{B}^3$. In particular, we chose  for \eqref{phi_burman} the couples for $\mathcal{B}^1$ as $\theta_1=0.1$  $\theta_2=0$,  for $\mathcal{B}^2$ as $\theta_1=0.01$  $\theta_2=0$ and $\mathcal{B}^3$ as $\theta_1=0.0001$  $\theta_2=0.001$.
The activation of the detection in 2D can be seen in Figure \ref{shockVortex_B1} for the very first iteration in time for $\mathcal{B}^1$, with the actual solution for the density displayed on the left and the detection displayed on the right.  For Figure \ref{shockVortex_B1_detection}, orange represents the least dissipative scheme with s=2, corresponding to the stabilized Galerkin schema of \eqref{GalerkinJump}, red represents s=1, i.e.  the PSI scheme with the stabilization, c.f. \eqref{RusPsiJump}, and, finally, s=0 the first order scheme given by the local Lax-Friedrichs one, recast from \eqref{Rus}. 

\begin{figure}[H]
\begin{center}
\subfigure[$\mathcal{B}^2$]{\includegraphics[width=0.49\textwidth]{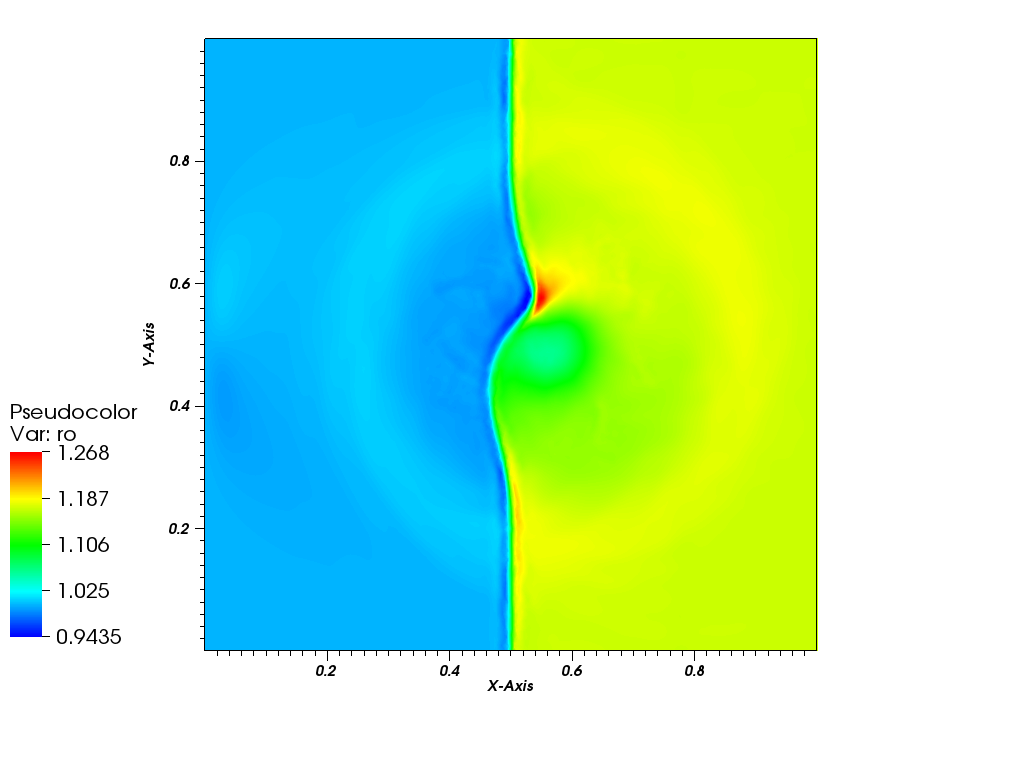}\label{shockVortex_B2_density}}
\hspace{-1.5cm} \subfigure[$\mathcal{B}^3$]{\includegraphics[width=0.49\textwidth]{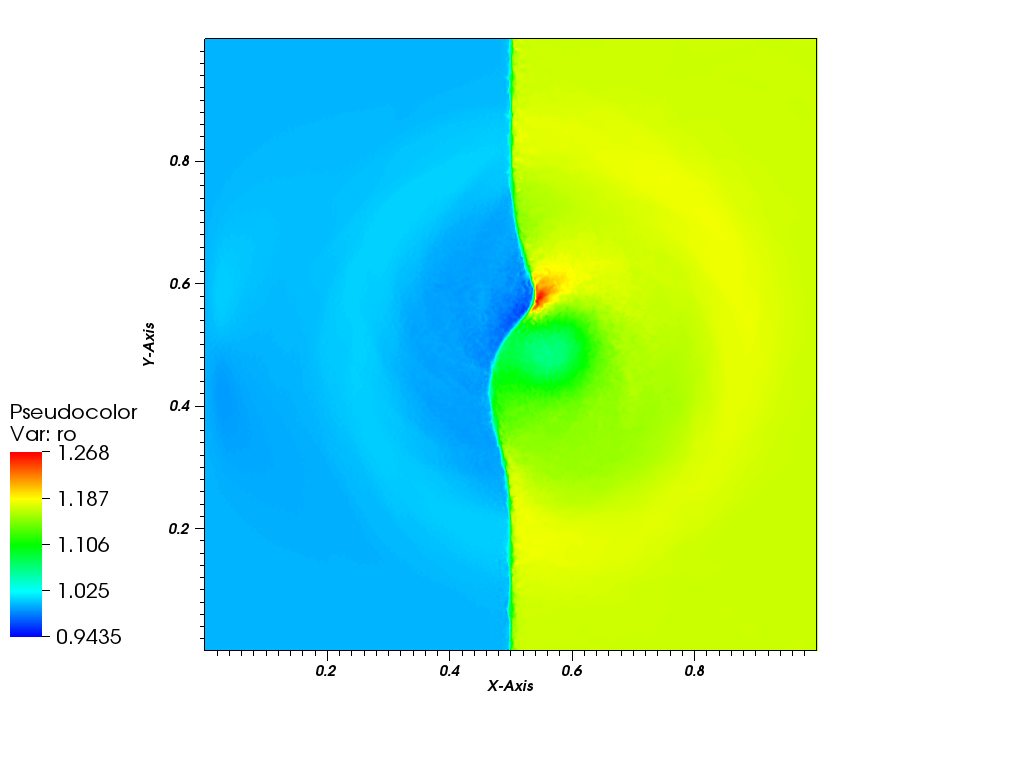}\label{shockVortex_B3_density}}
\caption{Shock Vortex Interaction. Mesh representation and results for the density at $T=0.25$ for different orders.}\label{shockVortex_B123}
\end{center}
\end{figure}

Finally, in Figure  \ref{shockVortex_B123} the density at the final time step is shown for $\mathcal{B}^2$ and $\mathcal{B}^3$. Here it is clearly visible how the solution increases in accuracy by displaying sharper contours when increasing the order.  

\subsubsection{2D Sod Problem - Structured Vs. Unstructured Meshes}\label{testSod2D}
Further, we have tested our high order RD scheme on a well-known 2D Sod benchmark problem. The initial conditions are given by
\begin{equation*}
(\rho_0,u_0,v_0,p_0) = 
\begin{cases}
[1, 0, 0, 1], \; &0 \leq r \leq 0.5, \\
[0.125, 0, 0, 0.1], \; &0.5 < r \leq 1,
\end{cases}
\end{equation*}
where $r = \sqrt{x^2+y^2}$ is the distance of the point $(x,y)$ from the origin. The stabilizing parameters have been set for $\mathcal{B}^1$, $\mathcal{B}^2$ to $\theta_1=0.01$ and $\theta_2=0$ and for $\mathcal{B}^3$ $\theta_1=0.02$ and $\theta_2=0$.

The idea behind this first test in 2D is to consider initially a structured mesh, which represents a straightforward extension of the 1D.
In particular, we show the results obtained with such a mesh for an arbitrary order of accuracy, such as $\mathcal{B}^2$ and then, successively, compare the results obtained for an unstructured mesh on the same test case.

\begin{figure}[H]
\centering
\subfigure[Structured grid, $N_3=16896$ elements]{\includegraphics[width=0.35\textwidth]{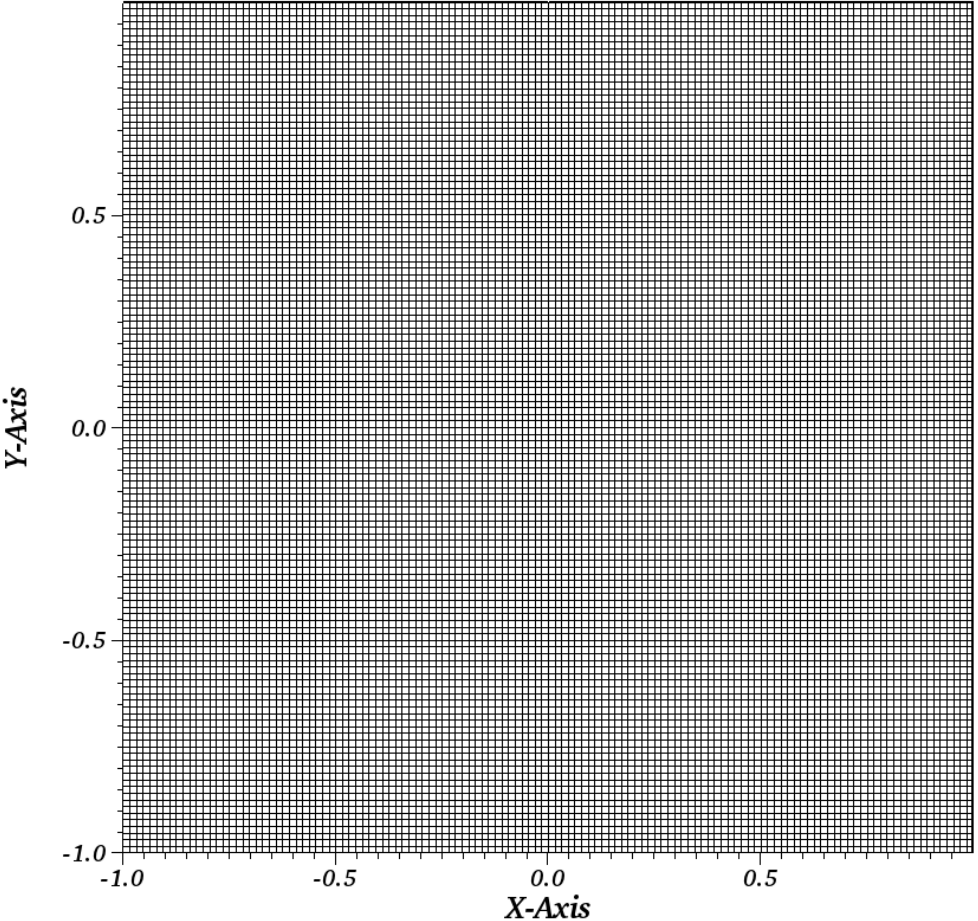}}
\hspace{0.4cm}\subfigure[Flux Indicator]{\includegraphics[width=0.35\textwidth]{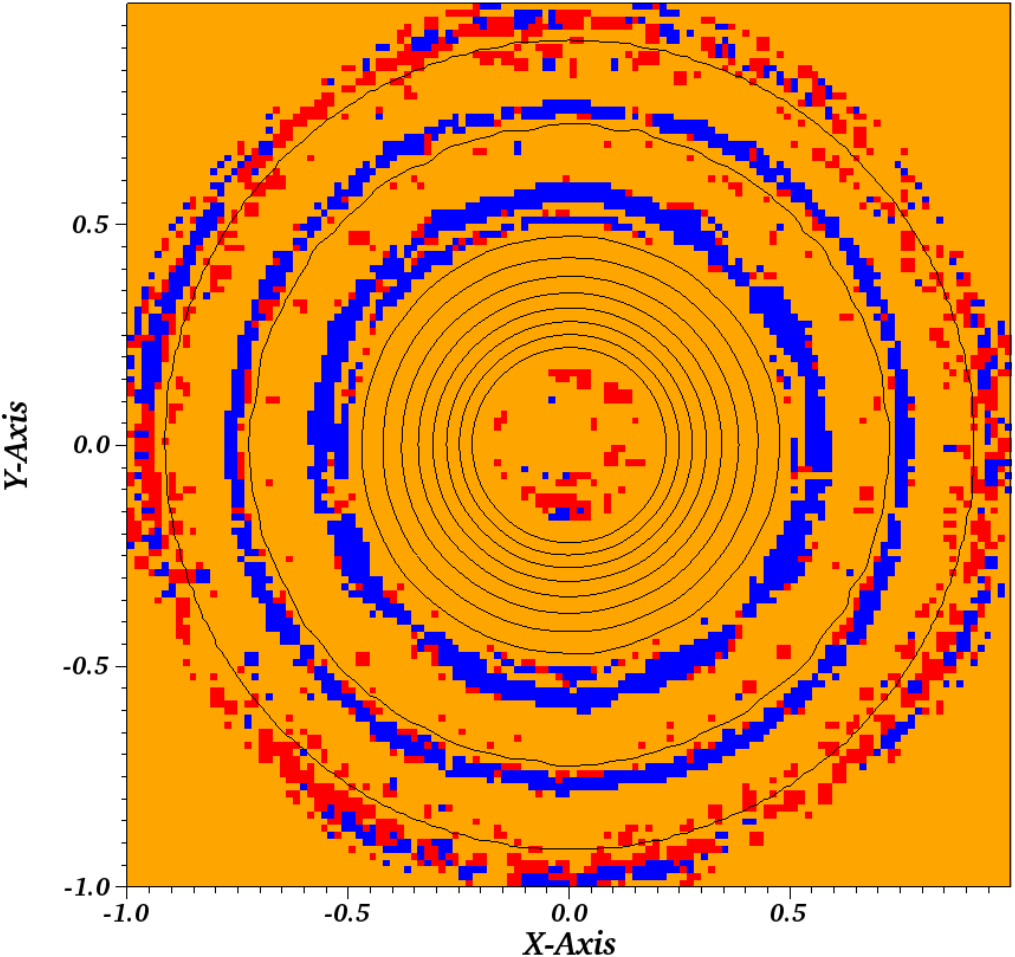}\label{Fig2_stucfluxi}}
\caption{2D Sod. Mesh (left) and flux indicator (right) for $\mathcal{B}^2$ at final time.
Orange corresponds to the GPJ scheme ($s=2$), red to RPJ ($s=1$) and blue to the Rusanov ($s=0$) scheme; Contour lines for the solution of the density are in black.}
\label{Sod_struct_fluxindic_Q2_ref3}
\end{figure}

\begin{figure}[H]
\centering
\subfigure[MOOD]{\includegraphics[width=0.45\textwidth]{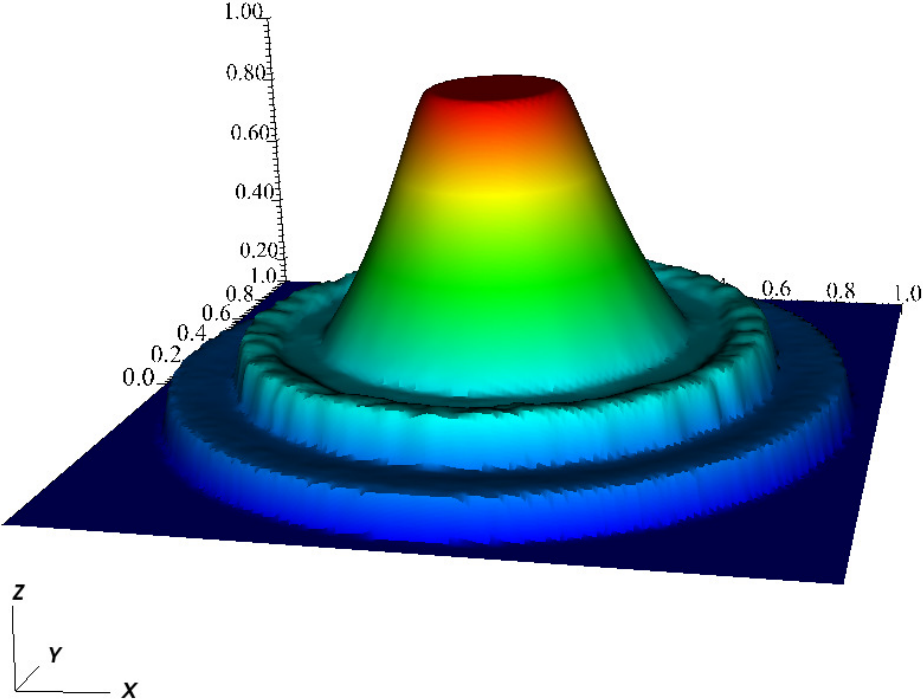}}
\subfigure[no MOOD]{\includegraphics[width=0.45\textwidth]{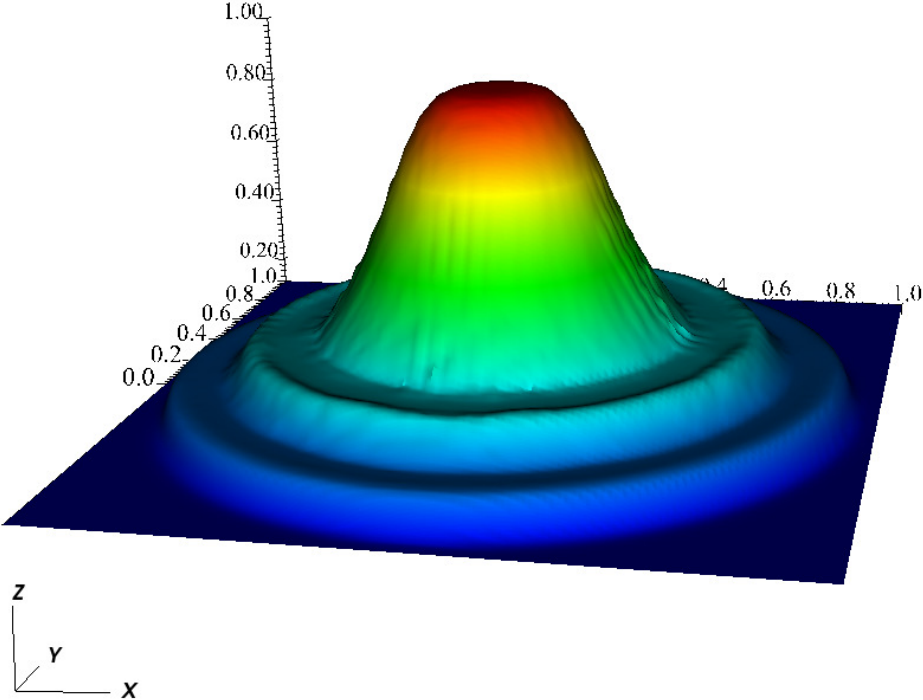}}
\caption{2D Sod. Results for the density on a structured grid with $N_3=16896$ elements and third order of accuracy.
Comparison between the ``MOOD'' strategy (left) and the ``N
no MOOD'' strategy (right).}
\label{Sod_struct_comparison_Q2_ref3_MoodNoMood}
\end{figure}

\begin{figure}[H]
\centering
\includegraphics[width=0.45\textwidth]{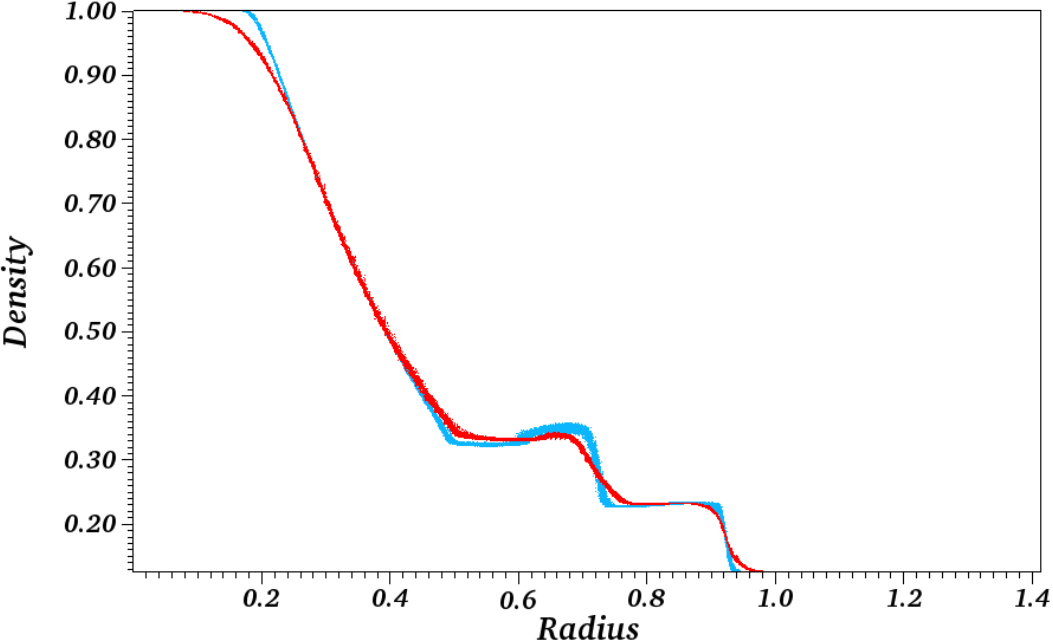}
\caption{2D Sod. Scatter plot of the density w.r.t. the radius on a structured grid with $N_3=16896$ elements and third order of accuracy.
Comparison between the ``MOOD'' strategy (light blue) w.r.t. the ``no MOOD'' approximation (red).}
\label{Sod_struct_scatter_Q2_ref3_MoodNoMood}
\end{figure}

\begin{figure}[H]
\centering
\includegraphics[width=0.45\textwidth]{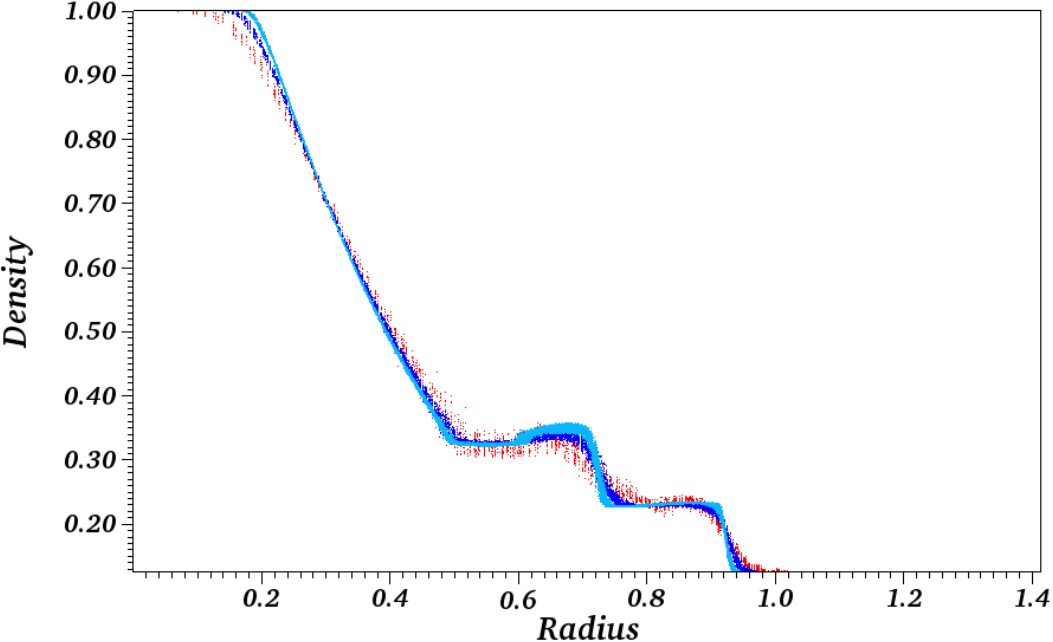}
\caption{2D Sod. Scatter plot of the density w.r.t. the radius on a structured grid.
Comparison of the ``MOOD'' strategy on different grid sizes, with $N_1=1152$ (red), $N_2=4352$ (blue) and
$N_3=16896$ (light blue) elements.}
\label{Sod_struct_scatter_Q2_refs_Mood}
\end{figure}

\begin{figure}[H]
\centering
\includegraphics[width=0.45\textwidth]{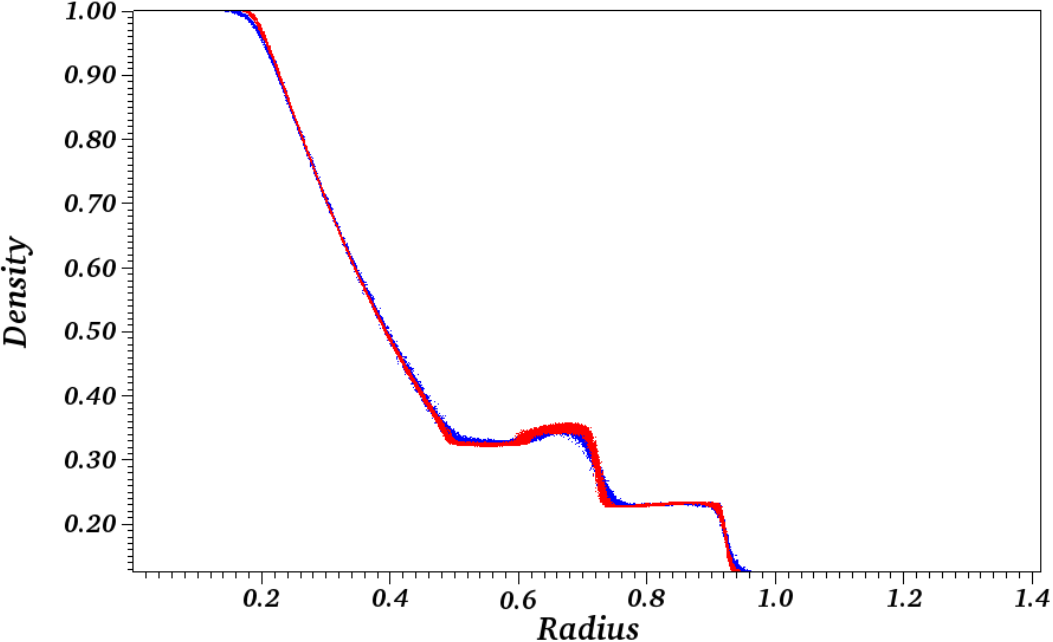}
\caption{2D Sod. Scatter plot of the density w.r.t. the radius for the ``MOOD'' strategy. Comparison of the scatter for different mesh types: structured with $16896$ elements (red) and unstructured with $13548$ elements (blue).}
\label{Sod_strucvsunsru_scatter_Q2_refs_Mood}
\end{figure}

As such, let us start by considering the detection procedure for 2D. To this extent we take a structured fine grid, with $N=16896$ elements. 
As shown in Figure \ref{Fig2_stucfluxi}, at the final time step, the contour lines (black lines) perfectly match the areas  where there is a flux change. This allows to see, that the detection criteria work effectively also in this case, and the mostly applied scheme throughout the computations is given, as expected, by the stabilized Galerkin (GPJ) approach.
Comparing in Figure \ref{Sod_struct_comparison_Q2_ref3_MoodNoMood} the obtained results for the density for the proposed scheme (denoted as ``MOOD'') and an approach without the ``a posteriori'' limiting strategy (i.e. ``no MOOD''), one can clearly see, that more structures appear to be outlined and, thus, allows for less dissipation.
Note, that the ``no MOOD'' approach corresponds to the sole RPJ scheme proposed in this chapter.
The improvement in accuracy, can also be seen in Figure \ref{Sod_struct_scatter_Q2_ref3_MoodNoMood}, where we have compared the density scatter plots for the ``MOOD'' and ``no MOOD'' approaches.
Furthermore, to guarantee mesh convergence, we have compared in Figure \ref{Sod_struct_scatter_Q2_refs_Mood} solutions obtained for different mesh refinements, i.e. $N_1=1152$,
$N_2=4352$ and
$N_3=16896$. 
As already observed for the 1D case, even for extremely coarse meshes the approximation results to be of high quality.
In Figure \ref{Sod_strucvsunsru_scatter_Q2_refs_Mood}, the scatter for the density shows an excellent overlap for the solutions obtained with a ``MOOD'' strategy for both structured and unstructured meshes, and, again, as in the previous tests in 2D, for $\mathcal{B}^2$. 
The grid size is in this test is of comparable order, as the one for the structured grid counts $16896$ elements, while the unstructured about $13548$. 

Having demonstrated the absence of any relevant difference between the structured and unstructured approach, we consider in the following Figure \ref{Test2D_unstruc_mesh_comparison_orders} the comparison between different order of accuracy. In particular, the focus is set on the comparison between the ``MOOD'' and ``no MOOD'' approaches, for a coarsed ($N_1=3576$ elements) and finer ($N_2=13548$ elements) unstructured mesh. One can note that more scatter is observed in  Figure \ref{Test2D_unstruc_mesh_comparison_orders}  for $\mathcal{B}^3$ compared to $\mathcal{B}^2$. This is due to the choice to plot all the degrees of freedom of each cell and thereforebeing there more points for higher order scheme the appearance is more dense and slightly wider.

\begin{figure}[H]
\centering
\subfigure[ $\mathcal{B}^1$ for $N=3576$]{\includegraphics[width=0.45\textwidth]{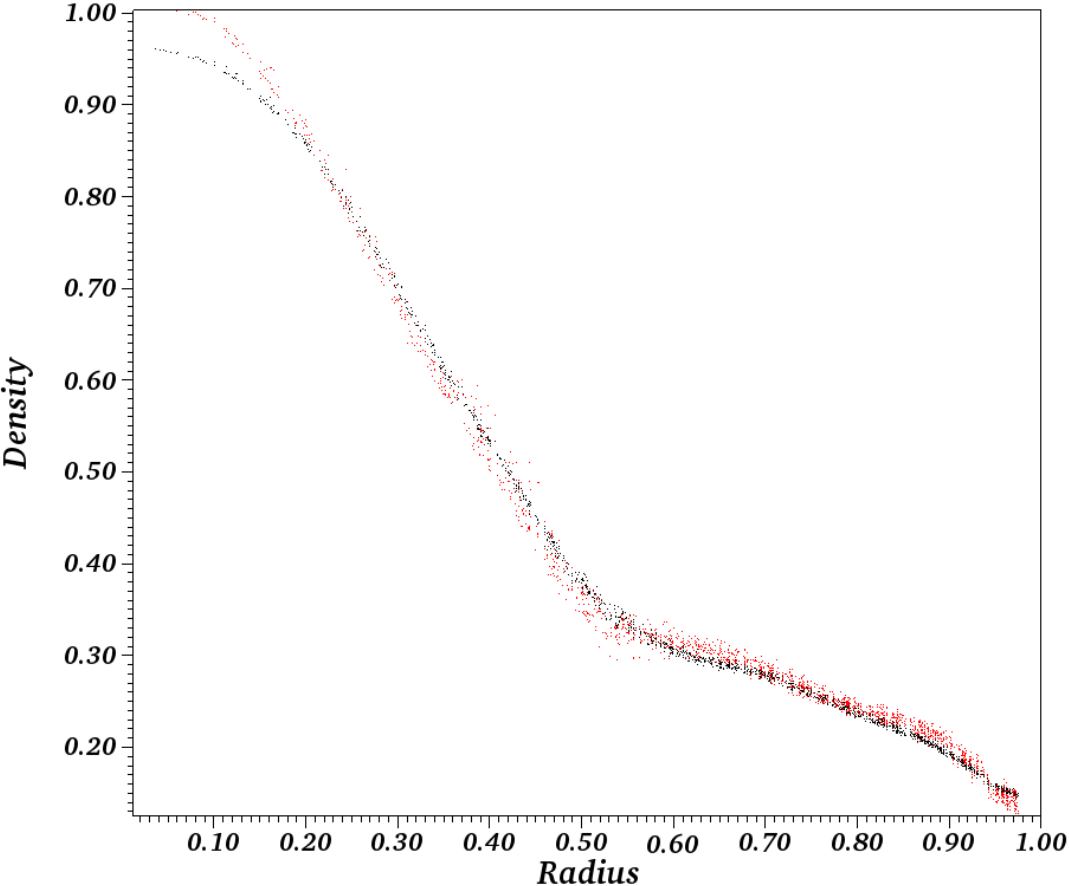}}
\subfigure[ $\mathcal{B}^1$ for $N=13548$]{\includegraphics[width=0.45\textwidth]{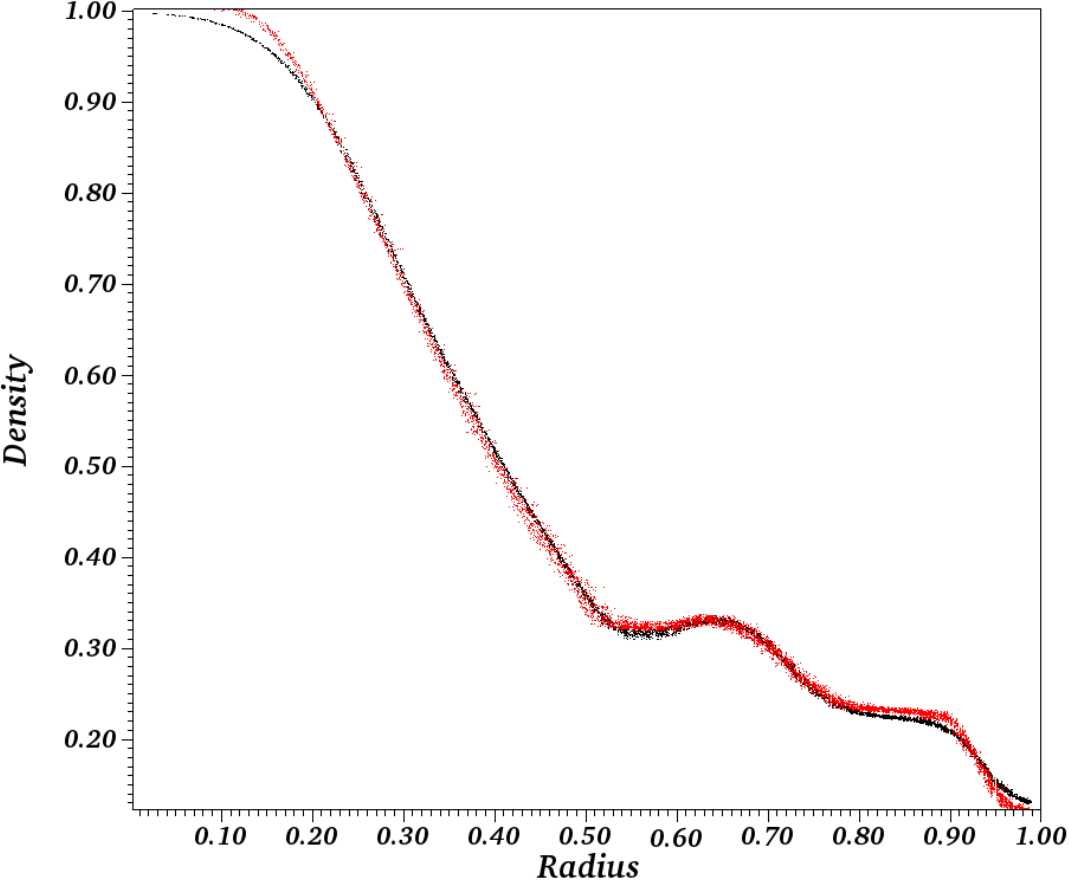}}\\
\subfigure[ $\mathcal{B}^2$ for $N=3576$]{\includegraphics[width=0.45\textwidth]{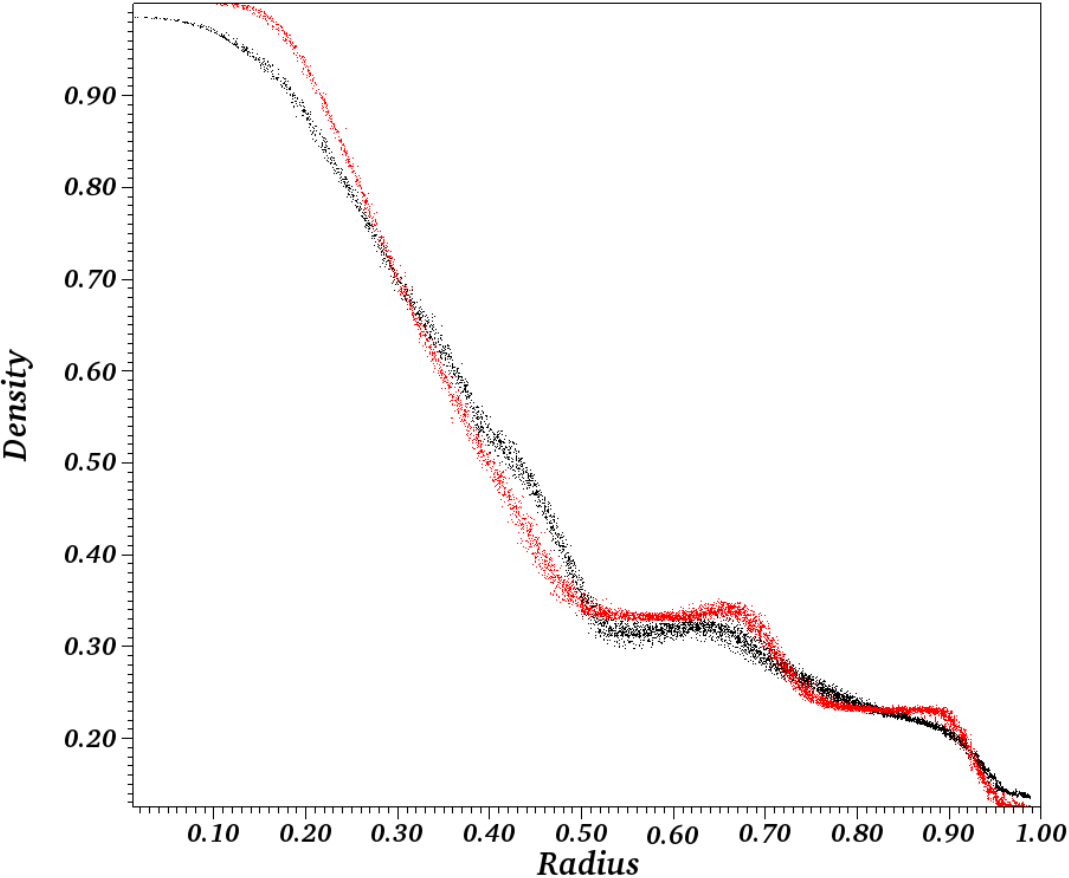}}
\subfigure[ $\mathcal{B}^2$ for $N=13548$]{\includegraphics[width=0.45\textwidth]{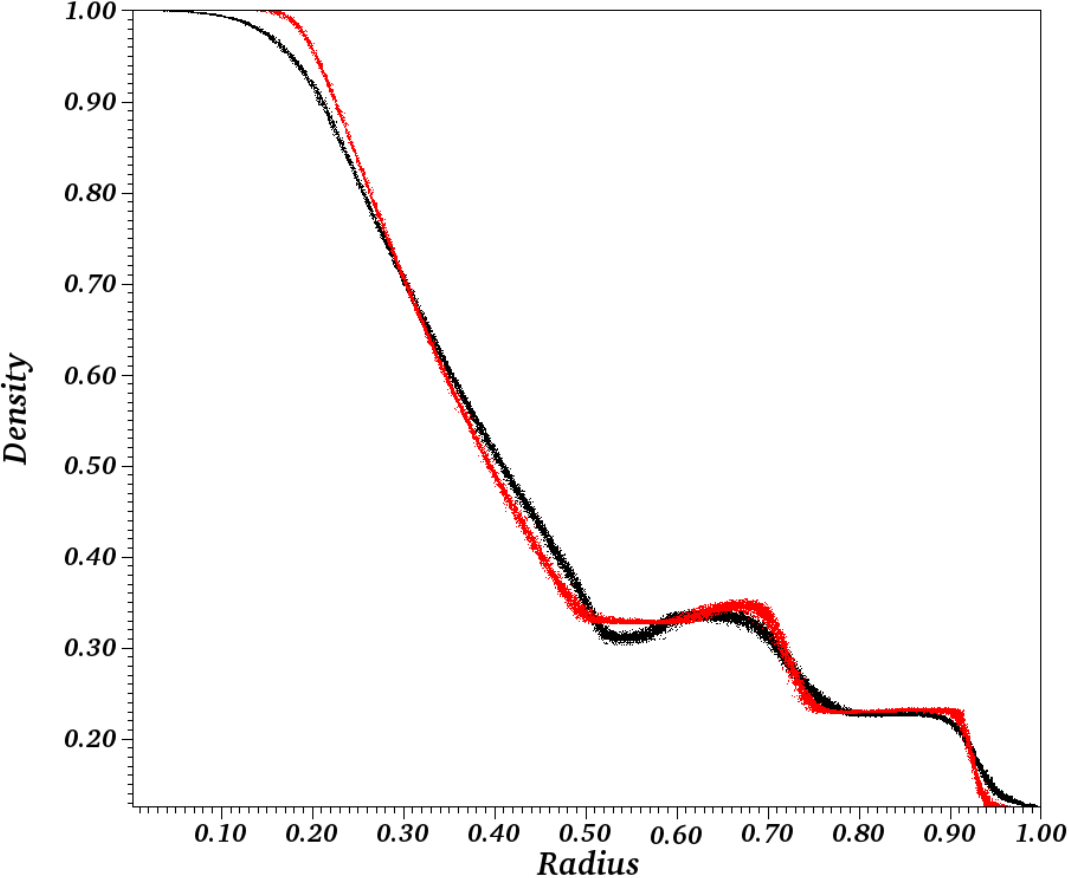}}\\
\subfigure[ $\mathcal{B}^3$ for $N=3576$]{\includegraphics[width=0.45\textwidth]{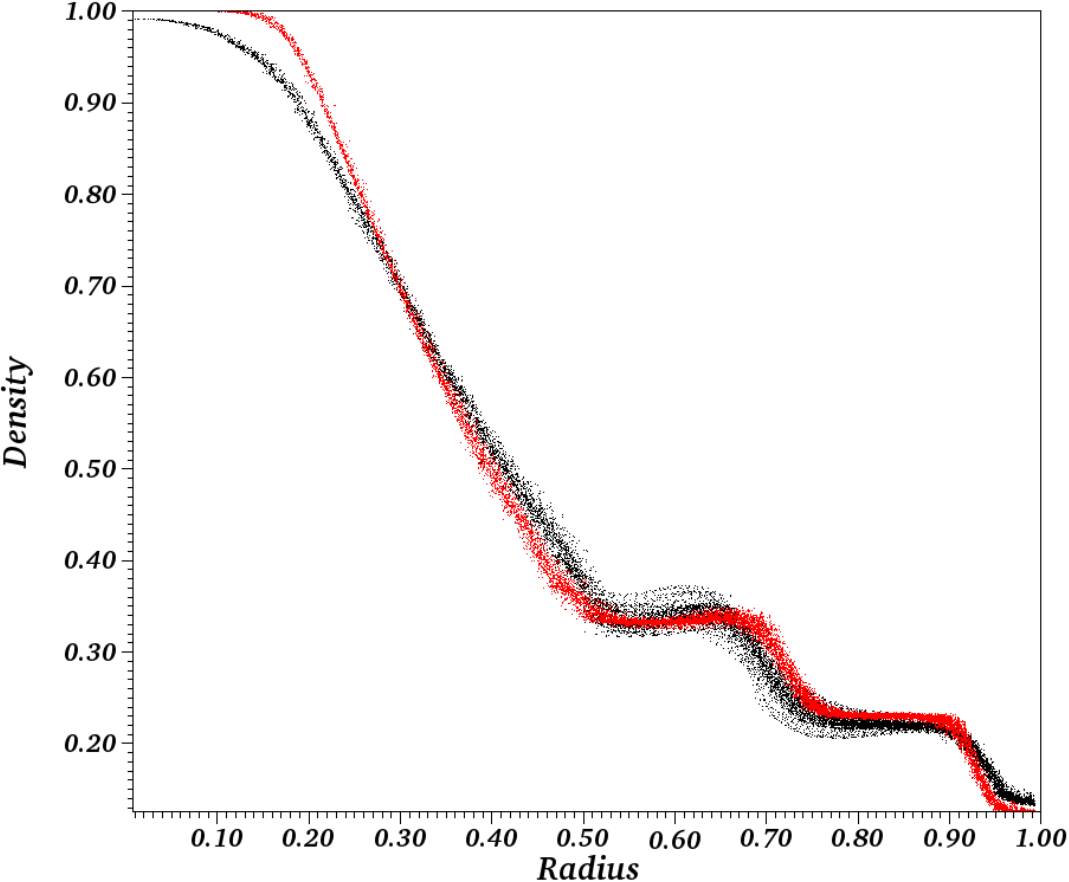}}
\subfigure[ $\mathcal{B}^3$ for $N=13548$]{\includegraphics[width=0.45\textwidth]{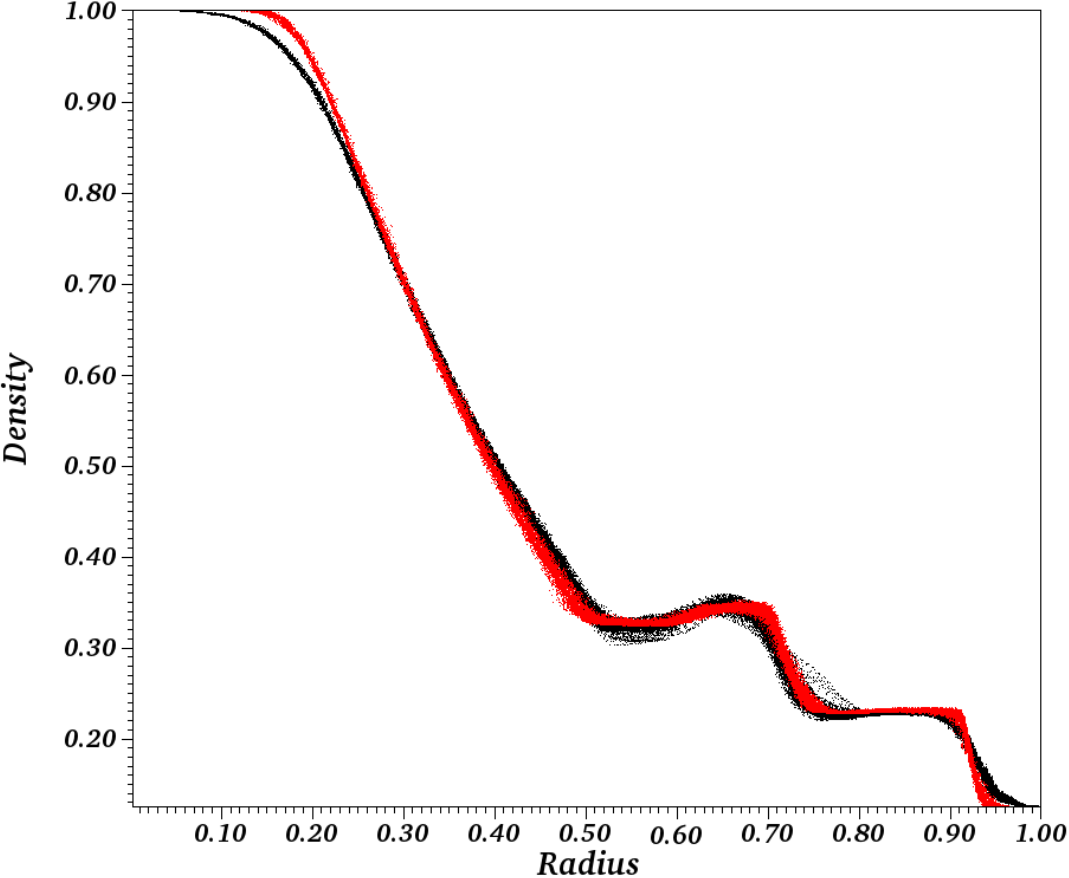}}
\caption{2D Sod. Density scatter plot w.r.t. the radius at the final time step for unstructured grids. Comparison between the ``MOOD (red) and ``no MOOD'' (black) strategy, for coarse (left) and fine (right) meshes.}
\label{Test2D_unstruc_mesh_comparison_orders}
\end{figure}

\newpage
\subsubsection{Mach 3 channel with forward-facing step}
\begin{figure}[H]
\begin{center}
\hspace{-0.9cm} 
\subfigure[$\mathcal{B}^2$ for $N_1$ - no MOOD]{\includegraphics[width=0.475\textwidth]{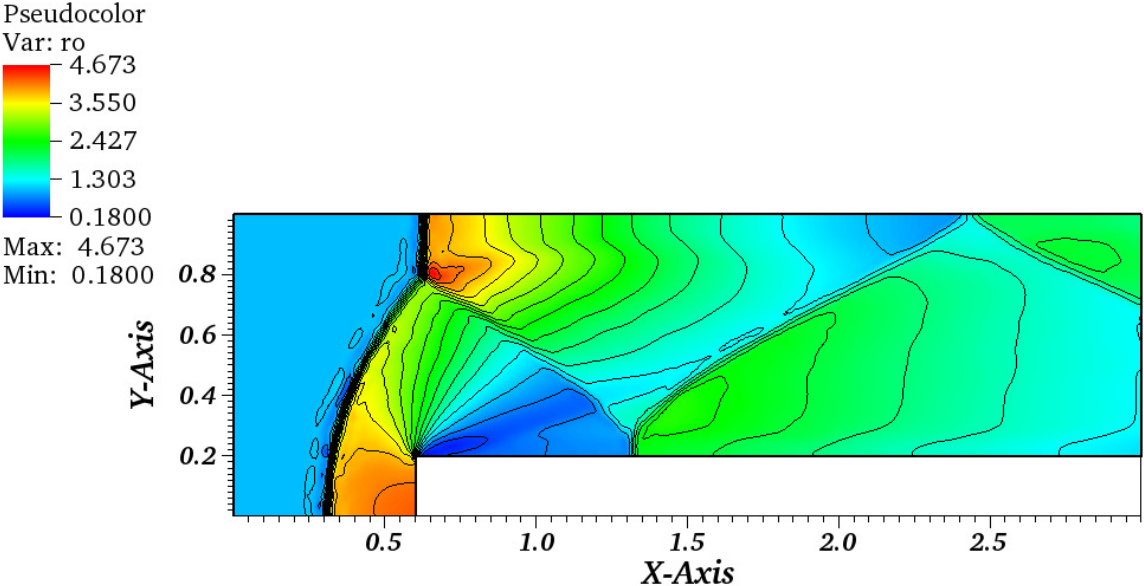}}
\subfigure[$\mathcal{B}^2$ for $N_1$ -  MOOD]{\includegraphics[width=0.47\textwidth]{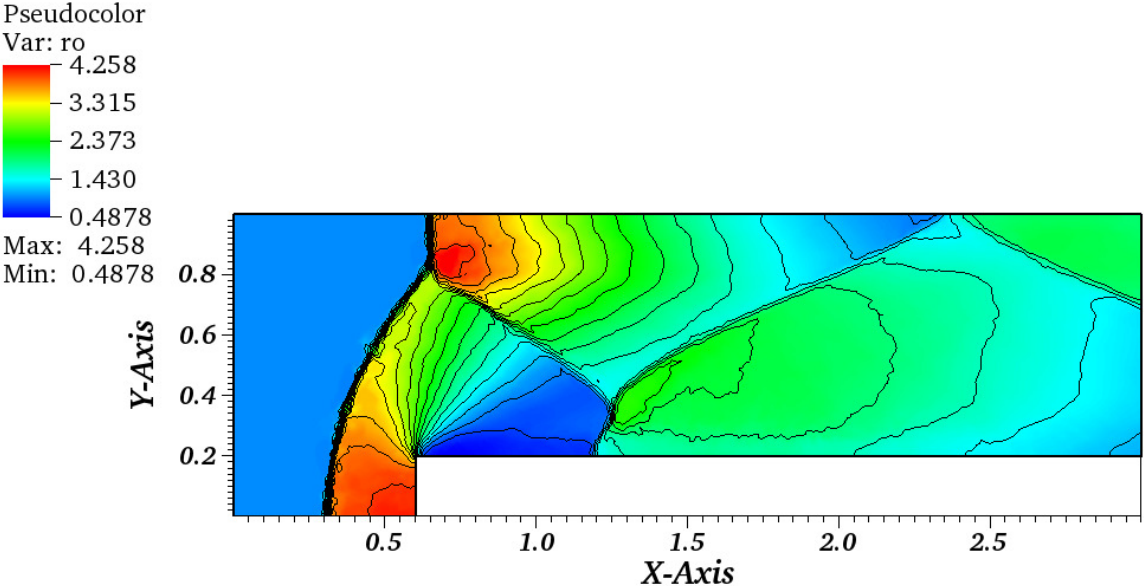}}
\caption{Mach 3 channel with step. Results at $T=4.0$. Comparison on a fine mesh for the approximations obtained for  non-MOOD (left) and MOOD (right) schemes.}\label{step_2D_2}
\end{center}
\end{figure}

To assess the robustness of the proposed scheme in multidimensional problems involving strong shock waves, the Mach $3$ channel with a forward-facing step \cite{Woodward1984} test case has been used with $\mathcal{B}^2$ 
elements on a mesh having  $N=11072$ cells\footnote{corresponds roughly to $60\times 200$ grid points}) (see Fig.~\ref{step_2D_2}). The stabilizing parameters have been set as in Section \ref{testSod2D}. 
In particular, we have compared the solution of the MOOD approach, with respect to a non-MOOD approach. It is possible to observe the gain in the quality of the approximation of shock waves when using the a posteriori limiting approach with respect to an a priori one. 
Note that no effort has been made to address the behaviour of the entropy in the expansion fan at the corner, this results in an anomalous behaviour of the reflected shock wave: there should be no lambda shock. Our purpose is not to address this issue.

\subsubsection{Double Mach Reflection problem}
\begin{figure}[H]
\begin{center}

\subfigure[$\mathcal{B}^2$ no MOOD]{\includegraphics[width=0.395\textwidth]{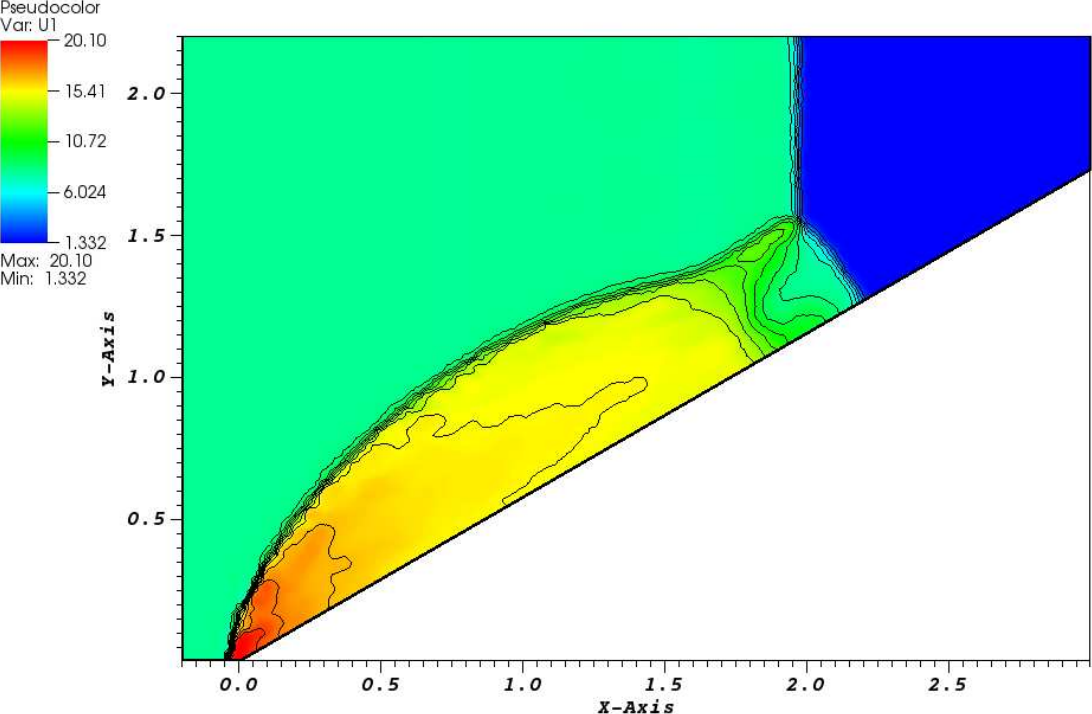}}
\hspace{0.6cm}
\vspace{0.4cm}\subfigure[$\mathcal{B}^2$ MOOD]{\includegraphics[width=0.41\textwidth]{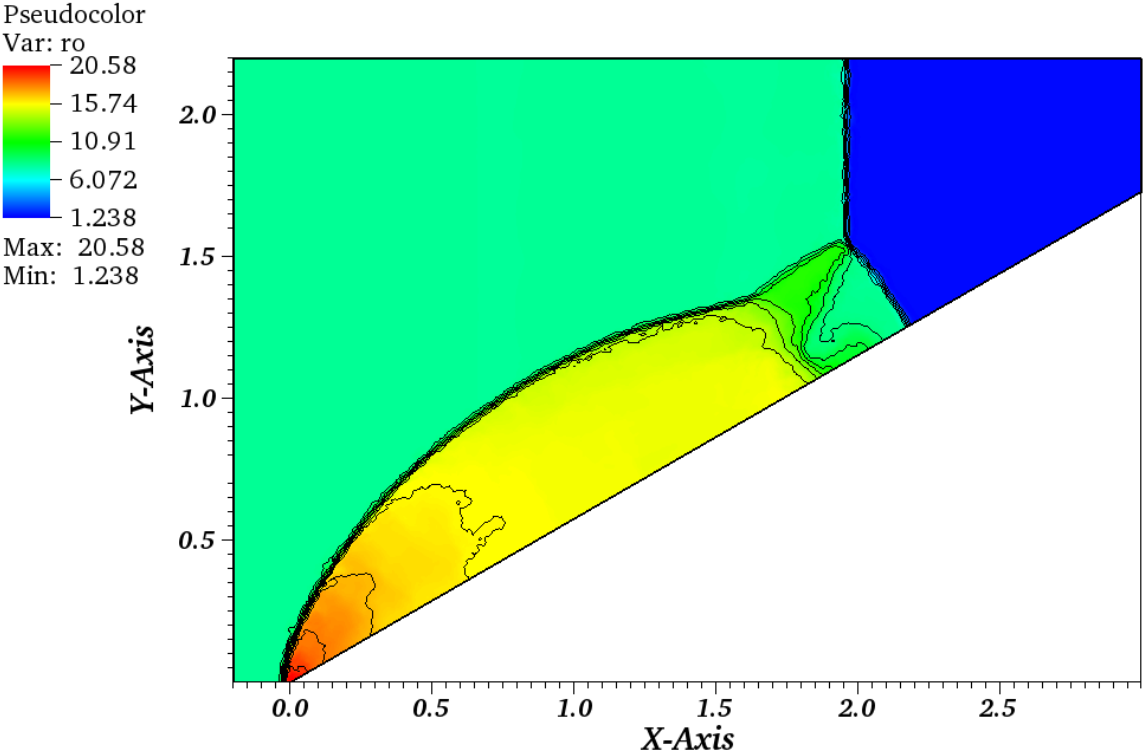}}\\
\caption{Double Mach Reflection problem. Results at $T=0.2$}
\label{Fig:DMR}
\end{center}
\end{figure}

Finally, we present a widely used benchmark problem of a double Mack reflection problem as described in \cite{Woodward1984}. In this case, $\mathcal{B}^2$ elements have been computed on a mesh having $N=19248$ cells\footnote{corresponds roughly to $60\times 200$ grid points}) (see Fig.~\ref{Fig:DMR}). The stabilizing parameters have been set as in Section \ref{testSod2D}. Also here, as expected, the quality of the solution increases when an a posteriori limiting approach is chosen. 



\section{Conclusion and Perspectives} \label{sec:conclusion}

A novel explicit high order residual distribution scheme with an ``a posteriori'' blending strategy has been designed in the context of the Euler equations in gas dynamics. 
Our strategy has been to approximate with a high order scheme in time and a least dissipative scheme in space a solution displaying strong interacting discontinuies, while guaranteeing through the design of  an ``a posteriori'' limiting strategy the detection of numerical/computational misbehaviour. This allows to locally intervene on affected areas with a more dissipative scheme for the spatial discretization.
The considered benchmark problems have thoroughly validated the proposed methodology, assessing for its capability to provide an accurate and robust numerical method in 1D and 2D.
For smooth solutions tested via the isentropic flow benchmark problem, we have seen that the effective optimal accuracy is attained, while 
for non smooth flows, as in case of the considered shock tubes, double Mach reflection and forward facing step, we have observe that the proposed method is able to provide non oscillatory and accurate numerical solutions. 
Moreover, we have tested two different strategies found in literature, which consider two approaches to distinguish natural oscillations from numerical ones, and have established that one strategy provides higher quality results within this novel designed approximation strategy of residual distribution schemes.
We have further shown, how the number of troubled cells detected by the considered detection methods 
is monitored, with the aim to display the feature that not  many cells are flagged, and, as such,  the number of cascade iterations is limited, and as such the extra-cost.
 Overall, the ``a posteriori'' treatment renders the high order residual distribution scheme robust and positivity preserving  when notoriously difficult tests are simulated. \\
Extensions to other models, such as multiphase flows or Lagrangian hydrodynamics, and further investigations of high order residual distribution schemes will be considered in  forthcoming papers. Along this line, we are currently extending the proposed approach to viscous problems by combining it with the discretisation technique as explained in \cite{DeSantis2015}. This last extension might require some modification in the time-stepping as the time step would be very small, and thus an implicit approach is needed with the challenge to have, nevertheless, a diagonal 'mass matrix'.

%
%

%
%
%
 
%



\section*{Acknowledgments}
P.B. has been funded by SNSF project 200021\_153604 ``High fidelity simulation for compressible materials''.  R.A. has been funded in part by the same project. The authors would like to thank Rapha{\"e}l Loub{\`e}re (CNRS and Universit{\'e} de Bordeaux, France) for the very helpful discussions about MOOD and Fran\c{c}ois Vilar (Universit{\'e} de Montpellier, France) for the helpful discussions for the LSE detection strategy of this paper.


\bibliographystyle{plain}
\bibliography{biblio}

\appendix
\section{Positivity property of the Rusanov residuals}
\label{positivity}

The question whether the designed Residual Distribution scheme with the ``a posteriori'' limiting guarantees, in case of the Rusanov parachute scheme for the spatial discretization, the preservation of positivity is, formally, still an open issue and will be the topic of a forthcoming work.
Even though a formal proof has not been carried out, the performed numerical experiments have widely confirmed this propriety (see cf. Section \ref{sec:numerics}).

The first order scheme, using Rusanov residuals, writes whatever the interpretation of the degrees of freedom, as 
\begin{equation}
\label{remi:1}
\U_\sigma^{n+1}=\U_{\sigma}^n-\dfrac{\Delta t}{|S_\sigma|}\sum_{K, \sigma\in K} \phi_{\sigma,\mathbf{x}}^{K,Rus}(\U^n).
\end{equation}
Using definition \eqref{S_sigma}, $\U_\sigma^{n+1}$ is rewritten as, setting $ |K_\sigma|=\dfrac{|K|}{N_{DoF}}$ as
\begin{equation}
\label{remi:2}
\U_\sigma^{n+1}=\sum_{K, \sigma\in K} \dfrac{|K_\sigma|}{|S_\sigma|} \U_\sigma^{K,\star}
\text{ with }\U_\sigma^{K,\star}=\U_\sigma^n-\dfrac{\Delta t}{|K_\sigma|}\phi_{\sigma,\mathbf{x}}^{K,Rus}(\U^n)
\end{equation}
Let us note, that the Rusanov's residual can be interpreted in two different ways
\begin{itemize}
\item Version 1:
$$\phi_{\sigma,\mathbf{x}}^{K,Rus}=\int_K \varphi_\sigma \text{ div } \F(\U) \;d\xx+\alpha_K\big (\U_\sigma-\overline{\U}\big )$$
\item Version 2:
$$\phi_{\sigma,\mathbf{x}}^{K,Rus}=\dfrac{1}{N_{DoF}}\int_{K}\text{ div } \F(\U) \;d\xx+\alpha_K\big (\U_\sigma-\overline{\U}\big )$$
\end{itemize}
The $\overline{\U}$ is the arithmetic average of the $\U_\sigma$, as specified in the previous subsection. The second variation is about the nature of the degrees of freedom: do we use Lagrange of Bernstein approximation.

The purpose is to estimate a minimal value of $\alpha_K$ which guarantees, for the compressible Euler system, that if the densities and pressure are positive at $t_n$, they will stay positive at the next time. 

In the following, we will first consider the case of the Lagrange interpolation. Then,  we will extend  it to the case of Bernstein approximation, which is our true target in this paper.

\subsection{Set of thermodynamical states}
The aim is to have positive densities $\rho$ and positive internal energies $e= E-\frac{1}{2}\rho u^2$. 
This set is convex under standard assumptions on the thermodynamics variables which hold for the equations of state for standard perfect and stiffened gas.

One can observe that in case of Lagrangian polynomials, where one takes the nodal values, the positivity of the density and of the internal energy are straightforward.

In case of Bernstein polynomials, the analysis is a bit more involved because the degree of freedom do not correspond in general to point values.
The first remark is that 
\begin{equation}
\begin{split}
\mathcal{K}_{th}'=\{ (\rho_\sigma, m_\sigma, E_\sigma)_{\sigma\in K} \text{ s. t. } & \rho=\sum_{\sigma\in K} \rho_\sigma B_\sigma\geq 0 \text{ on } K
\text{ and }  E-\frac{1}{2}\frac{m^2}{\rho}\geq 0,\\& \text{ with } E=\sum_{\sigma} E_\sigma B_\sigma \text{ and } m=\sum_{\sigma} m_\sigma B_\sigma,\},
\end{split}
\label{set_convex_Bernstein}
\end{equation}
where we have denoted the momentum by $m=\rho\,u$ is convex. We note also that instead testing the inequalities for all $x\in K$, one can test them only for a finite set of points, for example the Lagrange points, and the resulting set that we still denote by $\mathcal{K}_{th}'$ is also convex.

The proof is as follows.
\begin{proof}
If the functions $\rho=\sum\limits_{\sigma} \rho_\sigma B_\sigma$ and $\rho'=\sum\limits_{\sigma'} \rho_{\sigma'} B_\sigma$ defined similarly are positive, and hence for any $\lambda\in [0,1]$,  the densities defined from $\U=\lambda \U+(1-\lambda)\U'$ are positive on the simplex $K$.
The internal energy is a rational function of the conserved quantities.
one recasts a rational function, i.e. there is a division between two polynomials, and as such, some further considerations need to be done.
Let us consider the mapping $\varphi_e: (\rho, m, E)\mapsto E-\dfrac{1}{2}\dfrac{m^2}{\rho}.$
The internal energy $E-\dfrac{1}{2}\dfrac{m^2}{\rho}$ is a concave function of $(\rho,m,E)$, as its Hessian is 
$$\begin{pmatrix}
-\dfrac{m^2}{\rho^3} & \dfrac{m}{\rho^2}&0\\
\dfrac{m}{\rho^2}&-\dfrac{1}{\rho}&0\\
0&0&0
\end{pmatrix},
$$
with eigenvalues $0$ (twice) and $-\frac{m^2}{\rho^3}-\frac{1}{\rho}<0$ if $\rho\geq 0$.

Hence, if $U$ and $U'$ belong to $\mathcal{K}_{th}'$, and $\lambda\in [0,1]$, the density function associated to $\lambda U+(1-\lambda)U'$ will be positive, and the internal energy is $\varphi_e(\lambda U+(1-\lambda)U')$. Since $\varphi_e$ is concave,
$$\varphi_e(\lambda U+(1-\lambda)U')\geq \lambda\varphi_e( U)+(1-\lambda)\varphi_e( U')\geq 0,$$
and hence $\lambda U+(1-\lambda)U'\in \mathcal{K}_{th}'$. This shows that this set is also concave
\end{proof}

However, it is difficult to characterize $\mathcal{K}_{th}'$ and we introduce a stronger condition.

We consider $\overline{\mathcal{K}}'_{th}$ the set
$$ \overline{\mathcal{K}}'_{th}=\bigg \{ \text{for all DoF }\sigma, (\rho_\sigma, m_\sigma, E_\sigma ), \rho_\sigma\geq 0, E_\sigma-\frac{1}{2}\dfrac{m_\sigma^2}{\rho_\sigma}\geq 0\bigg \},$$
and we notice that $$\overline{\mathcal{K}}'_{th}\subset \mathcal{K}'_{th}.$$

\begin{proof}
Thanks to the positivity of the Bernstein polynomials, and the Cauchy-Schwarz inequality, we have
$$\bigg ( \sum\limits_{\sigma\in K} m_\sigma B_\sigma\bigg )^2 \leq \bigg ( \sum\limits_{\sigma\in K}\rho_\sigma B_\sigma\bigg ) \bigg ( \sum\limits_{\sigma\in K}\dfrac{m_\sigma^2}{\rho_\sigma} B_\sigma\bigg ),$$
and then
$$\sum\limits_{\sigma\in K} E_\sigma B_\sigma -\frac{1}{2}\dfrac{\bigg ( \sum\limits_{\sigma\in K} m_\sigma B_\sigma\bigg )^2}{\sum\limits_{\sigma\in K}\rho_\sigma B_\sigma}\geq \sum\limits_{\sigma\in K} \bigg ( E_\sigma-\frac{1}{2}\dfrac{m_\sigma^2}{\rho_\sigma}\bigg ) B_\sigma\geq 0$$
if all the $(\rho_\sigma, m_\sigma, E_\sigma)\in \overline{\mathcal{K}}'_{th}.$
\end{proof}

\subsection{Case of Lagrange interpolation}
Let us consider the Lagrange interpolation, so that $\U_\sigma=\U(\sigma)$ is the evaluation of the solution at the Lagrange degree of freedom. Remind that they are defined in a simplex by their barycentric coordinates which are, for the degree $n$ and for triangles, $(\frac{i_1}{n+1},\frac{i_2}{n+1}, \frac{i_3}{n+1})$ with $i_1+i_2+i_3=n+1$. In the 3D case, they would be defined in the same way. For quad or hex, they are simply obtained by tensorisation of the 1D Lagrange points.

\subsubsection{The one-dimensional case}
Rephrasing the proof of Perthame and Shu \cite{PerthameShu}, we have
\begin{equation}\label{Proofposi_lxf}
\U_i^{n+1}=\U_i^n-\frac{\Delta t}{\Delta x} \left[ \widehat{\F}(\U_{i+1},\U_i)-\widehat{\F}(\U_i,\U_{i-1})\right].
\end{equation}
Denoting by $K$ the interval $[x_i, x_{i+1}]$ and introducing the splitting in the equation
$$\U_t+\mathbf{F}(\U)_x=0$$
by
\begin{equation}
\label{remi:split:1}
\U_t+\big (\F(\U)+\eta\U\big )_x=0, ~\text{and}~ \U_t+\big (\F(\U)-\eta \U\big )_x=0,
\end{equation}
we see that if $\eta=\max\limits_{x\in K} ||u(x)||+c(x)$ (where $u$ is the velocity and $c$ is the sound speed), the Rusanov scheme is recast as combination of the Godunov scheme and the downwind scheme, i.e. the left and right equations in \eqref{remi:split:1} accordingly.
Hence the value $\U_i^{n+1}$ can be interpreted as the average of 
$$
\widetilde{\U}=\U_i^n-\frac{\Delta t}{\Delta x} \left[ \left( {\F}(\U_i^n)+\eta\U_i^n \right)-\left( \F(\U_{i-1}^n)+\eta\U_{i-1}^n \right) \right]$$
and
$$
\widetilde{ \widetilde{\U} }= \U_i^n-\frac{\Delta t}{\Delta x} \left[ \left({\F}(\U_{i+1}^n )-\eta\U_{i+1}^n \right)-\left( \F(\U_{i}^n) -\eta\U_{i}^n \right) \right].$$
If the states $\U_i^n$, $\U_{i+1}^n$ and $\U_{i-1}^n$ belong to the convex $\mathcal{K}_{th}$ defined as
$$\mathcal{K}_{th}= \left\{ (\rho, \rho u, E), \rho\geq 0, E-\frac{1}{2}\rho u^2\geq 0 \right\},$$
 then $\U_i^{n+1}$ will belong to the same convex $\mathcal{K}_{th}$.
 One can conclude that, for both the Lagrangian interpolation and the Bernstein reconstruction, the Rusanov scheme \eqref{Proofposi_lxf} for the one-dimensional case preserves the convex sets $\mathcal{K}_{th}$ for the Lagrange interpolation and $\mathcal{K}_{th}'$ for the Bernstein reconstruction.
 
\section{The multi-dimensional Case}\label{multi-D-appendix}
\begin{itemize}
\item Version 1:
The residuals can be written as:
\begin{equation}
\begin{split}
\phi_{\sigma,\mathbf{x}}^{K,Rus}&=\int_K\varphi_\sigma \; \text{ div }\F\; d\xx+\alpha \big ( \U_\sigma-\bar\U\big )
\\&=\sum_{\sigma'\in K} \left[ \left( \int_K \varphi_\sigma\nabla \varphi_{\sigma'}d\xx \right) \cdot \F_{\sigma'}+\frac{\alpha}{N_{DoF}} \left( \U_\sigma-\U_{\sigma'} \right) \right] \\ & =\sum_{\sigma'\in K} \left[  2\left( \int_K \varphi_\sigma\nabla\varphi_{\sigma'}d\xx \right) \cdot\dfrac{\F_\sigma+\F_{\sigma'}}{2}+\frac{\alpha_K}{N_{DoF}} \left( \U_\sigma-\U_{\sigma'} \right) \right]
\end{split}
\end{equation}
because $\sum\limits_{\sigma'\in K}\int_K \varphi_\sigma\nabla\varphi_{\sigma'}d\xx=\int_K \varphi_\sigma \nabla (1)\; d\xx=0$, 
so that \eqref{remi:2} writes:
\begin{equation*}
\begin{split}
\U_\sigma^\star&=\U_{\sigma}^n-\dfrac{\Delta t}{|K_\sigma|}\sum_{\sigma'\in K} \left[ \left( 2\int_K \varphi_\sigma\nabla\varphi_{\sigma'} d\xx \right) \cdot\dfrac{\F_\sigma+\F_{\sigma'}}{2}+\frac{\alpha}{N_{DoF}}  \left( \U_\sigma-\U_{\sigma'} \right) \right]\\
&=\dfrac{1}{N_{DoF}} \sum\limits_{\sigma'\in K}\left[ \U_{\sigma}^n-\dfrac{\Delta t}{|K_\sigma|} \;\omega_{\sigma\sigma'}\cdot \dfrac{\F_\sigma+\F_{\sigma'}}{2}+\alpha_K \left( \U_\sigma-\U_{\sigma'} \right) \right]
\end{split}
\end{equation*}
where $$\omega_{\sigma\sigma'}=2N_{DoF}\int_K \varphi_\sigma\nabla\varphi_{\sigma'}\; d\xx.$$
Then one can interpret the vector $\omega_{\sigma\sigma'}$ as a scaled normal\footnote{because $\int_K\varphi_\sigma=\frac{|K|}{N_{DoF}}$}, and the stability condition writes:
$$\alpha_K\geq \max_{\xx\in K} \rho(\mathbf{A(U(\xx))}\cdot\omega_{\sigma\sigma'})$$
where for any vector $\bn=(n_x,n_y)$, $\mathbf{A(U)}\cdot \bn=\dpar{F}{U}(U) n_x+\dpar{G}{U}(U) n_y$, and $F$ (resp $G$) is the $x$- (resp. $y$-) component of the flux $\F$.
\item Version 2:
The algebra is similar since 
$$\int_K \text{div }\F \; d\xx =\sum_{\sigma\in K} \left( \int_K\nabla \varphi_\sigma\; d\xx \right) \cdot \F_\sigma,$$
and we get a similar stability condition.
\end{itemize}

\subsection{Case of Bernstein approximation}
When using Bernstein degrees of freedom, we have,  at the Lagrange degrees of freedom denoted by $\sigma_L$ for this paragraph
$$U(\sigma_L)=\sum_{\sigma}U_\sigma B_\sigma(\sigma_L), \text{ with }B_\sigma(\sigma_L)\geq 0, \text{ and  }\sum_\sigma B_\sigma(\sigma_L)=1.$$
This means that there is a linear mapping
$$M=\bigg(B_\sigma(\sigma_L)\bigg )$$ that maps the Bernstein Dof to the Lagrange ones.

If we have 
$$\U_\sigma^{K,\star}=\U_\sigma^n-\dfrac{\Delta t}{|K_\sigma|}\phi_{\sigma,\mathbf{x}}^{K,Rus}(\U^n)$$ for the Bernstein degrees of freedom,
then we have the following relation between the vectors, which components are the estimation of $\U$ at the Lagrange points:
$$\begin{pmatrix}
\U(\sigma_L)^{K,\star}\end{pmatrix}=\begin{pmatrix}\U^n(\sigma_L)\end{pmatrix}-\dfrac{\Delta t}{|K_\sigma|}M\begin{pmatrix}\phi_{\sigma,\mathbf{x}}^{K,Rus}(\U^n)\end{pmatrix},
$$
but since we approximate the flux as:
$$\F=\sum\limits_{\sigma\in K} \F_\sigma B_\sigma,$$ we see that 
$M\begin{pmatrix}\phi_{\sigma,\mathbf{x}}^{K,Rus}(\U^n)\end{pmatrix}$ is nothing more than the Rusanov residuals computed with the Lagrange interpolation
$$\F=\sum\limits_{\sigma_L\in K}\F(\sigma_L)\varphi_{\sigma_L}$$ where $\varphi_{\sigma_L}$ is the Lagrange polynomial for $\Sigma_L$. Hence we can use what we have done above and under the same condition, we have the positivity of the density and the internal energy at the Lagrange points.

\subsection{A less dissipative version of Rusanov}\label{Appendix:Subcells}
\begin{figure}[h]
\begin{center}
\subfigure[$P_2$]{\includegraphics[width=0.4\textwidth]{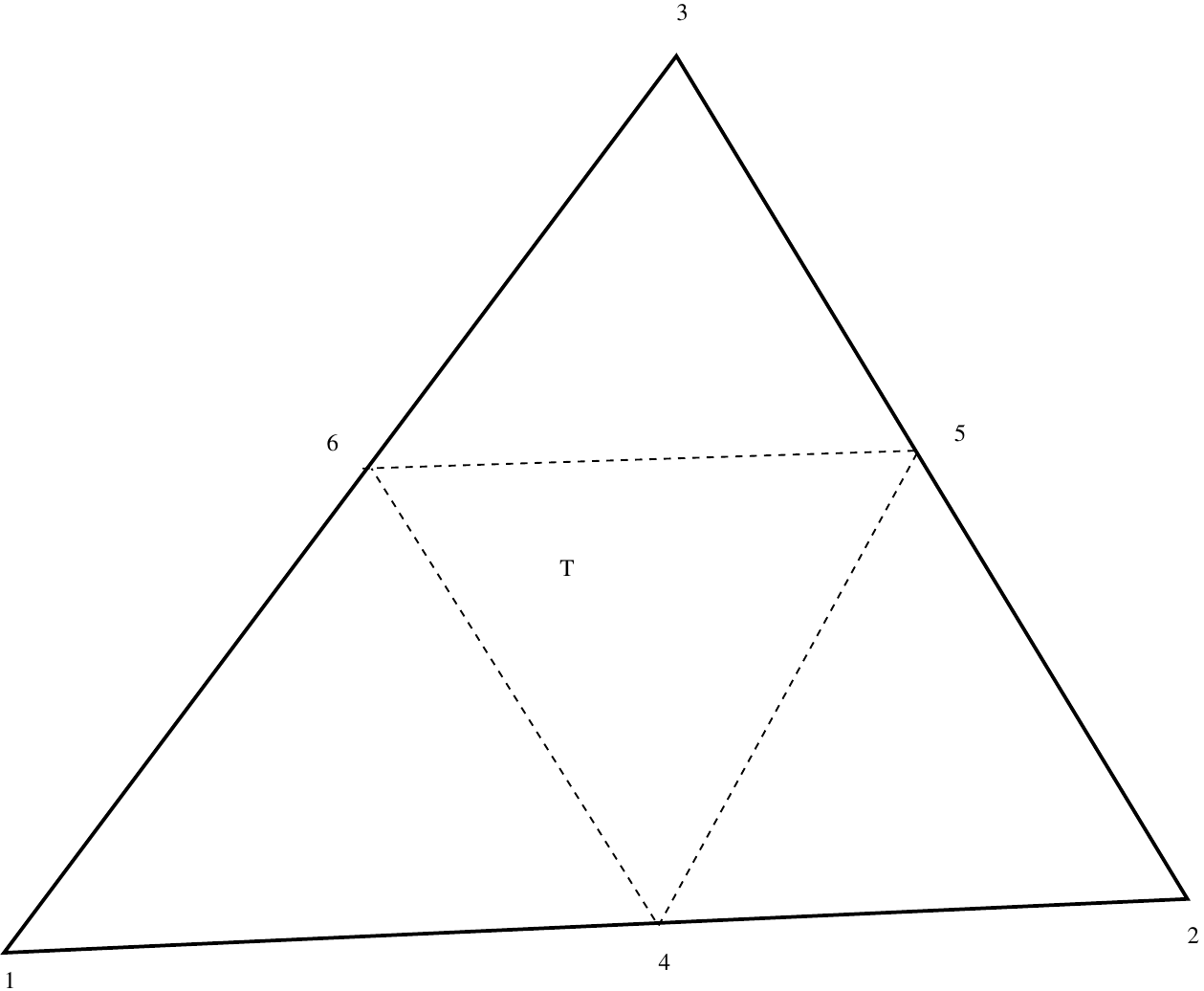}}
\subfigure[$P_3$]{\includegraphics[width=0.4\textwidth]{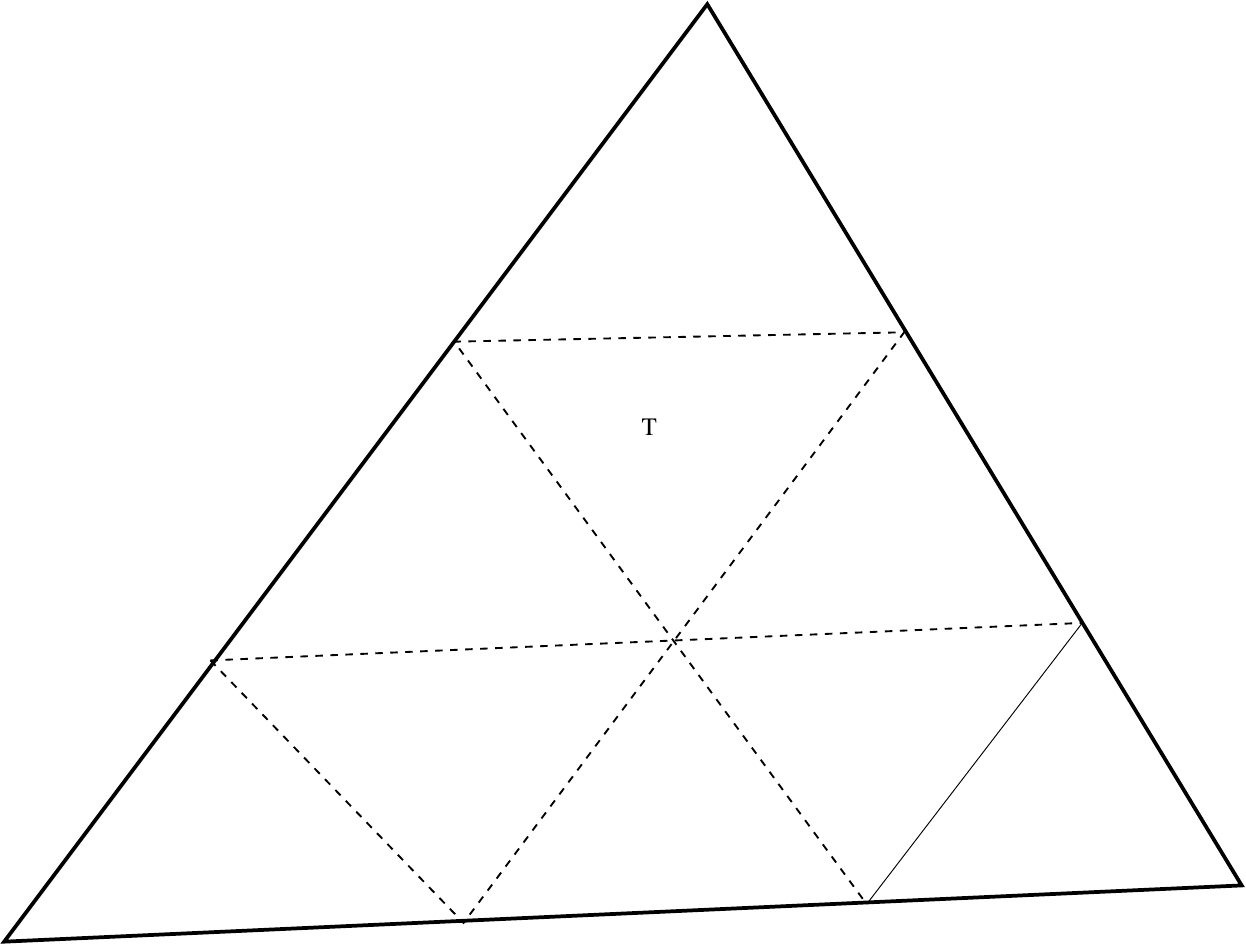}}
\end{center}
\caption{\label{sous triangles} Subdivision of a triangle into sub-triangles.}
\end{figure}

Consider $K$ a triangle (or a tetrahedron in 3D), and let us consider the Lagrange points on this element, and then the graph constructed from these points. See Figure \ref{sous triangles} for an illustration for $\mathbb{P}^2$ and $\mathbb{P}^3$ approximation.

In \cite{abg2001d} in the scalar case, and then in \cite{AbgrallViville2017} for the Euler equations (in two and three dimensions), we have noticed that if $\F$ is approximated by its Lagrange interpolant of degree $k$, $\F^{(k)}$, then
$$\int_K \text{ div }\F^{(k)}(\xx)\; d\xx=\sum_{K'\subset K}\omega_{K'}\int_{K'}\text{ div }\F^{(1)}(x)\; d\xx.
$$

For example, in the quadratic case,
$$\int_K \text{ div }\F^{(2)}(\xx)\; d\xx=\frac{1}{3} \int_{T_1}\text{ div }\F^{(1)}(x)\; d\xx+\frac{1}{3} \int_{T_2}\text{ div }\F^{(1)}(x)\; d\xx+\frac{1}{3} \int_{T_3}\text{ div }\F^{(1)}(x)\; d\xx+ \int_{T_4}\text{ div }\F^{(1)}(x)\; d\xx,$$
with
\begin{equation*}
\begin{split}
\int_{T_1}\text{ div }\F^{(1)}(x)\; d\xx&= \F_1\cdot \bn_1+\F_4\cdot \bn_2+\F_6\cdot \bn_3\\
\int_{T_2}\text{ div }\F^{(1)}(x)\; d\xx&= \F_4\cdot \bn_1+\F_2\cdot \bn_2+\F_5\cdot \bn_3\\
\int_{T_3}\text{ div }\F^{(1)}(x)\; d\xx&= \F_5\cdot \bn_1+\F_6\cdot \bn_2+\F_3\cdot \bn_3\\
\int_{T_4}\text{ div }\F^{(1)}(x)\; d\xx&= -\F_5\cdot \bn_1-\F_6\cdot \bn_2-\F_4\cdot \bn_3\\
\end{split}
\end{equation*}
If we use the Bernstein approximation, we have 
$$
\int_K \text{ div }\F^{(2)}(\xx)\; d\xx=\frac{2}{3} \int_{T_1}\text{ div }\F^{(1)}(x)\; d\xx+\frac{2}{3} \int_{T_2}\text{ div }\F^{(1)}(x)\; d\xx+\frac{2}{3} \int_{T_3}\text{ div }\F^{(1)}(x)\; d\xx,$$
\begin{equation*}
\begin{split}
\int_{T_1}\text{ div }\F^{(1)}(x)\; d\xx&= \F_1\cdot \bn_1+\F_4\cdot \bn_2+\F_6\cdot \bn_3\\
\int_{T_2}\text{ div }\F^{(1)}(x)\; d\xx&= \F_4\cdot \bn_1+\F_2\cdot \bn_2+\F_5\cdot \bn_3\\
\int_{T_3}\text{ div }\F^{(1)}(x)\; d\xx&= \F_5\cdot \bn_1+\F_6\cdot \bn_2+\F_3\cdot \bn_3\\
\end{split}
\end{equation*}
where here the $\F_i$s are the Bernstein DoFs (and not the Lagrange ones). We can easily check that the two formula are identical.
In each sub-triangle, one can define a Rusanov residual, for example in the Bernstein case, we have for $\sigma \in T_j$,
$$\phi_{\sigma}^{T_j}=\frac{1}{3}\int_{T_j}\text{ div }\F^{(1)}(x)\; d\xx+\alpha_j \big (\U_\sigma-\bar \U_{T_j}\big )$$
and then, introducing the residual
$$\Phi_\sigma=\sum_{T_j, \sigma\in T_j} \omega_{T_j}\Phi_\sigma^{T_j},$$
and using the same technique as above, we can show that the positivity of the density and the internal energy is guaranteed, provided that 
$$\alpha_{T_j}\geq \max\limits_{\xx\in T_j}\rho(A(\xx)).$$
The stencil of this residual is smaller than the Rusanov residual defined above.

We have not used these residual in the numerical experiments, as shown in \cite{AbgrallHO2018}, and if applied, this could lead to even less diffused results.


\end{document}